\documentclass[11pt]{article}
\usepackage{amssymb,amscd,amsxtra,calc}
\usepackage{cmmib57}
\usepackage{multirow}
\usepackage{mathrsfs}
\usepackage{graphicx}
\usepackage{bm}
\usepackage{amsmath}
\usepackage{amssymb}
\usepackage{ntheorem}
\usepackage{graphicx}
\usepackage{accents}
\usepackage{tikz} 
\usepackage{tikz-cd}
\usepackage{subfigure}
\usepackage{color}
\usepackage{stmaryrd}
\usepackage{comment}
\usepackage[font=footnotesize,labelfont=bf]{caption} 
\usepackage[section]{placeins} 
\allowdisplaybreaks[3]
\usepackage{import} 
\usepackage{colonequals}

\newtheorem{theorem}{Theorem}[section]
\newtheorem{proposition}[theorem]{Proposition}

\newtheorem{lemma}[theorem]{Lemma}

\newtheorem{remark}[theorem]{Remark}

\def\Proof{\noindent{\bf Proof} \quad}
\def\qed{\hfill $\Box$ \smallskip}

\def\Proof{\noindent{\bf Proof} \quad}
\def\qed{\hfill $\Box$ \smallskip}

\renewcommand{\d}{\mathrm{d}}
\newcommand{\p}{\partial}

\newcommand{\bn}{\mathbf{n}}
\newcommand{\bX}{\mathbf{X}}

\renewcommand{\l}{\left}
\renewcommand{\r}{\right}

\DeclareMathOperator{\tr}{tr}

\def\R{{\mathbb R}}

\def\P{{\mathbb P}}

\def\U{{\mathbb U}}

\def\K{{\mathbb K}}

\def\X{{\mathbb X}}

\def\cT{{\mathcal T}}

\newcommand{\bD}{\mathbb D}

\renewcommand\bX{{\boldsymbol{X}}}
\newcommand\bu{{\boldsymbol{u}}}
\newcommand\bx{{\boldsymbol{x}}}
\newcommand\be{{\boldsymbol{e}}}
\newcommand\bg{{\boldsymbol{g}}}
\newcommand\bp{{\boldsymbol{p}}}
\newcommand\bq{{\boldsymbol{q}}}

\newcommand\bM{{\boldsymbol{M}}}

\newcommand{\cE}{\mathcal{E}}
\newcommand{\cH}{\mathcal{H}}

\newcommand{\cV}{\mathcal{V}}
\newcommand{\cL}{\mathcal{L}}

\renewcommand\bn{{\boldsymbol{n}}}
\newcommand{\bsigma}{\boldsymbol{\sigma}}
\newcommand{\bnu}{\boldsymbol{\nu}}

\newcommand{\I}{\mathbb{I}}

\begin{document}
\title{Structure-preserving parametric finite element methods for two-phase Stokes flow based on Lagrange multiplier approaches}
\author{Harald Garcke\thanks{Fakult{\"a}t f\"ur Mathematik, Universit{\"a}t Regensburg, Regensburg, Germany ({harald.garcke@ur.de})}, \ Dennis Trautwein\thanks{Fakult{\"a}t f\"ur Mathematik, Universit{\"a}t Regensburg, Regensburg, Germany ({dennis.trautwein@ur.de})}, \ and \ Ganghui Zhang\thanks{Department of Mathematics, University of Edinburgh, EH9 3FD Edinburgh, United Kingdom ({zhangganghui2025@gmail.com})}}

\maketitle

\abstract{We present a novel formulation for parametric finite element methods to approximate two-phase Stokes flow. The new formulation is based on the classical Stokes equation in the bulk and a novel choice of interface conditions with additional Lagrange multipliers. This new Lagrange multiplier approach ensures that the numerical methods exactly preserve two physical structures of two-phase Stokes flow at the fully discrete level: (i) the energy-decaying and (ii) the volume-preserving properties. Moreover, different types of higher-order time discretization methods are employed, including the Crank--Nicolson method and the second-order backward differentiation formula approach. The resulting schemes are nonlinear and can be efficiently solved by using the Newton method with a decoupling technique. Extensive numerical experiments demonstrate that our methods achieve the desired temporal accuracy while preserving the two physical structures of the two-phase Stokes system.} 

{\small {\bf MSCcodes}:  65M60, 65M12, 35K55, 53C44

}

\section{Introduction}

In this paper, we consider the unsteady, viscous, incompressible two-phase Stokes system \cite{BGN_2013_stokes,HansboLZ14}:
\begin{subequations}\label{Stokes}
	     \begin{align}
	     -\mu_{\pm}\Delta\bu+\nabla p &=\bg,\qquad\quad\ \text{in }\Omega_{\pm}(t),\\
		\nabla\cdot \bu &=0,\qquad \quad\ \text{in }\Omega_{\pm}(t),
	     \end{align}	
	\end{subequations}
where $\Omega_{+}(t)$ and $\Omega_{-}(t)$ are two distinct, time-dependent subdomains of a bounded domain $\Omega\in\mathbb{R}^d$, $d\in\{2,3\}$, representing the two fluid phases. The vector field $\bu$ denotes the fluid velocity, $p$ the pressure, $\bg$ a possible external force, and $\mu_{\pm}>0$ are the viscosities of the respective fluids. 

Numerous methods 
have been considered for simulating  two-phase Stokes flow. There are two main types of approaches, depending on how  the interface is represented. One is to capture the interface implicitly by a global function defined on $\Omega$, for instance, the volume of fluid method \cite{Popinet2009,AlkeB09}, the level set method
\cite{Sussman1994, Gross2007, HansboLZ14} and the phase field method \cite{Abels2012,Feng2006,Jacqmin1999,Grun2014,GarckeHK16, LiLiuShenZ25,Aland2010,Lowengrub1998} are examples of such approaches. On the other hand, one can explicitly track the interface through front tracking methods \cite{Unverdi1992}, immersed interface methods \cite{Leveque1997} and the immersed boundary method \cite{Peskin2002}. This paper focuses on a numerical method representing the interface by a parametrization and using the finite element method. 

It is well-known that the two-phase Stokes model satisfies two physical structures \cite{BGN_2013_stokes}: an a~priori energy equality and a volume-conservation property; that is, during the time evolution, it holds that
\begin{subequations}\label{Physical-structures}
\begin{align}
    \frac{\d }{\d t}\cE(t)
    &=-2\int_{\Omega}\mu\, \bD(\bu):\bD(\bu)\ \d\cL^d+\int_{\Omega}\bu\cdot \bg\ \d \cL^d, 
    \\
    \frac{\d }{\d t}V(t)&=0.
\end{align}
\end{subequations}
Here $\cE(t)=\gamma\int_{\Gamma(t)}1\ \d\cH^{d-1}$ is the interfacial energy at time $t>0$, with $\Gamma(t)=\p\Omega_{-}(t) \cap \Omega$ being the interface, and $V(t)$ denotes the volume of $\Omega_-(t)$. The other quantities are defined below in Section~\ref{sec:stokes-model}.
In particular, the total energy is decreasing in the setting without external forces, i.e., $\bg\equiv 0$. Above, $\cL^d$ denotes the Lebesgue measure in $\mathbb{R}^d$ and $\cH^{d-1}$ is the $(d-1)$-dimensional Hausdorff measure.

While the construction of numerical approximations that preserve a discrete energy equality (or inequality) is well studied, achieving exact volume conservation remains challenging in practice, mainly due to surface tension forces.
In \cite{BGN_2013_stokes}, the authors proposed a fully linear finite element method for two-phase Stokes flow that satisfies a discrete energy inequality. Their method uses an unfitted bulk mesh, meaning that the bulk mesh and the interface mesh are completely independent. In practice, this approach shows good volume conservation properties, especially when the pressure space is enriched by the characteristic function of the enclosed subdomain $\Omega_-(t)$. 
This enrichment is known as the $\text{XFEM}_\Gamma$ approach, since the additional terms in the finite element system can be expressed as integrals over the interface $\Gamma(t)$, see \cite{BGN_2013_stokes, BGN_2015_navierstokes}. However, a limitation of this method is that exact volume conservation only holds at the time-continuous level, and not after time discretization. Despite this, the $\text{XFEM}_\Gamma$ method handles pressure jumps across the interface very well.
The ideas from \cite{BGN_2013_stokes} were later extended to fitted bulk meshes in \cite{agnese_nurnberg_2016_stokes}, where the bulk mesh adapts to the moving interface. This variant retains the main properties of the unfitted method, including the same drawback, i.e., the lack of exact volume conservation after time discretization. For other extended finite element methods ($\text{XFEM}$) which are also called cut finite element methods we refer to \cite{HansboLZ14, BurmanSHLM15, KirchhartGrossReusken16}.

To address the issue that on the time-discrete level the volume is not conserved, a nonlinear scheme was proposed in \cite{Garcke-Nurnberg-Zhao2023} for the two-phase Navier--Stokes model. Their approach achieves exact volume conservation even after time discretization, both for unfitted meshes and for the arbitrary Lagrangian-Eulerian (ALE) method with fitted meshes. The key idea builds on earlier work by \cite{Bao-Zhao, Nurnberg} and involves using time-weighted approximations of the interface normal vector $\bnu(t)$ in the finite element system. This leads to a nonlinear discrete formulation, but it allows the method to rigorously preserve both the volume and a discrete energy inequality at each time step. While this comes with higher computational cost due to the nonlinearity, the scheme represents a significant step forward in designing accurate and structure-preserving numerical methods for two-phase flows. Thus, it is desirable in the numerical community to design schemes that preserve both volume and energy structures at the fully discrete level. 

The construction of such structure-preserving schemes has been extensively studied in the context of parametric finite element methods applied to various other models, including surface diffusion \cite{Jiang21, Bao-Zhao, Bao-Li2024}, two-phase Navier--Stokes flow \cite{Garcke-Nurnberg-Zhao2023}, their viscoelastic extensions \cite{GNT_2024}, and the Mullins--Sekerka problem \cite{Nurnberg}.
A Lagrange multiplier approach was recently proposed to construct structure-preserving schemes for surface diffusion \cite{Garcke2024}. The key idea is to introduce two Lagrange multipliers to enforce two important physical structures: an energy equality and a volume conservation law for the region enclosed by the evolving surface. The main objective of this paper is to extend this Lagrange multiplier framework to the setting of two-phase Stokes flow, ensuring that both physical structures in \eqref{Physical-structures} are preserved while maintaining an accurate approximation of the interface.

In the literature, there are numerous results concerning Lagrange multiplier methods for the Stokes or Navier--Stokes equations. However, these approaches typically focus either on the single-phase flow setting, without any interface coupling \cite{DoanHJL25, YangTK22}, or on diffuse interface models such as the Cahn--Hilliard--Navier--Stokes system, see, e.g. \cite{LiLiuShenZ25, DoanHJL25b, cheng_shen_wang_2025_lagrange}. In contrast, our focus lies on sharp interface models for two-phase flow, where the interface is treated explicitly and its evolution is directly coupled to the surrounding fluid dynamics. The extension of a structure-preserving Lagrange multiplier framework to this setting presents new challenges and, to the best of our knowledge, has not been addressed in the existing literature.

While many studies have addressed the spatial discretization of one-phase or multi-phase flow models, often using front-tracking methods, the development of higher-order time discretizations for sharp interface models remains relatively limited. To our knowledge, the only works that consider higher-order schemes for two-phase Navier--Stokes flow are \cite{Bansch2012, Weller2017, Zhu2020, Frei2017}, where time discretizations based on backward differentiation formulas (BDF) are explored together with discontinuous Galerkin methods.
A further goal of this paper is to combine these types of higher-order time discretization methods with the Lagrange multiplier framework in the setting of two-phase Stokes flow, where our main focus lies on preserving the physical structures \eqref{Physical-structures} at the fully discrete level.

In summary, the two main novelties of our proposed method are:
\begin{itemize}
\item Compared to earlier front-tracking methods for two-phase flows, such as \cite{BGN_2013_stokes, BGN_2015_navierstokes, Garcke-Nurnberg-Zhao2023}, our fully discrete scheme exactly preserves both the energy equality (not only an inequality) and the volume identity at each time step.

\item Our formulation naturally allows the application of second-order time discretization methods. Instead of extrapolating the interface directly to the new time step, we use a lower-order predictor to ensure good mesh quality of the evolving interface.  
\end{itemize}   
We emphasize that, unlike previously studied higher-order time discretization schemes for two-phase flow models \cite{Bansch2012,Weller2017,Zhu2020,Frei2017}, our method is able to preserve both geometric structures at the fully discrete level. In addition to that, our numerical experiments confirm that the proposed method achieves the expected experimental order of convergence in several numerical tests.

The remainder of this work is organized as follows. In Section 2, we first review the classical model of two-phase Stokes flow and present our Lagrange multiplier approach. In addition, we introduce a semi-discrete approximation that is continuous in time and discrete in space. Section 3 is devoted to three different examples of structure-preserving approximations at the fully-discrete level. In Section 4, several numerical results are reported to demonstrate the accuracy and practicality of all three Lagrange multiplier schemes. Finally, in Section 5 we draw some conclusions.

\section{Model equations and semi-discrete approximation}

In this section, we first review the classical model for two-phase Stokes flow considered, e.g., in
\cite{BGN_2013_stokes,HansboLZ14,BaenschDGP23}. After that, we incorporate two Lagrange multipliers and introduce a semi-discrete approximation.

\subsection{The two-phase Stokes model} \label{sec:stokes-model}
\begin{figure}
\centering
\includegraphics[scale=0.6]{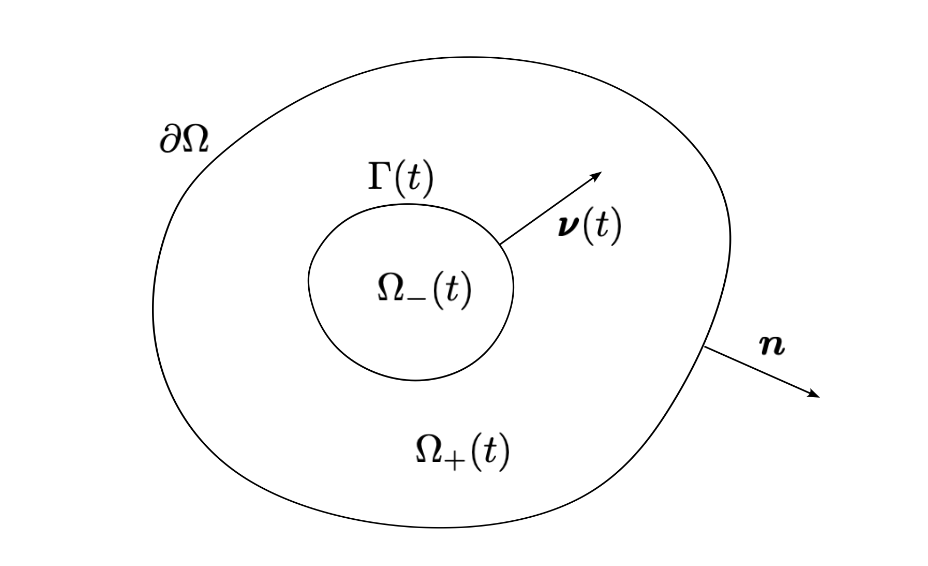}
\caption{A schematic sketch of the domain $\Omega$ in two space dimensions.}
\label{fig:domain}
\end{figure}

We let $\Omega\subset\R^d$, $d\in\{2,3\}$, be a bounded domain with boundary $\partial\Omega$ and let $t\in(0,T)$ with $T>0$. We assume that the domain $\Omega$ is decomposed into two distinct time-dependent regions $\Omega_+(t)$ and $\Omega_-(t)\colonequals \Omega\setminus \overline{\Omega_+(t)}$ for all $t\in[0,T]$. The two subdomains are separated by a time-evolving interface $\{\Gamma(t)\}_{t\in[0,T]}$ that is fully contained in $\Omega$ and does not intersect $\partial\Omega$. We further suppose that $\{\Gamma(t)\}_{t \in [0,T]}$ is an orientable and smooth $(d-1)$-dimensional hypersurface without boundary which is globally parametrized over a given reference manifold $\mathfrak{X}\subset \R^d$ by a mapping $\bX(\cdot,t) \colon\, \mathfrak{X} \to \R^d$, i.e., $\Gamma(t) = \bX(\mathfrak{X},t)$. 
We denote the normal velocity of the evolving hypersurface by $\cV=\p_t\bX\cdot \bnu$, where $\bnu$ is the normal vector on $\Gamma(t)$ pointing into $\Omega_{+}(t)$. 
We further denote the mean curvature on $\Gamma(t)$ by $\kappa(\cdot,t)$, with the sign convention $\kappa(\cdot,t)<0$ where $\Omega_{-}(t)$ is locally convex.
Moreover, $\bn$ is the outer unit normal on $\partial\Omega$. For a schematic sketch in two dimensions, we refer to Figure \ref{fig:domain}.

The two-phase Stokes model is formulated as follows. Find functions
\begin{alignat*}{2}
    &\bu\colon\, \Omega\times (0,T]\rightarrow \R^d,\qquad 
    &&p\colon\, \Omega\times (0,T]\rightarrow \R,\\
    &\bX\colon\, \mathfrak{X} \times (0,T]\rightarrow \R^d,\qquad 
    &&\kappa\colon\, \mathfrak{X}\times (0,T]\rightarrow \R,
\end{alignat*}
that satisfy the following equations for all $t\in(0,T]$.
\begin{itemize}
\item The Stokes equations in the bulk regions $\Omega_{\pm}(t)$:
\begin{subequations}\label{BGN:bulk}
\begin{alignat}{2}
    -\mu_\pm \nabla\cdot(\nabla\bu+(\nabla\bu)^\top )+\nabla p &=\bg,\qquad\quad &&\text{in }\Omega_{\pm}(t),
    \\
    \nabla\cdot \bu &=0,\qquad \quad\ &&\text{in }\Omega_{\pm}(t).
\end{alignat}	
\end{subequations}

\item The jump conditions at the interface $\Gamma(t)$:
\begin{subequations}\label{BGN:jump}
\begin{alignat}{2}
    \llbracket  \bu    \rrbracket  &=0,\qquad \qquad &&\text{on } \Gamma(t),\label{BGN:jump1}
    \\
    \llbracket \bsigma\bnu \rrbracket &= -\gamma \kappa\bnu,\qquad &&\text{on } \Gamma(t),\label{BGN:jump2}
\end{alignat}	
\end{subequations}
where $\gamma>0$ is a constant related to surface tension, 
\begin{align*}
    \bsigma=\mu(\nabla\bu+(\nabla\bu)^\top)-p\I=2\mu\bD(\bu)-p\I
\end{align*}
is the viscous stress tensor, $\I\in \R^{d\times d}$ denotes the identity matrix, and $\bD(\bu)=\frac{1}{2}(\nabla\bu+(\nabla\bu)^\top)$ is the rate of deformation tensor (or symmetrized velocity gradient). 
For future reference, we define the viscosity function $$\mu(\cdot,t) = \mu_+ (1-\bm{\chi}_{\Omega_-(t)}) + \mu_- \bm{\chi}_{\Omega_-(t)},$$ where $\bm{\chi}_{\Omega_-(t)}$ is the characteristic function of $\Omega_-(t)$.
Moreover, the terms $\llbracket \bu \rrbracket \colonequals \bu_+-\bu_-$ and $\llbracket \bsigma\bnu \rrbracket\colonequals\bsigma_+\bnu-\bsigma_-\bnu$ are defined as the jumps of the fluid velocity and normal part of stress tensor across the interface.

\item The interface flow equations:
\begin{subequations}\label{BGN:interface}
\begin{alignat}{2}
    \cV &= \bu \cdot \bnu ,\qquad \qquad &&\text{on } \Gamma(t),
    \label{Normal velocity} 
    \\
    \Delta_{\Gamma}\mathrm{Id}&=\kappa\bnu,\qquad \quad \quad &&\text{on } \Gamma(t), \label{Curvature identity} 
\end{alignat}
\end{subequations}
where $ \Delta_{\Gamma}$ is the surface Laplacian on $\Gamma(t)$.

\item The initial and boundary conditions:
\begin{subequations}\label{BGN:initial_and_boundary}
\begin{align}
    \Gamma(0)&=\Gamma^0,
    \\
    \bu &=0 \qquad \qquad \text{on }\p\Omega.
\end{align}
\end{subequations}
\end{itemize}

\smallskip
The curvature identity \eqref{Curvature identity} is  a well-established tool in the analysis of surface evolution problems \cite{Deckelnick-Dziuk-Elliott,Dziuk1991} and has also been widely used in the study of interface dynamics \cite{Bansch2001,Ganesan2007}. Additionally, the choice of the interface evolution law \eqref{BGN:interface} leads to favourable mesh properties in numerical simulations. Specifically, in two spatial dimensions, it results in an equidistribution of the interface nodes, while in three dimensions, the evolving hypersurface $\Gamma(t)$ forms a so-called conformal polyhedral surface \cite{BGN20}. Such properties are highly beneficial in practice, as they contribute to maintaining good mesh quality over time. For a detailed discussion and further applications, we refer to \cite{BGN20,BGN_2015_navierstokes,BGN2015b,BGN2015c}.

\smallskip

The two-phase Stokes model satisfies the well-known physical identities \eqref{Physical-structures}, for a proof we refer to \cite{BGN_2013_stokes, BaenschDGP23}. Therefore, it is desirable to design discrete schemes that preserve these two structures, as we will present in the next part.

\smallskip

\subsection{The Lagrange multiplier method} 
We now introduce our Lagrange multiplier approach. In addition to $(\bu,p,\bX,\kappa)$, we also seek two scalar Lagrange multipliers
\begin{align*}
	&\lambda_{\cE}:(0,T]\rightarrow \R,\qquad \lambda_{V}:(0,T]\rightarrow \R
\end{align*}
that satisfy the Stokes equations \eqref{BGN:bulk} in the bulk, the jump conditions \eqref{BGN:jump}, the initial and boundary conditions \eqref{BGN:initial_and_boundary}, and the interface equations \eqref{BGN:interface} with \eqref{Normal velocity} replaced by 
\begin{equation}\label{Lag:jump3}
    \cV = \bu \cdot \bnu +\lambda_{\cE}\kappa+\lambda_{V} 
    ,\qquad \qquad \text{on } \Gamma(t),
\end{equation}
where the Lagrange multipliers are used to enforce the energy and volume identities
\begin{subequations}\label{Lage:Enforce}
\begin{align}
    \frac{\d }{\d t}\cE(t)
    &=-2\int_{\Omega}\mu\, \bD(\bu):\bD(\bu)\ \d\cL^d+\int_{\Omega}\bu\cdot \bg\ \d \cL^d,\label{Lag:Enforce1} 
    \\
    \frac{\d }{\d t}V(t)&=0,\label{Lag:Enforce2}
\end{align}
\end{subequations}
where $ A:B=\tr ( A^T\, B)$ 
denotes the Hilbert--Schmidt inner product for matrices $A$ and $B$ in $\mathbb{R}^{d\times d}$.

\smallskip

For later use, we denote by $(\cdot, \cdot)$ the $L^2$-inner product over $\Omega$ and by $\l<\cdot, \cdot\r>_{\Gamma}$ the $L^2$-inner product over $\Gamma(t)$, respectively. 
When the context is clear, we may omit the explicit reference to the surface and simply write $\l<\cdot, \cdot\r>$ instead of $\l<\cdot, \cdot\r>_{\Gamma}$.

The two-phase Stokes model and the Lagrange multiplier method are closely related by the following theorem.

\begin{theorem}\label{Thm:relation:continuous}
Let $(\bu,p,\bX,\kappa)$ be a solution to the two-phase Stokes model \eqref{BGN:bulk}--\eqref{BGN:initial_and_boundary}. Then, $(\bu,p,\bX,\kappa,\lambda_{\cE},\lambda_{V})$ with $\lambda_{\cE}=\lambda_{V}=0$ is a solution to the Lagrange multiplier method. Conversely, if $(\bu,p,\bX,\kappa,\lambda_{\cE},\lambda_{V})$ is a solution to the Lagrange multiplier method and the curvature $\kappa(\cdot,t)$ is not constant with respect to the spatial variable, then it holds $\lambda_{\cE}=\lambda_{V}=0$ and hence $(\bu,p,\bX,\kappa)$ is a solution to the two-phase Stokes model \eqref{BGN:bulk}--\eqref{BGN:initial_and_boundary}.
\end{theorem}

\Proof
It suffices to show $\lambda_{\cE}=0$ and $\lambda_{V}=0$ whenever $(\bu,p,\bX,\kappa,\lambda_{\cE},\lambda_{V})$ is a solution to the Lagrange multiplier method at time $t>0$ and $\kappa(\cdot,t)$ is not constant. By using a version of Reynold's transport theorem \cite[Thm.~32]{BGN20} together with \eqref{Lag:jump3}, \eqref{BGN:bulk}, \eqref{BGN:jump} and \eqref{Normal velocity}, we compute 
\begin{align*}
	\frac{\d }{\d t}\cE(t)
	&=\gamma\frac{\d }{\d t}\l<1,1\r>=-\gamma \l<\kappa,\cV\r>\\
	&=-\gamma \l<\kappa,\bu\cdot \bnu\r>-\gamma\l<\kappa,\kappa\r>\lambda_{\cE}-\gamma \l<\kappa,1\r>\lambda_{V}\\
	&=-\l(2\mu\bD(\bu), \bD(\bu)\r)+\l(\bg,\bu\r) -\gamma\l<\kappa,\kappa\r>\lambda_{\cE}-\gamma \l<\kappa,1\r>\lambda_{V}.
\end{align*} 
Thus, combining with \eqref{Lag:Enforce1}, we obtain
\[
\gamma\l<\kappa,\kappa\r>\lambda_{\cE}+\gamma \l<\kappa,1\r>\lambda_{V}=0.
\]
Similarly, noting a surface  version of Reynold's transport theorem \cite[Thm.~33]{BGN20} and by combining it with the Stokes equations, \eqref{Lag:jump3}, \eqref{Normal velocity},
and \eqref{Lag:Enforce2}, we compute
\begin{align*}
	0=\frac{\d }{\d t}V(t)
	&=\frac{\d }{\d t}\int_{\Omega_-(t)}1\ \d \cL^{d}=\l<\cV,1 \r>\\
	&=\l<\bu\cdot \bnu,1\r>+\l<\kappa\lambda_{\cE}+\lambda_{V},1\r>\\
	&=\int_{\Omega_-(t)}\nabla\cdot \bu\, \d \cL^{d}+\l<\kappa,1\r>\lambda_{\cE}+\l<1,1\r>\lambda_{V}\\
	&=\l<\kappa,1\r>\lambda_{\cE}+\l<1,1\r>\lambda_{V}.
\end{align*}
Therefore, $\lambda_{\cE}$ and $\lambda_{V}$ satisfy a time-dependent algebraic system
\begin{equation*}
	\bM(t)\begin{bmatrix}
		\lambda_{\cE}\\
		\lambda_{V}
	\end{bmatrix}=0,\qquad \bM(t)=\begin{bmatrix}
		\gamma \l<\kappa,\kappa\r> & \gamma\l<\kappa,1\r>\\
		\l<\kappa,1\r> & \l<1,1\r>
	\end{bmatrix}.
\end{equation*}
Noticing that $\gamma$ is a positive constant, by the Cauchy-Schwarz inequality, we have
\[
	\det\l(\bM(t)\r)
	=\gamma\l[\l<\kappa,\kappa\r>\l<1,1\r>-
\l<\kappa,1\r>^2\r]\ge 0,
\]
 and the equality holds if and only if $\kappa(\cdot, t)$ is constant. Hence $\lambda_{\cE}=0$ and $\lambda_{V}=0$ if $\kappa(\cdot,t)$ is not constant and we completed the proof.
\qed

\smallskip

To introduce a weak formulation of our Lagrange multiplier method, we define the function spaces 
\begin{align*}
   & \U\colonequals H^1_0(\Omega,\R^d),\quad \P\colonequals L^2(\Omega),\quad \widehat{\P}\colonequals\l\{f\in \P:\int_{\Omega}f \ \d\cL^{d}=0 \r\},\\
    &\qquad\qquad \X\colonequals H^{1}(\mathfrak{X},\R^d),\qquad  \K\colonequals L^2(\mathfrak{X},\R^d).
\end{align*}
Then by integrating the equations for the Lagrange multiplier method over the respective spatial domains and by parts, a natural weak formulation for our Lagrange multiplier method is given as follows: Let the initial surface $\Gamma(0)$ be given. Then, for all $t>0$, we aim to find a solution 
\[
(\bu(\cdot,t), p(\cdot,t), \bX(\cdot,t), \kappa(\cdot,t)) \in \U\times \widehat{\P}\times \X\times\K
\]
and Lagrange multipliers $(\lambda_{\cE}(t),\lambda_{V}(t))\in \R\times \R$, such that 
\begin{subequations}\label{Lag:weak}
\begin{align}
    0&=\l(2\mu\bD(\bu),\bD(\bm{\varphi}) \r)-\l(p,\nabla\cdot \bm{\varphi}\r)-\gamma\l<\kappa\bnu,\bm{\varphi}\r>-\l(\bg,\bm{\varphi} \r),
    \\
    0&=\l(\nabla\cdot \bu,q \r),
    \\
    0&=\l<\p_{t}\bX\cdot \bnu,\eta\r>-\l<\bu\cdot\bnu, \eta\r>-\lambda_{\cE}\l<\kappa,\eta \r>-\lambda_{V}\l<1,\eta \r>,\label{Lag:weak3} 
    \\
    0&=\l<\kappa\bnu,\bm{\omega}\r>+\l<\nabla_{\Gamma}\bX,\nabla_{\Gamma}\bm{\omega}\r>,
    \\
    0&=\frac{\d }{\d t}\cE(t)+\l(2\mu\bD(\bu),\bD(\bu) \r)-\l(\bu, \bg\r), 
    \\
    0&=\frac{\d }{\d t}V(t),
\end{align}
\end{subequations}
for any $(\bm{\varphi},q,\eta,\bm{\omega})\in \U\times \widehat{\P}\times \K\times \X$,
where $\cE(t)$ and $V(t)$ are the energy and the enclosed volume, respectively, defined as in \eqref{Physical-structures}.

Similarly to \cite{BGN_2013_stokes, BGN_2015_navierstokes}, we note that we employ a slight abuse of notation in \eqref{Lag:weak}, in the sense that here and throughout this work, functions defined on the reference manifold $\mathfrak{X}$ are identified with functions defined on $\Gamma(t)$.

\subsection{Space discretization}\label{sec:discretization}

We now fix the discrete setting of our work.
We first recall the notation concerning finite elements from \cite{BGN_2013_stokes}. Assume that $\Omega^h \subset\mathbb{R}^d$, $d\in\{2,3\}$, is a bounded Lipschitz domain with polygonal (or polyhedral, respectively) boundary, and $\Gamma^h\subset \Omega^h$ is a $(d-1)$-dimensional polyhedral surface, fully contained in the interior of $\Omega$. 
We suppose there exist decompositions 
\[
\overline{\Omega^h }= \bigcup_{j=1}^{J_{\Omega^h}} \overline{o_j},
\qquad \Gamma^h = \bigcup_{j=1}^{J_{\Gamma^h}}\overline{\sigma}_j,
\]
where $\{o_j\}_{j=1}^{J_{\Omega^h}}$ and $\{\sigma_j\}_{j=1}^{J_{\Gamma^h}}$ are families of mutually disjoint open $d$-simplices and $(d-1)$-simplices, respectively. We denote the vertices of $\Gamma^h$ by $\{\bq_k\}_{k=1}^{K_{\Gamma^h}}$. 
We assume throughout that $\Gamma^h$ is non-degenerate in the following sense: firstly, $\cH^{d-1}(\sigma_j)>0$ for all $j=1,\ldots, J_{\Gamma^h}$, and secondly, $\Gamma^h$ has no self-intersections.
We denote by $\cT^h$ the regular partitioning of $\Omega^h$, where 
\[h\colonequals\max_{j=1,\ldots,J_{\Omega^h}}\mathrm{diam}(o_j).
\] 
We define the finite element spaces on $\Omega^h$ by
\begin{align*}
    S_k&\colonequals\{\chi\in C(\overline{\Omega}^h): \chi|_{o}\in P_{k}(o),\, \forall\, o\in \cT^h  \}\subset W^{1,\infty}(\Omega^h),\quad k\in \mathbb{N}_{\ge 1},
    \\
    S_0&\colonequals\{\chi\in L^\infty(\Omega^h): \chi|_{o}\in P_{0}(o),\, \forall\, o\in \cT^h  \}\subset L^\infty(\Omega^h),
\end{align*}
where $P_{k}(o)$ denotes the space of polynomials of degree less than or equal to $k\in\mathbb{N}_{\geq0}$ on $o\in \cT^h$.
We choose the discrete function space tuple $(\U^h,\widehat{\P}^h)\subset (\U,\widehat{\P})$ for velocity/pressure such that the usual inf-sup condition is satisfied, meaning there exists a constant $C_0>0$ independent of $h$ such that
\begin{equation}
    \inf_{\varphi\in \widehat{\P}^h}\sup_{\bm{\xi}\in \U^h}\frac{(\varphi,\nabla\cdot \bm{\xi})}{\|\varphi\|_{L^2}\|\bm{\xi}\|_{H^1}}\ge C_0>0.
\end{equation}
For instance, we can choose the lowest order Taylor-Hood element $\U^h=[S_2]^d\cap \U$, $\widehat{\P}^h=S_1\cap \widehat{\P}$ for $d=2,3$, also called the P2-P1 element. Other inf-sup stable examples are the P2-P0 element $\U^h=[S_2]^d\cap \U$, $\widehat{\P}^h= S_0\cap \widehat{\P}$ for $d=2$, or the P2-(P1+P0) element $\U^h=[S_2]^d\cap \U$, $\widehat{\P}^h= (S_1 + S_0) \cap \widehat{\P}$ for $d=2$. We note that proofs can be found in, e.g., \cite{boffi_2012_infsup, ern_guermond_2004}. 
The inf-sup inequality is fulfilled for these examples when the family of meshes $\{\cT^h \}_{h>0}$ is shape regular. Sometimes, an additional (but very mild) assumption on the triangulation is posed due to technical reasons, but usually it can be dropped if the triangulation $\cT^h$ is fine enough, see \cite[Chap.~8.8]{BBF2013_fem}. We note that the results of this work are not restricted to the above stable pairs. However, other inf-sup stable examples such as the (P1+bubble)-P1 element (MINI element) and the (P1-iso-P2)-P1 element \cite{ern_guermond_2004} are not considered here.

Similar to \cite{BGN_2013_stokes}, we define the finite element space on $\Gamma^h$ by
\[
W^h\colonequals\{\chi\in C(\Gamma^h): \chi|_{\sigma_j}\in P_1(\sigma_j),\  j=1,\ldots,J_{\Gamma^h}\}\subset H^1(\Gamma^h),
\]
where $P_1(\sigma_j)$ denotes the space of polynomials of degree $1$ on $\sigma_j$, $j=1,\ldots,J_{\Gamma^h}$. Following \cite{BGN_2013_stokes}, we define the $L^2$-inner product and its lumped mass version on polyhedral surface $\Gamma^h$ as follows 
\[
\l<u,v \r>_{\Gamma^h}\colonequals\int_{\Gamma^h}u\cdot v\ \d\cH^{d-1},\quad \l<u,v \r>_{\Gamma^h}^h\colonequals\frac{1}{d}\sum_{j=1}^{J_{\Gamma^h}}\cH^{d-1}(\sigma_j)\sum_{k=1}^d (u\cdot v)(\bq_{j_k}^{-}),
\] 
where $\{\bq_{j_k}\}_{k=1}^d$ are the vertices of $\sigma_j$, and we define $u(\bq_{j_k}^{-})\colonequals\lim_{\sigma_j\ni\bp\rightarrow\bq_{j_k}}u(\bp)$.

\begin{figure}
\centering
\includegraphics[width=0.35\linewidth,trim={2cm 3.8cm 2cm 3.8cm},clip]{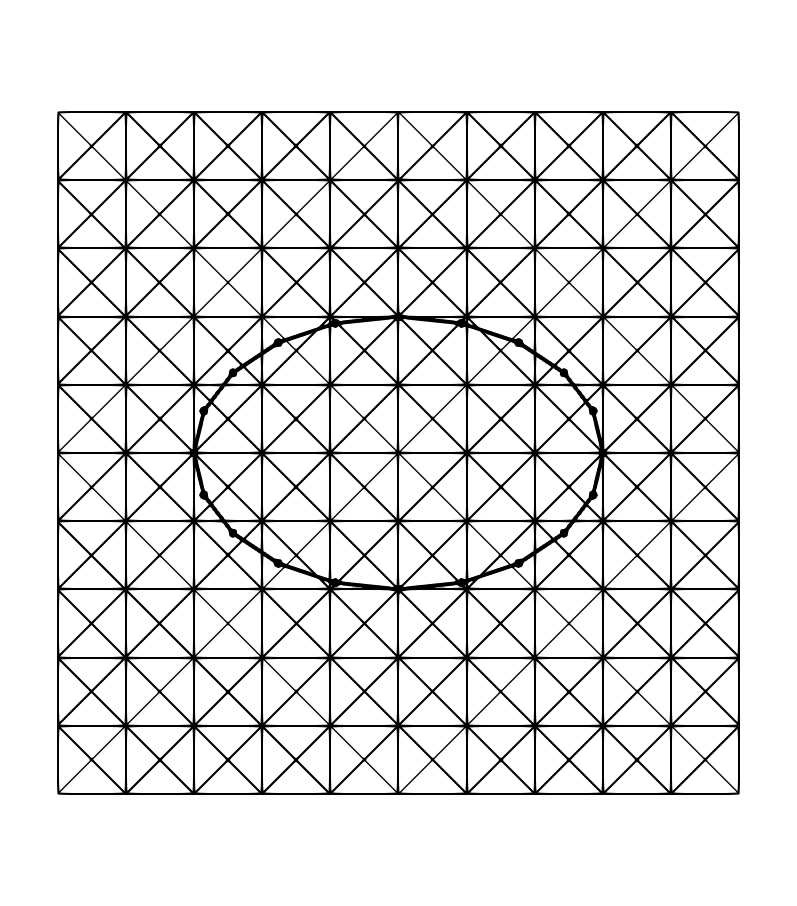}
\includegraphics[width=0.35\linewidth,trim={2cm 3.8cm 2cm 3.8cm},clip]{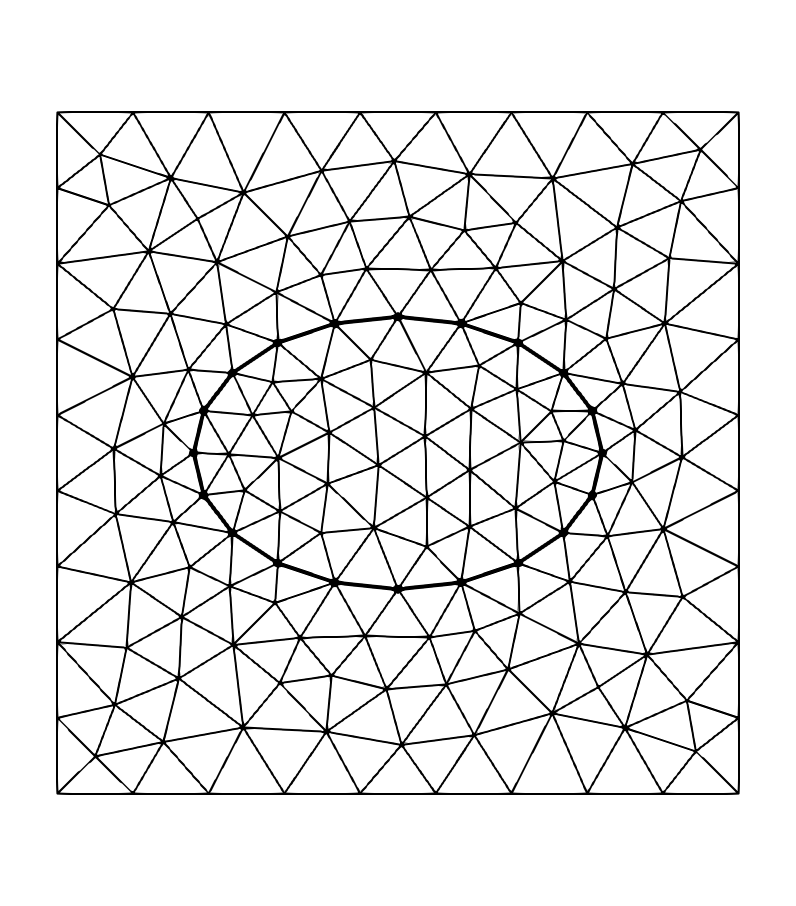}
\caption{Unfitted and fitted bulk meshes together with an interface mesh.}
\label{fig:meshes}
\end{figure}

In general, there are two common strategies for handling the interface and the bulk meshes in numerical simulations. In the unfitted mesh approach \cite{BGN_2013_stokes}, the interface and the bulk meshes are completely independent of each other. This allows for greater flexibility, since one can usually avoid remeshing or moving of the bulk mesh at every time step, but can complicate the evaluation of certain coupling terms between interface and bulk variables of the form $\l<\bu \cdot \bnu, \eta\r>$, see \eqref{Lag:weak3}.
In contrast, the fitted mesh approach ensures that the bulk mesh is aligned with the interface, for example by remeshing or adapting it to the evolving interface mesh at every time step, see \cite{agnese_nurnberg_2016_stokes}. This simplifies the numerical integration of terms involving both bulk and interface quantities. 
In this work, we allow both the fitted and the unfitted approaches. We refer to Figure \ref{fig:meshes}, where typical meshes for both approaches are shown.

Given a polyhedral surface $\Gamma^h$ and a decomposition $\cT^h$ of $\Omega^h$, the unfitted approach makes the bulk mesh being decomposed into three disjoint parts: the elements in the interior $\Omega^h_-$ of $\Gamma^h$, the elements in the exterior $\Omega^h_+$ of $\Gamma^h$, and the elements with non-empty intersection with the interface $\Gamma^h$, i.e.,
\begin{align*}
    \cT^h_-&\colonequals\{o\in \cT^h:o\subset \Omega^h_-\},\\
    \cT^h_+&\colonequals\{o\in \cT^h:o\subset \Omega^h_+\},\\
    \cT^h_{\Gamma^h}&\colonequals\{o\in \cT^h:o\cap \Gamma^h \neq \emptyset\}.
\end{align*}
Note that the set $\cT^h_{\Gamma^h}$ is empty for the fitted approach, i.e., when the bulk mesh is fitted to the interface mesh.
As previously stated, we always assume that $\Gamma^h$ has no self-intersections and we define the piecewise constant unit normal $\bnu^h$ to $\Gamma^h$ such that $\bnu^h$ points into $\Omega^h_+$. 
Similarly to \cite{BGN_2013_stokes}, we define the discrete viscosity function $\mu^h\in S_0$ as
\begin{equation}\label{Discrete_mu}
    \mu^h|_{o}=\begin{cases}
        \mu_- & o\in \cT^h_-,\\
        \mu_+ & o\in \cT^h_+,\\
        \frac{1}{2}(\mu_-+\mu_+) & o\in \cT^h_{\Gamma^h}.
    \end{cases}
\end{equation}
For other possible definitions we refer to, e.g., \cite{BGN_2015_navierstokes}.

\subsection{A semi-discrete scheme}

We now consider the following semi-discrete approximation of \eqref{Lag:weak} based on the above notations. Here, we allow for time-dependent discrete function spaces, since in general the bulk and interface meshes can change in time.

Let $\Gamma^h(0)$ be a given initial polyhedral surface. Then, for all $t>0$, we aim to find a solution 
\[
\l(\bu^h(\cdot,t), p^h(\cdot,t), \bX^h(\cdot,t), \kappa^h(\cdot,t) \r) \, \in \, \U^h(t)\times \widehat{\P}^h(t) \times [W^h(t)]^d\times W^h(t)
\]
and Lagrange multipliers $(\lambda_{\cE}^h(t), \lambda_{V}^h(t))\in \R\times \R$, such that
\begin{subequations}\label{Lag:Semi}
\begin{align}
    0 &= \l(2\mu^h\bD(\bu^h),\bD(\bm{\varphi}^h) \r)-\l(p^h,\nabla\cdot \bm{\varphi}^h\r)-\gamma\l<\kappa^h\bnu^h,\bm{\varphi}^h\r>_{\Gamma^h} -\l(\bg^h,\bm{\varphi}^h \r),
    \\
    0 &= \l(\nabla\cdot \bu^h,q^h \r), \label{Lag:Semi2} 
    \\
    0 &= \l<\p_{t}\bX^h\cdot \bnu^h,\eta^h\r>_{\Gamma^h}^h-\l<\bu^h\cdot\bnu^h, \eta^h\r>_{\Gamma^h}-\lambda_{\cE}^{h}\l<\kappa^{h},\eta^h \r>_{\Gamma^h}^h-\lambda_{V}^{h}\l<1,\eta^{h} \r>_{\Gamma^h}^h, \\
    0 &= \l<\kappa^h\bnu^h,\bm{\omega}^h\r>_{\Gamma^h}^h+\l<\nabla_{\Gamma^h}\bX^h,\nabla_{\Gamma^h}\bm{\omega}^h\r>_{\Gamma^h},
    \\
    0 &= \frac{\d }{\d t}\cE^h(t)+\l(2\mu^h\bD(\bu^h),\bD(\bu^h) \r)-\l(\bu^{h}, \bg^h\r),\label{Lag:Semi5} 
    \\
    0 &= \frac{\d }{\d t}V^h(t),\label{Lag:Semi6}
\end{align}
\end{subequations}
for any $(\bm{\varphi}^h,q^h,\eta^h,\bm{\omega}^h)\in \U^h(t) \times \widehat{\P}^h(t) \times W^h(t) \times [W^h(t)]^d$,
where $V^h(t)$ is the volume of $\Omega^h_{-}(t)$, and $\cE^h(t)$ is defined as 
\[
\cE^h(t)=\gamma\int_{\Gamma^h(t)}1\ \d\cH^{d-1}.
\]
Moreover, we identify the polyhedral surface with the parametrization $\Gamma^h(t)=X^h(\mathfrak{X},t)$, and define $\bg^h\colonequals I^h_2\bg$, where $I^h_2$ is the standard interpolation operator onto $[S_2]^d$.

\smallskip

To relate our semi-discrete scheme \eqref{Lag:Semi} with the semi-discrete scheme from \cite{BGN_2013_stokes} for the two-phase Stokes model, we recall their $\text{XFEM}_\Gamma$ approach. Their main feature is to include the characteristic function of $\Omega^h_-(t)$ as a test function in \eqref{Lag:Semi2}. We observe that $\bm{\chi}_{\Omega^h_-(t)} \in \mathbb{P}^h(t)$ is automatically satisfied in the case of fitted bulk meshes, provided that $\mathcal{S}_0(t) \subset \mathbb{P}^h(t)$, where $\P^h(t) := \widehat{\mathbb{P}}^h(t) \oplus \mathrm{span}\{1\}$. 
In contrast, for unfitted bulk meshes, this cannot be guaranteed a priori, and thus the condition $\bm{\chi}_{\Omega^h_-(t)} \in \mathbb{P}^h(t)$ needs to be imposed explicitly.
Indeed, this $\text{XFEM}_\Gamma$ approach has been proposed to ensure that the volume $V^h(t)$ is preserved at the semi-discrete level \cite{BGN_2013_stokes}. The following theorem relates the semi-discrete scheme with $\text{XFEM}_\Gamma$ from \cite{BGN_2013_stokes} with our semi-discrete scheme \eqref{Lag:Semi} for the Lagrange multiplier method.

\begin{theorem}\label{Thm:relation:semi}
Let $(\bu^h,p^h,\bX^h,\kappa^h,\lambda_{\cE}^h,\lambda_{V}^h)$ be a solution to the semi-discrete Lagrange multiplier method \eqref{Lag:Semi} at time $t>0$ and assume that $\bm{\chi}_{\Omega^h_-(t)} \in \mathbb{P}^h(t)$. If the discrete curvature $\kappa^h(\cdot,t)$ is not constant (with respect to the spatial variable), then it holds $\lambda_{\cE}^h(t)=\lambda_{V}^h(t)=0$.
\end{theorem}

\Proof
Taking in \eqref{Lag:Semi} $\bm{\varphi}^h=\bu^h$, $q^h=p^h$, $\eta^h=\kappa^h$ 
$\bm{\omega}^h=\p_t\bX^h$, combining with \eqref{Lag:Semi5} and using the identity (cf.~\cite[Theorem~70 and Lemma~9]{BGN20}) 

\begin{equation}
    \frac{\d }{\d t}\cE^h(t)=\gamma\l<\nabla_{\Gamma^h}\bX^h,\nabla_{\Gamma^h}\p_t\bX^h\r>_{\Gamma^h},
\end{equation}
we obtain 
\begin{align*}
    &\lambda_{\cE}^{h}\gamma\l<\kappa^{h},\kappa^h \r>_{\Gamma^h}+\lambda_{V}^{h}\gamma\l<1,\kappa^h \r>_{\Gamma^h}^h
    \\
    &=\gamma\l<\p_{t}\bX^h\cdot \bnu^h,\kappa^{h}\r>_{\Gamma^h}^h-\gamma\l<\bu^h\cdot\bnu^h, \kappa^{h}\r>_{\Gamma^h}
    \\
    &=-\gamma\l<\nabla_{\Gamma^h}\bX^h,\nabla_{\Gamma^h}\p_{t}\bX^h\r>_{\Gamma^h}-\l[\l(2\mu^h\bD(\bu^h),\bD(\bu^h) \r)-\l(p^h,\nabla\cdot \bu^h\r)-(\bg^h,\bu^h)\r]
    \\
    &= -\frac{\d }{\d t}\cE^h(t)-\l[\l(2\mu^h\bD(\bu^h),\bD(\bu^h) \r)-(\bg^h,\bu^h)\r]=0.
\end{align*}
On the other hand, choosing in \eqref{Lag:Semi} $q^h=\bm{\chi}_{\Omega^h_-}-\frac{\cL^d(\Omega^h_-)}{\cL^d(\Omega^h)}$, $\eta^h=1$, combining with \eqref{Lag:Semi6}, and noting $\bm{\chi}_{\Omega^h_-(t)} \in \mathbb{P}^h(t)$ and Reynold's transport theorem (cf.~\cite[Theorem~71]{BGN20}), so that
\begin{equation}
    \frac{\d}{\d t}V^h(t)=\l<\p_t \bX^h\cdot \bnu^h,1 \r>_{\Gamma^h},
\end{equation}
we get 
\begin{align*}
    0
    &=\frac{\d }{\d t}V^h(t)
    =\l<\p_{t}\bX^h\cdot \bnu^h,1\r>_{\Gamma^h}=\l<\p_{t}\bX^h\cdot \bnu^h,1\r>_{\Gamma^h}^h\\
    &=\l<\bu^h\cdot\bnu^h, 1\r>_{\Gamma^h}+\lambda_{\cE}^{h}\l<\kappa^{h},1 \r>_{\Gamma^h}^h+\lambda_{V}^{h}\l<1,1 \r>_{\Gamma^h}^h\\
    &=\int_{\Omega_-^h}\nabla\cdot \bu^h\ \d \cL^d+\lambda_{\cE}^{h}\l<\kappa^{h},1 \r>_{\Gamma^h}^h+\lambda_{V}^{h}\l<1,1 \r>_{\Gamma^h}^h\\
    &=\l(\nabla\cdot \bu^h, \bm{\chi}_{\Omega^h_-}-\frac{\cL^d(\Omega^h_-)}{\cL^d(\Omega^h)}\r)+\lambda_{\cE}^{h}\l<\kappa^{h},1 \r>_{\Gamma^h}^h+\lambda_{V}^{h}\l<1,1 \r>_{\Gamma^h}^h\\
    &=\lambda_{\cE}^{h}\l<\kappa^{h},1 \r>_{\Gamma^h}^h+\lambda_{V}^{h}\l<1,1 \r>_{\Gamma^h}^h.
\end{align*}
Therefore, $\lambda^h_{\cE}$ and $\lambda^h_{V}$ satisfy the time-dependent algebraic system
\begin{equation*}
	\bM^h(t)\begin{bmatrix}
		\lambda^h_{\cE}\\
		\lambda^h_{V}
	\end{bmatrix}=0,\qquad \bM^h(t)=\begin{bmatrix}
		\gamma \l<\kappa^h,\kappa^h\r>_{\Gamma^h}^h & \gamma\l<\kappa^h,1\r>_{\Gamma^h}^h\\
		\l<\kappa^h,1\r>_{\Gamma^h}^h & \l<1,1\r>_{\Gamma^h}^h
	\end{bmatrix}.
\end{equation*}
An argument based on the Cauchy--Schwarz inequality, similar to the one used in Theorem~\ref{Thm:relation:continuous}, completes the proof.
\qed

\begin{remark} 
Theorems~\ref{Thm:relation:continuous} and~\ref{Thm:relation:semi} show that our Lagrange multiplier formulation is equivalent to the classical two-phase Stokes model and its semi-discrete approximation from \cite{BGN_2013_stokes}, respectively, before the equilibrium is reached, and therefore inherits its good mesh quality. However, our approach offers significant advantages at the fully discrete level, as it ensures both energy dissipation and exact volume conservation without requiring an XFEM method as in \cite{BGN_2013_stokes} in the case of unfitted bulk meshes. In addition, several standard time discretization schemes can be naturally applied within our Lagrange multiplier framework.
\end{remark}

\section{Time discretization}
Building on the semi-discrete formulation \eqref{Lag:Semi}, we now consider three time discretization strategies: Euler, Crank--Nicolson, and the second-order backward differentiation formula (BDF2). Each approach leads to a nonlinear fully discrete scheme, which we solve using Newton iterations. For this, we propose a simple and efficient strategy that exploits the structure of the equations.

In the following, we assume that the time interval $[0,T]$ is partitioned into equidistant subintervals $[t_{m-1}, t_m)$ with $t_m=m \Delta t$, $\Delta t = \frac{T}{M}$, $M\in\mathbb{N}$ and $m \in \{0, \ldots, M\}$.
We use the super-index $m$ for time discrete quantities at time $t_m = m\Delta t$, including finite element spaces such as $\U^m$ and $\widehat{\P}^m$, functions and meshes.
Moreover, we assume throughout that the appearing polyhedral surfaces $\Gamma^{m-1}$, $\Gamma^m$, $\widetilde{\Gamma}^{m+\frac12}$ and $\widetilde{\Gamma}^{m+1}$, respectively, satisfy the non-degeneracy property stated in Section~\ref{sec:discretization}, i.e., they have no self-intersections and all interface elements have positive Hausdorff measure.

\subsection{Fully discrete numerical schemes}

\subsubsection{Euler scheme}
We first consider an Euler type discretization of \eqref{Lag:Semi}.

Let $m\in\{0,\ldots,M-1\}$. Given $\Gamma^m$, we aim to find a solution
\[
(\bu^{m+1},p^{m+1},\bX^{m+1},\kappa^{m+1})\in \U^m\times \widehat{\P}^m\times [W^m]^d\times W^m
\]
and Lagrange multipliers $(\lambda_{\cE}^{m+1},\lambda_{V}^{m+1})\in \R\times \R$, such that 
\begin{subequations}\label{Lag:Euler}
\begin{align}
    0&=\l(2\mu^{m}\bD(\bu^{m+1}),\bD(\bm{\varphi}^h) \r) -\l(p^{m+1},\nabla\cdot \bm{\varphi}^h\r)-\gamma\l<\kappa^{m+1}\bnu^m,\bm{\varphi}^h\r>_{\Gamma^m}-\l(\bg^{m} ,\bm{\varphi}^h \r),
    \\
    0&=\l(\nabla\cdot \bu^{m+1},q^h \r),
    \\
     0&=\l<\l(\frac{\bX^{m+1}-\bX^m}{\Delta t} \r)\cdot \bnu^m,\eta^h\r>_{\Gamma^m}^h -\l<\bu^{m+1}\cdot \bnu^m,\eta^h\r>_{\Gamma^m} \notag 
     \\
     &\qquad\qquad -\lambda_{\cE}^{m+1}\l<\kappa^{m+1},\eta^h \r>_{\Gamma^m}^h-\lambda_{V}^{m+1}\l<1,\eta^{h} \r>_{\Gamma^m}^h, \label{Lag:Euler3}
     \\
    0&=\l<\kappa^{m+1}\bnu^m,\bm{\omega}^h\r>_{\Gamma^m}^h+\l<\nabla_{\Gamma^m}\bX^{m+1},\nabla_{\Gamma^m}\bm{\omega}^h\r>_{\Gamma^m},\qquad\quad\qquad\quad\quad\quad
    \\
    0&=\frac{\cE^{m+1}-\cE^{m}}{\Delta t}+\l(2\mu^m\bD(\bu^{m+1}),\bD(\bu^{m+1}) \r)-\l(\bu^{m+1}, \bg^{m}\r),\label{Lag:Euler5} 
    \\
    0&=V^{m+1}-V^m,	\label{Lag:Euler6}
\end{align}
\end{subequations}
for any $(\bm{\varphi}^h,q^h,\eta^h,\bm{\omega}^h)\in \U^m\times \widehat{\P}^m\times W^m\times [W^m]^d$.
Here, $V^m$ is the volume of $\Omega_-^m$ and $\cE^m$ is defined as 
\[
\cE^{m}=\cE(\bX^{m})=\gamma\int_{\Gamma^m}1\ \d \cH^{d-1}.
\]
We note that the scheme \eqref{Lag:Euler} corresponds to an Euler-type discretization in time, where, similar to \cite{BGN_2013_stokes}, the surface integrals are evaluated over the previous interface $\Gamma^m$. Furthermore, the discrete viscosity $\mu^m$ and the discrete normal vector $\bnu^m$ are determined by $\Gamma^m$, and the discrete body force is given by $\bg^{m} = I^h_2\bg(t_{m})$. The nonlinearity of the system arises in the equations \eqref{Lag:Euler3}, \eqref{Lag:Euler5}, and \eqref{Lag:Euler6}. Once the system has been solved, the updated polyhedral surface is obtained as a parametrization over the previous surface, i.e., $\Gamma^{m+1} = \bX^{m+1}(\Gamma^m)$.

\begin{remark}
For later comparison, we recall the fully linear scheme for two-phase Stokes flow proposed by Barrett, Garcke and N\"urnberg \cite{BGN_2013_stokes}, also abbreviated as the BGN scheme. Let $m\in\{0,\ldots,M-1\}$. Given $\Gamma^m$, we aim to find a solution
\[
(\bu^{m+1},p^{m+1},\bX^{m+1},\kappa^{m+1})\in \U^m\times \widehat{\P}^m\times [W^m]^d\times W^m
\]
such that 
\begin{subequations}\label{BGN}
\begin{align}
    0&=\l(2\mu^{m}\bD(\bu^{m+1}),\bD(\bm{\varphi}^h) \r) -\l(p^{m+1},\nabla\cdot \bm{\varphi}^h\r)-\gamma\l<\kappa^{m+1}\bnu^m,\bm{\varphi}^h\r>_{\Gamma^m}-\l(\bg^{m},\bm{\varphi}^h \r),
    \\
    0&=\l(\nabla\cdot \bu^{m+1},q^h \r),
    \\
     0&=\l<\l(\frac{\bX^{m+1}-\bX^m}{\Delta t} \r)\cdot \bnu^m,\eta^h\r>_{\Gamma^m}^h -\l<\bu^{m+1}\cdot \bnu^m,\eta^h\r>_{\Gamma^m},\\
    0&=\l<\kappa^{m+1}\bnu^m,\bm{\omega}^h\r>_{\Gamma^m}^h+\l<\nabla_{\Gamma^m}\bX^{m+1},\nabla_{\Gamma^m}\bm{\omega}^h\r>_{\Gamma^m},
    \qquad\quad\qquad\quad\quad\quad
\end{align}
\end{subequations}
for any $(\bm{\varphi}^h,q^h,\eta^h,\bm{\omega}^h)\in \U^m\times \widehat{\P}^m\times W^m\times [W^m]^d$, and then set $\Gamma^{m+1} = \bX^{m+1}(\Gamma^m)$.

For the scheme \eqref{BGN}, the unconditional energy stability was shown in \cite{BGN_2013_stokes}, by proving a discrete energy inequality corresponding to \eqref{Lag:Euler5}. Unique solvability under a mild assumption was also established. However, the temporal accuracy of the scheme is limited to the first order, and it does not preserve the volume at the fully discrete level.
\end{remark}

\smallskip

\subsubsection{Crank--Nicolson scheme}
Subsequently, we consider the Crank--Nicolson discretization based on a predictor-corrector method.

Let $m\in\{0,\ldots,M-1\}$. Let $\Gamma^m$ and an approximation $\widetilde{\Gamma}^{m+\frac12}$ of $\Gamma(t_{m+\frac12})$ be given, where $t_{m+\frac12} = (m+\frac12)\Delta t$. Then, we aim to find a solution
\[
(\bu^{m+1},p^{m+1}, \bX^{m+1},\kappa^{m+1})\in \U^m\times \widehat{\P}^m\times [\widetilde{W}^{m+\frac12}]^d\times \widetilde{W}^{m+\frac12}
\]
and Lagrange multpliers $(\lambda_{\cE}^{m+1},\lambda_{V}^{m+1})\in \R\times \R$, such that 
\begin{subequations}\label{Lag:CN}
\begin{align}
    0&=\l(2\widetilde{\mu}^{m+\frac{1}{2}}\bD(\bu^{m+\frac{1}{2}}),\bD(\bm{\varphi}^h) \r)-\l(p^{m+\frac{1}{2}},\nabla\cdot \bm{\varphi}^h\r) \notag\\
    &\qquad\qquad -\gamma\l<\kappa^{m+\frac{1}{2}}\widetilde{\bnu}^{m+\frac{1}{2}},\bm{\varphi}^h\r>_{\widetilde{\Gamma}^{m+\frac{1}{2}}}-\l(\widetilde{\bg}^{m+\frac{1}{2}},\bm{\varphi}^h \r),\\
    0&=\l(\nabla\cdot \bu^{m+\frac{1}{2}},q^h \r),\\
    0&=\l<\l(\frac{\bX^{m+1}-\bX^m}{\Delta t} \r)\cdot \widetilde{\bnu}^{m+\frac{1}{2}},\eta^h\r>_{\widetilde{\Gamma}^{m+\frac{1}{2}}}^h-\l<\bu^{m+\frac{1}{2}}\cdot \widetilde{\bnu}^{m+\frac{1}{2}},\eta^h\r>_{\widetilde{\Gamma}^{m+\frac{1}{2}}} \notag\\
    &\qquad\qquad-\lambda_{\cE}^{m+\frac{1}{2}}\l<\kappa^{m+\frac{1}{2}},\eta^h \r>_{\widetilde{\Gamma}^{m+\frac{1}{2}}}^h-\lambda_{V}^{m+\frac{1}{2}}\l<1,\eta^{h} \r>_{\widetilde{\Gamma}^{m+\frac{1}{2}}}^h, \\
    0&=\l<\kappa^{m+\frac{1}{2}}\widetilde{\bnu}^{m+\frac{1}{2}},\bm{\omega}^h\r>_{\widetilde{\Gamma}^{m+\frac{1}{2}}}^h+\l<\nabla_{\widetilde{\Gamma}^{m+\frac{1}{2}}}\bX^{m+\frac{1}{2}},\nabla_{\widetilde{\Gamma}^{m+\frac{1}{2}}}\bm{\omega}^h\r>_{\widetilde{\Gamma}^{m+\frac{1}{2}}},\\
    0&=\frac{\cE^{m+1}-\cE^{m}}{\Delta t}+\l(2\widetilde{\mu}^{m+\frac{1}{2}}\bD(\bu^{m+\frac{1}{2}}),\bD(\bu^{m+\frac{1}{2}}) \r)-\l(\,\bu^{m+\frac{1}{2}}, \widetilde{\bg}^{m+\frac{1}{2}}\r),\label{Lag:CN5} \\
    0&=V^{m+1}-V^m,	\label{Lag:CN6}
\end{align}
\end{subequations}
for any $(\bm{\varphi}^h,q^h,\eta^h,\bm{\omega}^h)\in \U^m\times \widehat{\P}^m\times \widetilde{W}^{m+\frac12}\times [\widetilde{W}^{m+\frac12}]^d$, and then set $\Gamma^{m+1} = \bX^{m+1}(\widetilde{\Gamma}^{m+\frac{1}{2}})$.
Here, $\widetilde{W}^{m+\frac12}$ refers to the linear finite element space on $\widetilde{\Gamma}^{m+\frac12}$.
Concerning suitable examples for the predictor surface $\widetilde{\Gamma}^{m+\frac{1}{2}}$, we refer to Section~\ref{sec:predictors}. We note that the normal vector $\widetilde{\bnu}^{m+\frac{1}{2}}$ and the discrete viscosity $\widetilde{\mu}^{m+\frac{1}{2}}$ are both computed based on the predicted surface $\widetilde{\Gamma}^{m+\frac{1}{2}}$ and are therefore treated as known quantities. As a result, they do not introduce any additional nonlinearity into the system, also compare with the first-order methods \eqref{Lag:Euler} and \eqref{BGN}, where similar quantities are defined using the interface from the previous time step. Moreover, we use the mid-time approximation in the following sense:
\begin{align*}
    &\bX^{m+\frac{1}{2}}=\frac{\bX^{m+1}+\bX^{m}}{2},\qquad \bu^{m+\frac{1}{2}}=\frac{\bu^{m+1}+\bu^{m}}{2},\qquad \kappa^{m+\frac{1}{2}}=\frac{\kappa^{m+1}+\kappa^{m}}{2},\\
     &p^{m+\frac{1}{2}}=\frac{p^{m+1}+p^{m}}{2},\qquad \lambda_{\cE}^{m+\frac{1}{2}}=\frac{\lambda_{\cE}^{m+1}+\lambda_{\cE}^{m}}{2},\qquad\lambda_{V}^{m+\frac{1}{2}}=\frac{\lambda_{V}^{m+1}+\lambda_{V}^{m}}{2},
\end{align*}
where $(\bu^{m},p^{m}, \bX^{m},\kappa^{m})\in \U^m\times \widehat{\P}^m\times [\widetilde{W}^{m+\frac12}]^d\times \widetilde{W}^{m+\frac12}$ as well as $\lambda_{\cE}^m, \lambda_{V}^m\in\mathbb{R}$ are given from the previous time step. We point out that if the bulk mesh at the previous time step $t_m$ differs from the current mesh at $t_{m+1}$, as it is typically the case in the fitted mesh approach, then the velocity $\bu^m$ and pressure $p^m$ have to be appropriately interpolated or projected onto the current bulk mesh in practice.
Besides, the discrete body force is defined as $\widetilde{\bg}^{m+\frac{1}{2}}=I_2^h\bg(t_{m+\frac{1}{2}})$. In the case that the body force depends on the geometry, we use the predictor step to compute the geometry needed to define $\widetilde{\bg}^{m+\frac{1}{2}}$.

\subsubsection{BDF2 scheme}
We now consider another second temporal order discretization based on the second-order backward differentiation formula (BDF2).

Let $m\in\{1,\ldots,M-1\}$. Given $\Gamma^{m-1},\Gamma^m$ and an approximation $\widetilde{\Gamma}^{m+1}$ of $\Gamma(t_{m+1})$, we aim to find a solution
\[
(\bu^{m+1},p^{m+1}, \bX^{m+1},\kappa^{m+1})\in \U^m\times \widehat{\P}^m\times [\widetilde{W}^{m+1}]^d\times \widetilde{W}^{m+1}
\]
and Lagrange multipliers $(\lambda_{\cE}^{m+1},\lambda_{V}^{m+1})\in \R\times \R$, such that 
\begin{subequations}\label{Lag:BDF2}
\begin{align}
    0&=\l(2\widetilde{\mu}^{m+1}\bD(\bu^{m+1}),\bD(\bm{\varphi}^h) \r)-\l(p^{m+1},\nabla\cdot \bm{\varphi}^h\r)\notag
    \\
    &\qquad\qquad -\gamma\l<\kappa^{m+1}\widetilde{\bnu}^{m+1},\bm{\varphi}^h\r>_{\widetilde{\Gamma}^{m+1}}-\l(\widetilde{\bg}^{m+1},\bm{\varphi}^h \r),
    \\
    0&=\l(\nabla\cdot \bu^{m+1},q^h \r),
    \\
    0&=\l<\l(\frac{3\bX^{m+1}-4\bX^m+\bX^{m-1}}{2\Delta t} \r)\cdot \widetilde{\bnu}^{m+1},\eta^h\r>_{\widetilde{\Gamma}^{m+1}}^h-\l<\bu^{m+1}\cdot\widetilde{\bnu}^{m+1}, \eta^h\r>_{\widetilde{\Gamma}^{m+1}}\notag
    \\
    &\qquad\qquad-\lambda_{\cE}^{m+1}\l<\kappa^{m+1},\eta^h \r>_{\widetilde{\Gamma}^{m+1}}^h-\lambda_{V}^{m+1}\l<1,\eta^{h} \r>_{\widetilde{\Gamma}^{m+1}}^h, 
    \\
    0&=\l<\kappa^{m+1}\widetilde{\bnu}^{m+1},\bm{\omega}^h\r>_{\widetilde{\Gamma}^{m+1}}^h+\l<\nabla_{\widetilde{\Gamma}^{m+1}}\bX^{m+1},\nabla_{\widetilde{\Gamma}^{m+1}}\bm{\omega}^h\r>_{\widetilde{\Gamma}^{m+1}},\\
    0&=\frac{3\cE^{m+1}-4\cE^{m}+\cE^{m-1}}{2\Delta t}+\l(2\widetilde{\mu}^{m+1}\bD(\bu^{m+1}),\bD(\bu^{m+1}) \r)-\l(\,\bu^{m+1}, \widetilde{\bg}^{m+1}\r), \label{Lag:BDF5}\\
    0&=3V^{m+1}-4V^m+V^{m-1},\label{Lag:BDF6}
\end{align}
\end{subequations}
for any $(\bm{\varphi}^h,q^h,\eta^h,\bm{\omega}^h)\in \U^m\times \widehat{\P}^m\times \widetilde{W}^{m+1}\times [\widetilde{W}^{m+1}]^d$, and then set $\Gamma^{m+1} = \bX^{m+1}(\widetilde{\Gamma}^{m+1})$. 
Similar to before, $\widetilde{W}^{m+1}$ refers to the linear finite element space on $\widetilde{\Gamma}^{m+1}$.
Here, we employ a BDF2-type discretization for the time derivatives, which requires information from the previous two time steps, namely the surfaces $\Gamma^{m-1}$ and $\Gamma^m$. 
To initialize the BDF2 scheme \eqref{Lag:BDF2}, two initial surfaces $\Gamma^0$ and $\Gamma^1$ are needed in order to compute the first update $\Gamma^2$. Typically, $\Gamma^0$ is prescribed as the initial condition, while $\Gamma^1$ is obtained by applying the first-order scheme \eqref{Lag:Euler} with time step size $\Delta t$. To improve temporal accuracy, $\Gamma^1$ can instead be computed using one step of the Crank--Nicolson scheme \eqref{Lag:CN} or another time integration method of second-order accuracy. 
Similarly, previous values (or suitable initial guesses) for the energy $\mathcal{E}^{m-1}$, $\mathcal{E}^m$ and the volume $V^{m-1}$, $V^m$ are required.
Also, the predictor surface $\widetilde{\Gamma}^{m+1}$ can be determined using the strategies in Section~\ref{sec:predictors}. 
Similar to the Crank--Nicolson scheme \eqref{Lag:CN}, the normal vector $\widetilde{\bnu}^{m+1}$ and the viscosity $\widetilde{\mu}^{m+1}$ are both computed based on the predicted surface $\widetilde{\Gamma}^{m+1}$ and are hence treated as known quantities.
Moreover, the discrete body force is defined by $\widetilde{\bg}^{m+1}=I_2^h\bg\l(t_{m+1}\r)$.

\subsubsection{Structure-preserving properties}
Our key motivation for introducing the Lagrange multiplier method is that it enables fully discrete schemes which exactly preserve the physical structures of the two-phase Stokes equation. In particular, the resulting discretizations inherit both the energy dissipation property (in the absence of external forces) and the conservation of the enclosed volume. These properties are summarized in the following theorem. 

\begin{theorem}
For all $m=0,\ldots,M-1$, let $(\bu^{m+1},p^{m+1},\bX^{m+1},\kappa^{m+1})$ and $(\lambda^{m+1}_{\cE},\lambda^{m+1}_V)$ be a solution of \eqref{Lag:Euler}, \eqref{Lag:CN}, or \eqref{Lag:BDF2}. Then, it holds
\begin{equation}
    \cE^{m+1}\le \cE^{m}\le \ldots\le  \cE^{0},\quad \text{if}\ \bg\equiv 0,
\end{equation}
and 
\begin{equation}
    V^{m+1}=V^{m}=\ldots=V^0.
\end{equation}	
\end{theorem}

\Proof
The claims for the Euler scheme \eqref{Lag:Euler} and the Crank--Nicolson scheme \eqref{Lag:CN} are direct consequences of the enforced equations \eqref{Lag:Euler5}--\eqref{Lag:Euler6} and \eqref{Lag:CN5}--\eqref{Lag:CN6}, respectively. For the BDF2 scheme \eqref{Lag:BDF2}, the statement follows by induction, where we refer to \cite[Theorem~3.4]{Garcke2024} for a similar proof.
\qed

\subsubsection{Choice of surface predictors} \label{sec:predictors}

To apply higher-order time discretization methods such as the Crank--Nicolson scheme \eqref{Lag:CN} or BDF2 schemes \eqref{Lag:BDF2}, respectively, we require a predictor for the interface geometry at the intermediate or new time levels $t_{m+\frac12}$ or $t_{m+1}$, respectively. In this section, we describe the construction of these predictors and explain the reasoning behind their use.

A predictor surface is necessary because computing certain geometric quantities, such as the surface normal vector or viscosity function, at the new or mid time level requires an approximation of the evolving interface. A naive approach might be to extrapolate the interface geometry from previous time steps. However, such extrapolations often lead to severe mesh distortions in practice, as observed for other surface evolution problems \cite{Jiang2024b,Jiang2025}. In contrast, predictors based on first-order evolution yield more stable and mesh-regular parametrizations of the interface. We refer to \cite{Garcke2024,Jiang2024a,Jiang2024b,Jiang2025} for more detailed discussion of the choice of surface predictors.

For the Crank--Nicolson scheme \eqref{Lag:CN}, we define the predictor surface $\widetilde{\Gamma}^{m+\frac{1}{2}}$ with the help of a first-order method with a half time step size $\frac12 \Delta t$. The associated unit normal vector $\widetilde{\bnu}^{m+\frac12}$ and discrete viscosity $\widetilde{\mu}^{m+\frac12}$ are both computed using this predicted surface. Simlarly, for the BDF2 scheme \eqref{Lag:BDF2}, the predictor surface $\widetilde{\Gamma}^{m+1}$ is generated by a first-order scheme with full time step size $\Delta t$, and the corresponding quantities $\widetilde{\bnu}^{m+1}$ and $\widetilde{\mu}^{m+1}$ are determined from this predicted surface.

In the context of the fitted mesh approach, it seems also natural to fit the bulk mesh to these predicted interfaces $\widetilde{\Gamma}^{m+\frac12}$ or $\widetilde{\Gamma}^{m+1}$, respectively, to ensure consistent geometric alignment between the bulk and interface meshes.

We propose two practical options for constructing the predictors: firstly, the nonlinear Euler scheme \eqref{Lag:Euler}, or secondly, the fully linear BGN scheme \eqref{BGN}. On the one hand, using the BGN scheme \eqref{BGN} for prediction can reduce computational cost, since it is fully linear and thus less expensive to solve. On the other hand, the nonlinear Euler scheme \eqref{Lag:Euler} may offer improved geometric accuracy, since it preserves the enclosed volume exactly.
However, in line with observations from \cite{Garcke2024} for the surface diffusion model, our numerical experiments suggest that the choice between the Euler and BGN schemes has only a minor influence on the overall accuracy of the method, including in convergence studies.

\subsection{Properties of discrete stationary solution}

In this part, we study the properties of the discrete stationary solution of our schemes based on the Lagrange multiplier method. We take the Euler scheme \eqref{Lag:Euler} as an example, and we note that the argument can be simply extended to other discretization methods \eqref{Lag:CN} and \eqref{Lag:BDF2}.

\begin{proposition}
Let $(\bu^{m+1},p^{m+1},\bX^{m+1},\kappa^{m+1})$ and $(\lambda^{m+1}_{\cE},\lambda^{m+1}_V)$ be a solution of \eqref{Lag:Euler} with $\bg\equiv 0$. We further assume the following conditions:
\begin{enumerate}
\item[(i)] As $m\rightarrow+\infty$, we suppose that $(\bu^{m+1},p^{m+1},\bX^{m+1},\kappa^{m+1})$ converges to a stationary state $(\bu^{e},p^{e},\bX^{e},\kappa^{e})\in \U^e\times \widehat{\P}^e\times [W^e]^d\times W^e$, and that $(\lambda^{m+1}_{\cE},\lambda^{m+1}_V)$ converges to some $(\lambda^{e}_{\cE},\lambda^{e}_V)$;
\item[(ii)] The stationary polygon/polyhedra $\Gamma^e$ is 
non-degenerate in the sense of Section~\ref{sec:discretization}, i.e., $\Gamma^e$ has no self-intersections and all its interface elements have positive Hausdorff measure.
\end{enumerate}
Then, exactly one of the following situations occurs: 
\begin{itemize}
\item The discrete curvature $\kappa^{e}$ is not constant. In this case, the Euler scheme \eqref{Lag:Euler} yields a stationary solution with $\lambda_{\cE}^{e} = \lambda_{V}^{e} = 0$, which coincides with a stationary solution of the BGN scheme \eqref{BGN}, and vice versa.

\item The discrete curvature $\kappa^{e}\in\mathbb{R}$ is constant. Then, the stationary interface $\Gamma^e$ is a conformal polyhedral surface in the sense of \cite[Def.~60]{BGN20}. If $d=2$, then $\Gamma^e$ is also weakly equidistributed, i.e., any two neighbouring elements of the curve $\Gamma^e$ have equal length or are parallel.
\end{itemize}
\end{proposition}

\Proof
Under the assumption (i), the scheme \eqref{Lag:Euler} converges to 
\begin{subequations}\label{stationary}
\begin{align}
    0 &= \l(2\mu^{e}\bD(\bu^{e}),\bD(\bm{\varphi}^h) \r) -\l(p^{e},\nabla\cdot \bm{\varphi}^h\r)-\gamma\l<\kappa^{e}\bnu^e,\bm{\varphi}^h\r>_{\Gamma^e}, \label{stationary1}
    \\
    0&=\l(\nabla\cdot \bu^{e},q^h \r),
    \\
    0 &=  -\l<\bu^{e}\cdot \bnu^e,\eta^h\r>_{\Gamma^e}  -\lambda_{\cE}^{e}\l<\kappa^{e},\eta^h \r>_{\Gamma^e}^h-\lambda_{V}^{e}\l<1,\eta^{h} \r>_{\Gamma^e}^h,
    \\
    0 &= \l<\kappa^{e}\bnu^e,\bm{\omega}^h\r>_{\Gamma^e}^h+\l<\nabla_{\Gamma^e}\bX^{e},\nabla_{\Gamma^e}\bm{\omega}^h\r>_{\Gamma^e}, \qquad\quad\qquad\quad\quad\quad
    \\
    0&
    =\l(2\mu^e\bD(\bu^{e}),\bD(\bu^{e}) \r),\label{stationary6}
\end{align}
\end{subequations}
where the interface $\Gamma^e$ is non-degenerate due to the assumption (ii).
By setting $\bm{\varphi}^h=\bu^{e}$, $q^h=p^e$ and $\eta^h=\kappa^e$, we obtain 
\[
\lambda_{\cE}^{e}\l<\kappa^{e},\kappa^e \r>_{\Gamma^e}^h+\lambda_{V}^{e}\l<1,\kappa^e \r>_{\Gamma^e}^h=0.
\]
Similarly, by taking $\eta^h=1$ and by using the stationary property, we have 
\[
\lambda_{\cE}^{e}\l<\kappa^{e},1 \r>_{\Gamma^e}^h+\lambda_{V}^{e}\l<1,1 \r>_{\Gamma^e}^h=-\l<\bu^{e}\cdot \bnu^e,1\r>_{\Gamma^e}=0.
\]
Combining both identities and applying a simple argument using the Cauchy--Schwarz inequality, one can show that exactly one of these two cases must occur: On the one hand, when $\kappa^{e}$ is not constant, then necessarily $\lambda_{\cE}^{e}=\lambda_{V}^{e}=0$, and the Euler scheme reduces to the BGN scheme. On the other hand, if $\kappa^{e} \in\mathbb{R}$ is constant, then we have
\[
0=\kappa^{e}\l<\bnu^e,\bm{\omega}^h\r>_{\Gamma^e}^h+\l<\nabla_{\Gamma^e}\bX^{e},\nabla_{\Gamma^e}\bm{\omega}^h\r>_{\Gamma^e}.
\]
Note that this equation directly implies that $\Gamma^e$ satisfies the definition of a conformal polyhedral surfaces in the sense of \cite[Def.~60]{BGN20}, since $\Gamma^e$ is non-degenerate.
Additionally, in the two dimensional case $d=2$, the weak equidistribution property follows from \cite[Theorem~62]{BGN20}.
\qed

The above argument can be easilty extended to the other time discretizations and the single Lagrange multiplier approach (see Section~\ref{sec:Single}).

\subsection{Newton iterations and solution technique} \label{sec:Newton}

In this section, we introduce an efficient Newton iteration method and a decoupling technique to solve the nonlinear systems \eqref{Lag:Euler}, \eqref{Lag:CN}, or \eqref{Lag:BDF2}, respectively.

We first recall the first variation formula of for the energy and volume contributions.

\begin{lemma}\label{Lemma:FirstVariation}
Let $\Gamma^h$ be globally parametrized by $\bX^h$. For $\epsilon>0$ small enough, we denote by $\Gamma^h(\epsilon)$ the closed polyhedral surface parametrized by $\bX^h(\epsilon) = \bX^h + \epsilon \bX^\delta$, for some $\bX^\delta$. Then, it holds
\begin{align}
    \lim_{\epsilon\to 0} \frac{\d}{\d\epsilon} 
    \int_{\Gamma^h(\epsilon)} 1\, \d \cH^{d-1}
    = \l<\nabla_{\Gamma^h}\mathrm{Id},\nabla_{\Gamma^h}\bX^\delta \r>_{\Gamma^h},
    \qquad 
    \lim_{\epsilon\to 0} \frac{\d}{\d\epsilon} 
    V^h(\bX^h (\epsilon))
    =\l<\bX^{\delta},\bnu^{h}\r>_{\Gamma^h}^h, 
\end{align}
where $V^h(\bX^h (\epsilon))$ denotes the volume enclosed by $\Gamma^h(\epsilon)$.
\end{lemma}
\Proof
This directly follows from the transport theorems \cite[Theorems~70 and~71]{BGN20} together with \cite[Lemma~9]{BGN20}).
\qed

\smallskip

We now take the Euler scheme \eqref{Lag:Euler} as an example to illustrate our Newton iterations. The other two schemes \eqref{Lag:CN} and \eqref{Lag:BDF2} can be treated similarly. At each $i$-th iteration step, we aim to find a Newton direction 
\[
(\bu^{\delta},p^{\delta},\bX^{\delta},\kappa^{\delta})\in \U^m\times \widehat{\P}^m\times [W^m]^d\times W^m
\]
and $(\lambda_{\cE}^{\delta},\lambda_{V}^{\delta})\in \R\times \R$, such that 
\begin{subequations}\label{Lag:Newton}
\begin{align}
    F_1(\bm{\varphi}^h)&=\l(2\mu^{m}\bD(\bu^{\delta}),\bD(\bm{\varphi}^h) \r) -\l(p^{\delta},\nabla\cdot \bm{\varphi}\r)-\gamma\l<\kappa^{\delta}\bnu^m,\bm{\varphi}^h\r>_{\Gamma^m}, \label{Lag:Newton1}
    \\
    F_2(q^h)&=\l(\nabla\cdot \bu^{\delta},q^h \r), \\F_3(\eta^h)&=\l<\frac{\bX^{\delta}}{\Delta t} \cdot \bnu^m,\eta^h\r>_{\Gamma^m}^h-\l<\bu^{\delta}\cdot\bnu^{m}, \eta^h\r>_{\Gamma^m}\notag
    \\
    &\qquad\qquad -\lambda_{\cE}^{\delta}\l<\kappa^{m+1,i},\eta^h \r>_{\Gamma^m}^h-\lambda_{\cE}^{m+1,i}\l<\kappa^{\delta},\eta^h \r>_{\Gamma^m}^h-\lambda_{V}^{\delta}\l<1,\eta^{h} \r>_{\Gamma^m}^h, 
    \\
    F_4(\bm{\omega}^h)&=\l<\kappa^{\delta}\bnu^m,\bm{\omega}^h\r>_{\Gamma^m}^h+\l<\nabla_{\Gamma^m}\bX^{\delta},\nabla_{\Gamma^m}\bm{\omega}^h\r>_{\Gamma^m}, \label{Lag:Newton4}
    \\
    F_5&=\frac{\gamma}{\Delta} \l<\nabla_{\Gamma^{m+1,i}}\mathrm{Id}, \nabla_{\Gamma^{m+1,i}}\bX^\delta\r>_{\Gamma^{m+1,i}} + \l(4\mu^m\bD(\bu^{\delta}),\bD(\bu^{m+1,i}) \r)-\l(\,\bu^{\delta}, \widetilde{\bg}^{m+1}\r),\label{Lag:Newton5}
    \\
    F_6&=\l<\bX^{\delta} ,\bnu^{m+1,i}\r>_{\Gamma^{m+1,i}}^h,\label{Lag:Newton6}
\end{align}
\end{subequations}
for any $(\bm{\varphi}^h,q^h,\eta^h,\bm{\omega}^h)\in \U^m\times \widehat{\P}^m\times W^m\times [W^m]^d$,
where the terms $F_1,\ldots,F_6$ on the left-hand sides in \eqref{Lag:Newton} are defined as follows:
\begin{subequations}
\begin{align}
    F_1(\bm{\varphi}^h)&\colonequals-\l(2\mu^{m}\bD(\bu^{m+1,i}),\bD(\bm{\varphi}^h) \r) +\l(p^{m+1,i},\nabla\cdot \bm{\varphi}^h\r)\notag
    \\
    &\qquad\qquad +\gamma\l<\kappa^{m+1,i}\bnu^m,\bm{\varphi}^h\r>_{\Gamma^m} +(\widetilde{\bg}^{m+1},\bm{\varphi}^h),
    \\
    F_2(q^h)&\colonequals-\l(\nabla\cdot \bu^{m+1,i},q^h \r), \\F_3(\eta^h)&\colonequals-\l<\l(\frac{\bX^{m+1,i}-\bX^m}{\Delta t} \r)\cdot \bnu^m,\eta^h\r>_{\Gamma^m}^h+\l<\bu^{m+1,i}\cdot\bnu^{m}, \eta^h\r>_{\Gamma^m}\notag
    \\
    &\qquad\qquad +\lambda_{\cE}^{m+1,i}\l<\kappa^{m+1,i},\eta^h \r>_{\Gamma^m}^h+\lambda_{V}^{m+1,i}\l<1,\eta^{h} \r>_{\Gamma^m}^h, 
    \\
    F_4(\bm{\omega}^h)&\colonequals-\l<\kappa^{m+1,i}\bnu^m,\bm{\omega}^h\r>_{\Gamma^m}^h-\l<\nabla_{\Gamma^m}\bX^{m+1,i},\nabla_{\Gamma^m}\bm{\omega}^h\r>_{\Gamma^m},
    \\
    F_5&\colonequals-\frac{\cE^{m+1,i}-\cE^{m}}{\Delta t}-\l(2\mu^m\bD(\bu^{m+1,i}),\bD(\bu^{m+1,i}) \r)+\l(\bu^{m+1,i}, \widetilde{\bg}^{m+1}\r),
    \\
    F_6&\colonequals-(V^{m+1,i}-V^m).
\end{align}
\end{subequations}
We note that the surface integrals in \eqref{Lag:Newton5}--\eqref{Lag:Newton6} are evaluated over the updated interface $\Gamma^{m+1,i} = \bX^{m+1,i}(\Gamma^m)$, with the corresponding discrete normal vector $\bnu^{m+1,i}$. Likewise, $\mathcal{E}^{m+1,i}$ denotes the energy associated with $\Gamma^{m+1,i}$, and $V^{m+1,i}$ refers to the volume of the region $\Omega^{m+1,i}_-$ enclosed by this interface. In contrast, the surface integrals in \eqref{Lag:Newton1}--\eqref{Lag:Newton4} are always computed on the fixed interface $\Gamma^m$.

We choose the initial guess for the Newton iterations as
\[
\bu^{m+1,0}=\bu^m,\ p^{m+1,0}=p^m,\ \bX^{m+1,0}=\bX^m,\ \kappa^{m+1,0}=\kappa^m,\ \lambda_{\cE}^{m+1,0}=\lambda_{\cE}^m,\ \lambda_{V}^{m+1,0}=\lambda_{V}^m,
\]
and the iteration is updated as
\begin{align*}
    \bu^{m+1,i+1}&= \bu^{m+1,i}+\bu^\delta,\qquad\ \  p^{m+1,i+1}= p^{m+1,i}+p^\delta,\\
     \bX^{m+1,i+1}&= \bX^{m+1,i}+\bX^\delta,\qquad \kappa^{m+1,i+1}= \kappa^{m+1,i}+\kappa^\delta,\\
     \lambda_{\cE}^{m+1,i+1}&= \lambda_{\cE}^{m+1,i}+\lambda_{\cE}^\delta, \qquad\ \ \lambda_{V}^{m+1,i+1}= \lambda_{V}^{m+1,i}+\lambda_{V}^\delta,
\end{align*}
until the stopping criterium
\[
\max\{\|\bu^{\delta}\|_{L^{\infty}},\|p^{\delta}\|_{L^{\infty}},\|\bX^{\delta}\|_{L^{\infty}},\|\kappa^{\delta}\|_{L^{\infty}},|\lambda_{\cE}^{\delta}|,|\lambda_{V}^{\delta}|\}\le \mathrm{tol}
\]
is met for some given tolerance $\mathrm{tol}=10^{-10}$. Once the stopping criterium is satisfied, the (approximate) numerical solution of the nonlinear system \eqref{Lag:Euler} is defined by
\begin{align*}
    \bu^{m+1}&=\bu^{m+1,i+1},\qquad p^{m+1}=p^{m+1,i+1},\\ \bX^{m+1}&=\bX^{m+1,i+1},\qquad \kappa^{m+1}=\kappa^{m+1,i+1},\\ \lambda_{\cE}^{m+1}&=\lambda_{\cE}^{m+1,i+1},\qquad \lambda_{V}^{m+1}=\lambda_{V}^{m+1,i+1}.
\end{align*}
In our numerical experiments, the stopping criterion was typically satisfied after only a few iterations.

While the linear systems in \eqref{Lag:Newton} can in principle be solved directly, finding effective preconditioners for iterative solvers may pose a challenge due to the coupled structure of the equations. To address this, we propose a decoupling technique motivated by, e.g, \cite[Section 3.2.1]{li_shen_2022_msav_CHNSt}, that simplifies the solution process. We decompose the variables $(\bu^\delta,p^\delta,\bX^\delta,\kappa^\delta)$ in terms of the scalars $(\lambda^\delta_{\cE},\lambda^\delta_{V})$ 
\begin{align*}
    \bu^\delta
    &=: \bu^{(0)}+\lambda^\delta_{\cE} \bu^{(1)} +\lambda^\delta_{V} \bu^{(2)}, \qquad\qquad 
    p^\delta =: p^{(0)}+ \lambda^\delta_{\cE} p^{(1)} +\lambda^\delta_{V} p^{(2)},
    \\
    \bX^\delta 
    &=: \bX^{(0)}+\lambda^\delta_{\cE} \bX^{(1)} +\lambda^\delta_{V} \bX^{(2)}, \qquad\quad
    \kappa^\delta =: \kappa^{(0)} + \lambda^\delta_{\cE} \kappa^{(1)} + \lambda^\delta_{V}\kappa^{(2)} ,
\end{align*}
and we compute $(\bu^{(j)},p^{(j)},\bX^{(j)},\kappa^{(j)})$, $j=0,1,2$, by three sparse linear systems that only involve the subsystem \eqref{Lag:Newton1}--\eqref{Lag:Newton4}. Afterwards, we reconstruct $(\lambda^\delta_{\cE},\lambda^\delta_{V})$ with the help of \eqref{Lag:Newton5}--\eqref{Lag:Newton6} using $(\bu^{(j)},p^{(j)},\bX^{(j)},\kappa^{(j)})$, $j=0,1,2$. The three linear systems for $(\bu^{(j)},p^{(j)},\bX^{(j)},\kappa^{(j)})$, $j=0,1,2$ are given as follows.
\begin{itemize}

\item Leading term of order $(1,1)$ in \eqref{Lag:Newton1}--\eqref{Lag:Newton4}:
\begin{subequations}\label{Decouple1}
\begin{align}
    F_1(\bm{\varphi}^h) &=\l(2\mu^{m}\bD(\bu^{(0)}),\bD(\bm{\varphi}^h) \r) -\l(p^{(0)},\nabla\cdot \bm{\varphi}^h\r)-\gamma\l<\kappa^{(0)}\bnu^m,\bm{\varphi}^h\r>_{\Gamma^m},
    \\
    F_2(q^h)&=\l(\nabla\cdot \bu^{(0)},q^h \r), \\F_3(\eta^h)&=\l<\frac{\bX^{(0)}}{\Delta t} \cdot \bnu^m,\eta^h\r>_{\Gamma^m}^h-\l<\bu^{(0)}\cdot\bnu^{m}, \eta^h\r>_{\Gamma^m} -\lambda_{\cE}^{m+1,i}\l<\kappa^{(0)},\eta^h \r>_{\Gamma^m}^h, 
    \\
    F_4(\bm{\omega}^h) &=\l<\kappa^{(0)}\bnu^m,\bm{\omega}^h\r>_{\Gamma^m}^h+\l<\nabla_{\Gamma^m}\bX^{(0)},\nabla_{\Gamma^m}\bm{\omega}^h\r>_{\Gamma^m} .
\end{align}
\end{subequations}
\item Leading term of order $(\lambda^\delta_{\cE},1)$ in \eqref{Lag:Newton1}--\eqref{Lag:Newton4}:
\begin{subequations}\label{Decouple2}
\begin{align}
    0&=\l(2\mu^{m}\bD(\bu^{(1)}),\bD(\bm{\varphi}^h) \r) -\l(p^{(1)},\nabla\cdot \bm{\varphi}^h\r)-\gamma\l<\kappa^{(1)}\bnu^m,\bm{\varphi}^h\r>_{\Gamma^m},
    \\
    0&=\l(\nabla\cdot \bu^{(1)},q^h \r), \\
    \l<\kappa^{m+1,i},\eta^h \r>_{\Gamma^m}^h& = \l<\frac{\bX^{(1)}}{\Delta t} \cdot \bnu^m,\eta^h\r>_{\Gamma^m}^h-\l<\bu^{(1)}\cdot\bnu^{m}, \eta^h\r>_{\Gamma^m}
    -\lambda^{m+1,i} \l<\kappa^{(1)},\eta^h \r>_{\Gamma^m}^h, 
    \\
    0&= \l<\kappa^{(1)}\bnu^m,\bm{\omega}^h\r>_{\Gamma^m}^h+\l<\nabla_{\Gamma^m}\bX^{(1)},\nabla_{\Gamma^m}\bm{\omega}^h\r>_{\Gamma^m} .
\end{align}
\end{subequations}
\item Leading term of order $(1,\lambda^\delta_{V})$ in \eqref{Lag:Newton1}--\eqref{Lag:Newton4}:
  \begin{subequations}\label{Decouple3}
\begin{align}
    0&=\l(2\mu^{m}\bD(\bu^{(2)}),\bD(\bm{\varphi}^h) \r) -\l(p^{(2)},\nabla\cdot \bm{\varphi}^h\r)-\gamma\l<\kappa^{(2)}\bnu^m,\bm{\varphi}^h\r>_{\Gamma^m},
    \\
    0&=\l(\nabla\cdot \bu^{(2)},q^h \r), \\
    \l<1,\eta^h\r>^h_{\Gamma^m} &=\l<\frac{\bX^{(2)}}{\Delta t} \cdot \bnu^m,\eta^h\r>_{\Gamma^m}^h-\l<\bu^{(2)}\cdot\bnu^{m}, \eta^h\r>_{\Gamma^m} -\lambda^{m+1,i} \l<\kappa^{(2)},\eta^h \r>_{\Gamma^m}^h, 
    \\
    0&= \l<\kappa^{(2)}\bnu^m,\bm{\omega}^h\r>_{\Gamma^m}^h+\l<\nabla_{\Gamma^m}\bX^{(2)},\nabla_{\Gamma^m}\bm{\omega}^h\r>_{\Gamma^m} .
\end{align}
\end{subequations}
\end{itemize}
Subsequently, by substituting $\bu^\delta=\bu^{(0)} + \lambda^\delta_{\cE} \bu^{(1)} + \lambda^\delta_{V} \bu^{(2)}$ and $\bX^\delta=\bX^{(0)}+ \lambda^\delta_{\cE} \bX^{(1)} + \lambda^\delta_{V}\bX^{(2)}$ into the enforced scalar equations \eqref{Lag:Newton5} and \eqref{Lag:Newton6}, we then obtain the following linear $(2\times 2)$-system for $(\lambda^\delta_{\cE},\lambda^\delta_{V})$:
\begin{equation}\label{Reconstruct}
    \begin{pmatrix}
        A_{11} & A_{12}\\
        A_{21} & A_{22}
    \end{pmatrix}\begin{pmatrix}
        \lambda^\delta_{\cE}\\
        \lambda^\delta_{V}
    \end{pmatrix}=\begin{pmatrix}
        g_1\\
        g_2
    \end{pmatrix},
\end{equation}
where
\begin{subequations}
\begin{align}
    A_{11}&\colonequals   \frac{\gamma}{\Delta t} \l<\nabla_{\Gamma^{m+1,i}}\mathrm{Id}, \nabla_{\Gamma^{m+1,i}}\bX^{(1)}\r>_{\Gamma^{m+1,i}}+ \l(4\mu^m\bD(\bu^{(1)}),\bD(\bu^{m+1,i}) \r) \notag
    \\
    &\qquad\qquad -\l(\,\bu^{(1)}, \widetilde{\bg}^{m+1}\r),
    \\
    A_{12}&\colonequals  \frac{\gamma}{\Delta t} \l<\nabla_{\Gamma^{m+1,i}}\mathrm{Id},\nabla_{\Gamma^{m+1,i}}\bX^{(2)}\r>_{\Gamma^{m+1,i}}+ \l(4\mu^m\bD(\bu^{(2)}),\bD(\bu^{m+1,i}) \r) \notag
    \\
    &\qquad\qquad  -\l(\,\bu^{(2)}, \widetilde{\bg}^{m+1}\r),
    \\
    A_{21}&\colonequals  \l<\bX^{(1)} ,\bnu^{m+1,i}\r>_{\Gamma^{m+1,i}}^h,
    \\
    A_{22}&\colonequals  \l<\bX^{(2)} ,\bnu^{m+1,i}\r>_{\Gamma^{m+1,i}}^h,
    \\
    g_{1}&\colonequals F_5- \frac{\gamma}{\Delta t} \l<\nabla_{\Gamma^{m+1,i}}\mathrm{Id},\nabla_{\Gamma^{m+1,i}}\bX^{(0)}\r>_{\Gamma^{m+1,i}} - \l(4\mu^m\bD(\bu^{(0)}),\bD(\bu^{m+1,i}) \r) \notag
    \\
    &\qquad\qquad +\l(\,\bu^{(0)}, \widetilde{\bg}^{m+1}\r),
    \\
    g_{2}&\colonequals F_6-\l<\bX^{(0)} ,\bnu^{m+1,i}\r>_{\Gamma^{m+1,i}}^h.
\end{align}
\end{subequations}
Solving the three linear systems \eqref{Decouple1}, \eqref{Decouple2}, and \eqref{Decouple3} effectively corresponds to solving a single linear system with three different right-hand sides. The overall computational strategy can be further simplified by employing a Schur complement approach, as proposed in \cite[Section 5]{BGN_2015_navierstokes}, which enables a decoupling between the bulk variables $\bu^{(j)}$, $p^{(j)}$ and the interface variables $\bX^{(j)}$, $\kappa^{(j)}$, $j=0,1,2$. This decoupling further reduces the problem to solving two smaller linear systems.

An important observation is that the only matrix block that changes within the Newton iterations corresponds to the term involving $\lambda_{\cE}^{m+1,i} \l<\kappa^{(j)}, \eta^h \r>_{\Gamma^m}^h$ for $j=0,1,2$. However, this contribution is typically small, since $\lambda_{\cE}^{m+1,i}$ remains close to zero in practice, which was also confirmed by our numerical experiments. All other matrix components remain constant throughout the Newton iterations.
As a result, the structure of the system matrix is nearly identical to that of the linear BGN scheme \eqref{BGN}. Therefore, this similarity allows to reuse the preconditioned iterative solution methods discussed in \cite[Section 5]{BGN_2015_navierstokes}.

In addition, we can observe that the matrix $(A_{k \ell})_{k \ell}$ in \eqref{Reconstruct} becomes singular when $\kappa^{m+1,i}$ is nearly constant. In this case, equations \eqref{Decouple2} and \eqref{Decouple3} become almost identical, which results in linearly dependent columns in $(A_{k \ell})_{k \ell}$. In fact, the only difference between \eqref{Decouple2} and \eqref{Decouple3} lies in the terms $\l<\kappa^{m+1,i},\eta^h \r>_{\Gamma^m}^h$ and $\l<1,\eta^h \r>_{\Gamma^m}^h$, respectively. However, this degeneracy appears to be a natural consequence of the Lagrange multiplier approach, since the equivalence results in Theorems~\ref{Thm:relation:continuous} and~\ref{Thm:relation:semi} similarly rely on the assumption of a non-constant curvature. At this point, we refer to Section~\ref{sec:Single} for an approach using a single Lagrange multiplier.

\smallskip

\subsection{Single Lagrange multiplier approaches}\label{sec:Single}

As observed in Section~\ref{sec:Newton}, the formulation involving two Lagrange multipliers can lead to a singular linear system within the Newton iterations when the curvature is constant, e.g., for a stationary sphere. We also refer to the discussions in \cite[Section 4 and Remark 3.7]{Garcke2024} for related observations. This motivates to consider an alternative formulation involving only a single Lagrange multiplier. More precisely, we present two possible variants.
\begin{itemize}
    \item The energy-stable model variant: by setting $\lambda_{V}=0$ in the jump condition \eqref{Lag:jump3} and omitting the volume identity \eqref{Lag:Enforce2}.  
    \item The volume-preserving model variant: by setting $\lambda_{\cE}=0$ in the jump condition \eqref{Lag:jump3} and omitting the energy identity \eqref{Lag:Enforce1}.
\end{itemize}
Hereafter, we only discuss the Euler scheme \eqref{Lag:Euler} applied to the energy-stable model variant. 
However, this idea easily translates to the other numerical approaches in this work, i.e., the semi-discrete approximation \eqref{Lag:Semi}, the Crank--Nicolson scheme \eqref{Lag:CN}, and the BDF2 scheme \eqref{Lag:BDF2}. Similarly, the beforementioned schemes can also be applied to the model variant with enforced volume preservation.

Now the energy-stable Euler scheme reads as follows.
Given $\Gamma^m$, we aim to find a solution
\[
(\bu^{m+1},p^{m+1},\bX^{m+1},\kappa^{m+1})\in \U^m\times \widehat{\P}^m\times [W^m]^d\times W^m
\]
and a Lagrange multiplier $\lambda_{\cE}^{m+1}\in \R$, such that 
\begin{subequations}\label{Lag:Euler:single}
\begin{align}
    0 &= \l(2\mu^{m}\bD(\bu^{m+1}),\bD(\bm{\varphi}^h) \r) -\l(p^{m+1},\nabla\cdot \bm{\varphi}^h\r)-\gamma\l<\kappa^{m+1}\bnu^m,\bm{\varphi}^h\r>_{\Gamma^m} -\l(\widetilde{\bg}^{m+1},\bm{\varphi}^h \r),
    \\
    0 &= \l(\nabla\cdot \bu^{m+1},q^h \r),
    \\
    0 &= \l<\l(\frac{\bX^{m+1}-\bX^m}{\Delta t} \r)\cdot \bnu^m,\eta^h\r>_{\Gamma^m}^h -\l<\bu^{m+1}\cdot \bnu^m,\eta^h\r>_{\Gamma^m}  
    - \lambda_{\cE}^{m+1}\l<\kappa^{m+1},\eta^h \r>_{\Gamma^m}^h,
    \\
    0 &= \l<\kappa^{m+1}\bnu^m,\bm{\omega}^h\r>_{\Gamma^m}^h+\l<\nabla_{\Gamma^m}\bX^{m+1},\nabla_{\Gamma^m}\bm{\omega}^h\r>_{\Gamma^m},\qquad\quad\qquad\quad\quad\quad
    \\
    0&=\frac{\cE^{m+1}-\cE^{m}}{\Delta t}+\l(2\mu^m\bD(\bu^{m+1}),\bD(\bu^{m+1}) \r)-\l(\bu^{m+1}, \widetilde{\bg}^{m+1}\r), \label{Lag:Euler5:single}
\end{align}
\end{subequations}
for any $(\bm{\varphi}^h,q^h,\eta^h,\bm{\omega}^h)\in \U^m\times \widehat{\P}^m\times W^m\times [W^m]^d$.

\smallskip

The scheme is unconditionally energy-stable by construction. In addition to that, it has the desirable property of eliminating spurious velocities in the case of steady-state solutions, similar to the approach in \cite{BGN_2013_stokes}.

\begin{lemma}
Let $(\bu^{m+1},p^{m+1},\bX^{m+1},\kappa^{m+1})$ and $\lambda^{m+1}_{\cE}$ be a solution of \eqref{Lag:Euler:single} with $\widetilde{\bg}^{m+1}=0$. If $\Gamma^{m+1}=\Gamma^m$, then $\bu^{m+1}=0$.
\end{lemma}

\Proof
It suffices to notice that $\Gamma^{m+1}=\Gamma^m$, $\widetilde{\bg}^{m+1}=0$ and \eqref{Lag:Euler5:single} lead to 
\[
\l(2\mu^m\bD(\bu^{m+1}),\bD(\bu^{m+1}) \r)=0,
\]
where $\mu^m(\cdot) \geq \min(\mu_+, \mu_-) > 0$.
Thus, by employing Korn's inequality we obtain $\bu^{m+1}=0$, as desired.
\qed 

We note that a corresponding result holds for the structure-preserving variant \eqref{Lag:Euler} involving two Lagrange multipliers, as the proof is based on the enforced energy identity.

\section{Numerical results}
In the following, we present three numerical scenarios including convergence tests to illustrate the practicability of our proposed numerical method. 
We restrict to the two-dimensional case $d=2$. Also, we always fix the function spaces for the discrete velocity and pressure as $\mathbb{U}^m = [S_2]^2 \cap \mathbb{U}$ and $\widehat{\mathbb{P}}^m = (S_1 + S_0) \cap \widehat{\mathbb{P}}$.
The numerical tests have been performed using the finite element library FEniCS \cite{fenics_book_2012}, which provides access to the linear algebra backend PETSc \cite{petsc-user-ref_2021}. The implementation builds upon the code developed in \cite{GNT_2024}. In most cases, we focus on the fitted approach, and similar experiments can be conducted using an unfitted bulk mesh.

\subsection{Expanding bubble}

To test our method, we consider a scenario similar to \cite{agnese_nurnberg_2016_stokes} where a nontrivial divergence-free and radially symmetric solution $\bu$ to \eqref{Lag:weak} can be constructed on a domain that does not contain the origin. In this case, we will set 
\[\Omega = (-1,1)^2 \setminus [-\tfrac15, \tfrac15]^2 ,\] 
and let $\alpha\geq0$ be a given constant. 
Let $\Gamma(t) = \{\bx \in \R^2\mid |\bx| = r(t)\}$ be a sphere of radius $r(t)$ and curvature $\kappa(t) = -\frac{1}{r(t)}$. Then, the expanding sphere with
\begin{gather*}
    r(t) = \l( r(0)^2 + 2\alpha t  \r)^{1/2},
    \quad
    \bu(\bx,t) = \alpha \frac{\bx}{|\bx|^2},
    \\
    p(\bx,t) = - \l( \gamma + 2 \alpha \frac{\mu_+ - \mu_-}{r(t)} \r) \cdot \kappa(t) \left( \bm{\chi}_{\Omega_-(t)} - \frac{\mathcal{L}^2(\Omega_-(t))}{\mathcal{L}^2(\Omega)} \right),
\end{gather*}
is an exact solution with forcing term $\bg(\bx) = \alpha \frac{\bx}{|\bx|^2}$ and non-homogeneous boundary data $\bu_D(\bx) = \alpha \frac{\bx}{|\bx|^2}$ on $\partial\Omega$. 

However the presence of non-homogeneous Dirichlet boundary conditions $\bu = \bu_D$ on $\partial\Omega$ leads to modified energy and volume identities compared to those in \eqref{Lage:Enforce}. As a result, the numerical schemes must be adapted accordingly.
In this setting, the corresponding energy and volume identities take the form:
\begin{align}
     \label{Lage:energy_expanding_sphere}
    0&=\frac{\d }{\d t}\cE(t)
    + \l(2\mu\bD(\bu),\bD(\bu-\tilde{\bu}) \r) 
    - \l(\bg, \bu-\tilde{\bu}\r)
    + \l( p, \nabla\cdot \tilde{\bu} \r)
    + \gamma \int_{\Gamma(t)} \kappa \bu_D \cdot \bnu \, \d\mathcal{H}^{1}, 
    \\
    \label{Lage:volume_expanding_sphere}
    0&=\frac{\d }{\d t}V(t) + \int_{\partial_2\Omega} \bu_D \cdot \bn \, \d\mathcal{H}^{1},
\end{align}
where $\tilde{\bu}$ is an arbitrary and smooth extension of $\bu_D$ to the interior of the domain, and $\partial_2\Omega \colonequals \partial\Omega \cap \partial\Omega_-(t)$ is non-empty since $\Omega = (-1,1)^2 \setminus [-\frac15, \frac15]^2$ contains a hole in this setting.
We note that by fixing $\tilde{\bu} = \bu$ the energy identity \eqref{Lage:energy_expanding_sphere} simplifies to
\begin{align}
    \label{Lage:energy2_expanding_sphere}
    0&=\frac{\d }{\d t}\cE(t)
    + \gamma \int_{\Gamma(t)} \kappa \bu_D \cdot \bnu \, \d\mathcal{H}^{1}.
\end{align}
We remark that the schemes \eqref{Lag:Euler}, \eqref{Lag:CN}, \eqref{Lag:BDF2}, as well as their variants involving a single Lagrange multiplier, and the BGN scheme \eqref{BGN} can be adapted to the scenario of an expanding bubble with only minor changes. For this reason, these modifications are not explicitly shown here.

Since the curvature of the expanding sphere $\Gamma(t)$ remains spatially constant at all times $t\geq 0$, the resulting system matrix in the Newton iterations may become singular, see Section~\ref{sec:Newton}. To avoid this issue, we fix the Lagrange multiplier $\lambda_{\cE}(t)= 0$ and omit the energy equation \eqref{Lage:energy_expanding_sphere}. As explained in Section~\ref{sec:Single}, this results in a single Lagrange multiplier formulation for the schemes \eqref{Lag:Euler}, \eqref{Lag:CN}, and \eqref{Lag:BDF2}, where only the volume growth identity \eqref{Lage:volume_expanding_sphere} is enforced.

For the expanding bubble experiment, the model parameters are fixed as 
\begin{align*}
    T = 0.1, \quad
    \mu_+ = 5, \quad
    \mu_- = 2, \quad
    \gamma = 0.1, \quad
    \alpha = 0.2 \, .
\end{align*}
Moreover, we consider the following values for the the time step size and the number of interface points:
\begin{align*}
    \Delta t \in \{T/4, T/8, \ldots, T/2048\}, \qquad 
    K_{\Gamma^h} \in \{16, 20, 32, 40, 64, 80, 120, 160\} \, .
\end{align*}
In this test case, we employ a fitted bulk mesh approach, following the strategy outlined in \cite{agnese_nurnberg_2016_stokes}. Accordingly, the spatial mesh size $h$ is chosen proportional to $1/K_{\Gamma^h}$, and for simplicity, we write $h = 1/K_{\Gamma^h}$.

The initial parametric surface $\Gamma^0$ consists of $K_{\Gamma^h}$ vertices and is defined as a polygonal approximation of a circle with center $(0,0)^\top$ and radius $r(0)=0.5$, respectively.

Starting from the initial polyhedral surface $\Gamma^0$, we generate a bulk triangulation $\mathcal{T}^0$ of the domain $\Omega$, that is fitted to $\Gamma^0$, using the open-source mesh generator Gmsh \cite{gmsh_library}. 
For each time step $m \geq 0$, after computing the updated interface $\Gamma^{m+1}$, the bulk mesh $\mathcal{T}^{m+1}$ must be adapted to fit $\Gamma^{m+1}$, ideally close to $\mathcal{T}^m$. To reduce the computational overhead associated with remeshing the full domain, we apply a mesh smoothing procedure, where a suitable displacement field is computed by solving an elliptic elasticity problem, which is then used to deform the mesh. However, if the resulting deformation becomes too large and elements get stretched too much, a full remeshing of the domain $\Omega$ becomes necessary. The precise remeshing and smoothing strategy is detailed in \cite[Sec.~3.7]{agnese_nurnberg_2016_stokes}.

\begin{figure}[h!]
\centering
\includegraphics[width=0.32\textwidth,trim={2cm 3.8cm 2cm 3.8cm},clip]{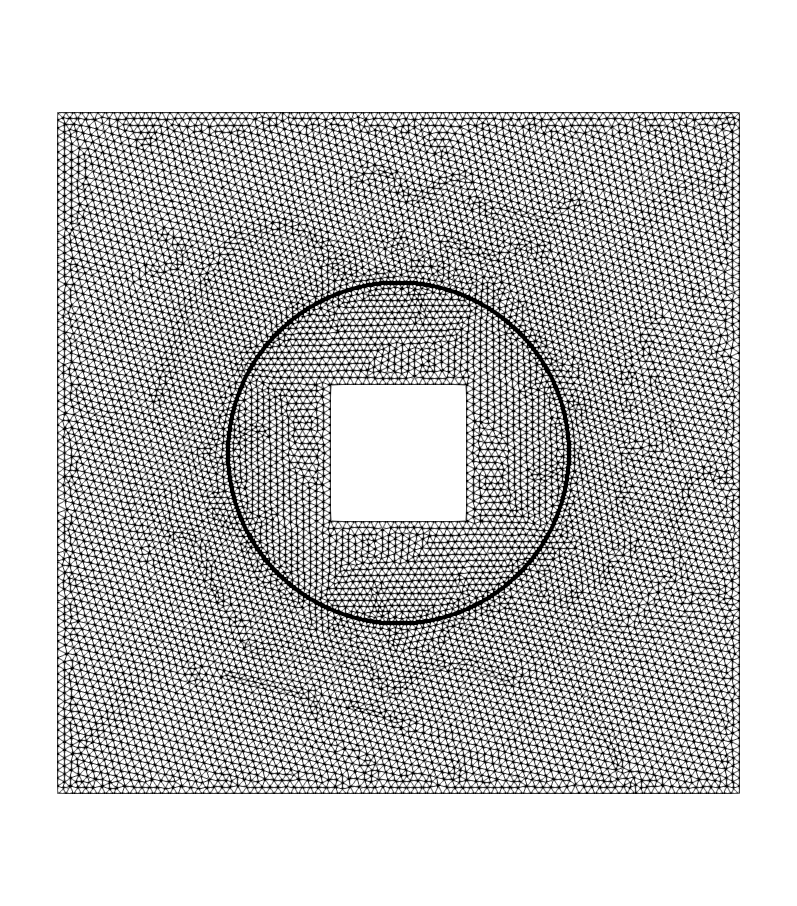} \,
\includegraphics[width=0.32\textwidth,trim={2cm 3.8cm 2cm 3.8cm},clip]{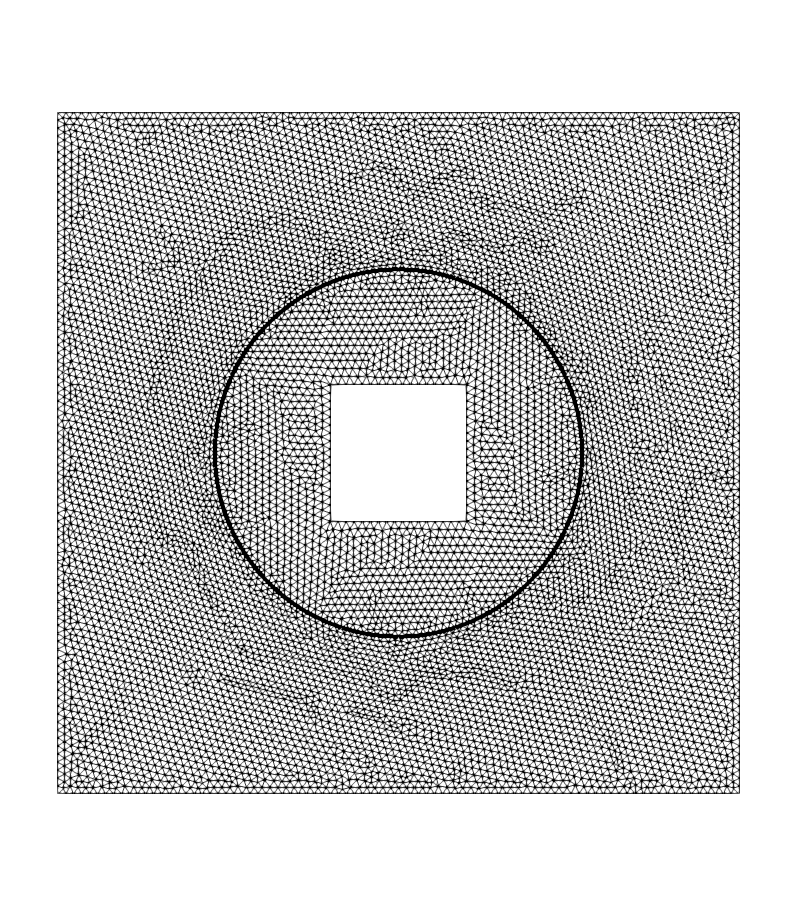}
\caption{
The surface and bulk meshes for the expanding bubble experiment at times $t\in\{0, 0.1\}$, using the volume-preserving Euler scheme with $K_{\Gamma^h} = 160$ and $\Delta t = T/2048$, where $T=0.1$.}
\label{fig:expanding2_bubbles}
\end{figure}

\begin{figure}[h!]
\centering
\includegraphics[width=0.28\textwidth]{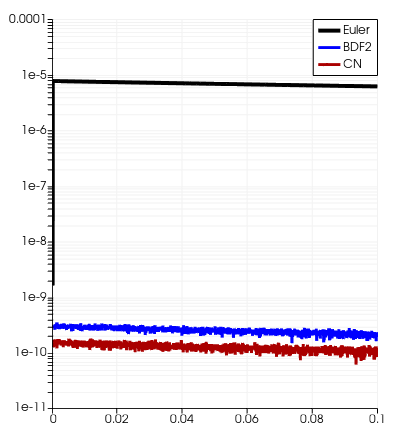} 
\,
\includegraphics[width=0.28\textwidth]{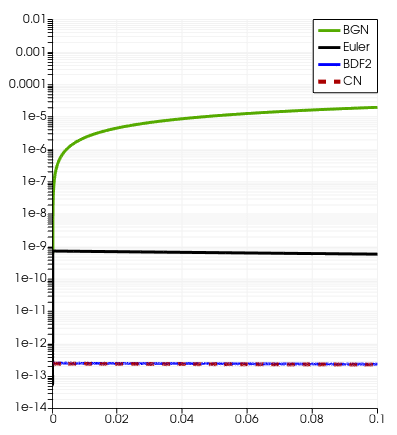}
\, 
\includegraphics[width=0.28\textwidth]{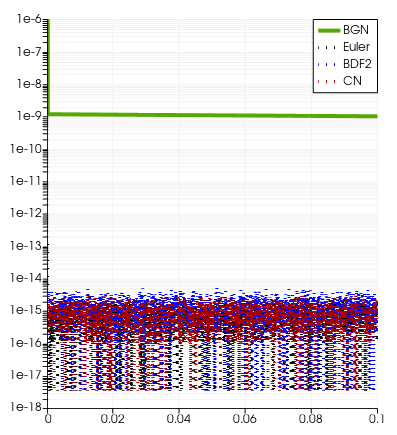}
\caption{Numerical results for the expanding bubble over the time interval $[0, 0.1]$, computed using the BGN scheme and the volume-preserving Euler, BDF2, and Crank--Nicolson (CN) schemes with $K_{\Gamma^h} = 160$ and $\Delta t = T/2048$, where $T=0.1$.
From left to right: time evolution of the magnitude of the Lagrange multiplier $\lambda_V$, and the numerical residuals of the energy identity \eqref{Lage:energy2_expanding_sphere} and volume identity \eqref{Lage:volume_expanding_sphere} adapted to the respective schemes \eqref{Lag:Euler}, \eqref{Lag:CN} and \eqref{Lag:BDF2}.}
\label{fig:expanding2_Lagr}
\end{figure}

In Figure~\ref{fig:expanding2_bubbles}, we display the interface and the fitted bulk mesh at times $t \in \{0, 0.1\}$ for the expanding bubble experiment, showing the growth of the bubble over time. The simulation was performed using the volume-preserving Euler scheme with $K_{\Gamma^h} = 160$ and $\Delta t = T/2048$. Since the meshes produced by the BGN, BDF2 and Crank--Nicolson schemes are visually indistinguishable, they are not shown here. 

Figure~\ref{fig:expanding2_Lagr} presents the time evolution of the Lagrange multiplier $\lambda_V$ for the volume equation \eqref{Lage:volume_expanding_sphere}, along with the numerical residuals of the energy identity \eqref{Lage:energy2_expanding_sphere} and the enforced volume identity \eqref{Lage:volume_expanding_sphere}, for all schemes. Note that the Lagrange multiplier is not shown for the BGN scheme. Also, the term numerical residuum refers to the residual of the respective discrete identities of \eqref{Lage:volume_expanding_sphere} and \eqref{Lage:energy2_expanding_sphere}, as adapted to the BGN, Euler, Crank--Nicolson, and BDF2 schemes in \eqref{BGN}, \eqref{Lag:Euler}, \eqref{Lag:CN}, and \eqref{Lag:BDF2}, respectively. 
For instance, the results show that the Lagrange multiplier $\lambda_V$ is of magnitude $10^{-5}$ for the Euler scheme, and around $10^{-10}$ for the second-order schemes, with the Crank--Nicolson method producing slightly smaller values. The volume identity is preserved up to machine accuracy across all volume-preserving schemes, while the residuum is at magnitude $10^{-9}$ for the BGN scheme. In contrast, the energy residual is of order $10^{-5}$ for the BGN scheme and of order $10^{-9}$ for the Euler scheme, while it is reduced to approximately $10^{-12}$ for the BDF2 and Crank--Nicolson schemes. The reason for this is that the energy identity is not enforced for the volume-preserving schemes based on Section \ref{sec:Single}.

Subsequently, we study the numerical convergence for the Euler, Crank--Nicolson and BDF2 schemes, and compare it to the results from the BGN scheme \eqref{BGN} based on \cite{agnese_nurnberg_2016_stokes}.
We consider the following error quantities for the parametrization $\bX$, the velocity $\bu$, the pressure $p$, and the Lagrange multiplier $\lambda_{V}$:
\begin{align}
\label{def:errors}
    \nonumber
    e_{\bX} &= \max_{m\in\{1,\ldots,M\}} d(\Gamma_{ref}(t_m), \Gamma^m),
    \qquad
    e_{\bu} = \left( \Delta t \sum_{m=1}^{M} 
    \|I_2^h \bu(\cdot,t_m) - \bu^m\|_{L^2}^2 \right)^{1/2},
    \\ 
    e_{p} &= \left( \Delta t \sum_{m=1}^{M} 
    \|p(\cdot,t_m) - p^m\|_{L^2}^2 \right)^{1/2},
    \qquad
    e_{\lambda, V} = \left( \Delta t \sum_{m=1}^{M} |\lambda_{V}^m|^2 \right)^{1/2},
\end{align}
where 
\begin{align*}
    d(\Gamma_1,\Gamma_2) = |(\Omega_1 \setminus \Omega_2) \cup (\Omega_2 \setminus \Omega_1)| = |\Omega_1| + |\Omega_2| - 2 |\Omega_1 \cap \Omega_2|
\end{align*}
is the manifold distance \cite{Jiang2024a,Zhao2021} between two surfaces $\Gamma_i$, $i\in\{1,2\}$, with enclosed volumes $\Omega_i$, $i\in\{1,2\}$.
For computational reasons, we approximate the exact interface $\Gamma(t)$ using a highly resolved polygonal curve $\Gamma_{ref}$ with $K_{ref} = 640$ vertices. This reference curve satisfies the condition $|\bX_{ref}(\bq_k, t)| = r(t)$ for all vertices $\bq_k$, $k = 1, \ldots, K_{ref}$, where $r(t)$ denotes the radius of the exact interface $\Gamma(t)$ for the expanding bubble problem.

To investigate the convergence behaviour more systematically, we study spatial and temporal convergence separately for each numerical scheme. Specifically, we plot the errors in the relevant quantities against the spatial mesh size $h$ and the time step size $\Delta t$ in two distinct sets of plots. In these plots, lines with different symbols correspond to varying time step sizes (for spatial convergence) and mesh sizes (for temporal convergence), respectively.

\begin{figure}[h!]
\centering
\includegraphics[width=0.22\textwidth]{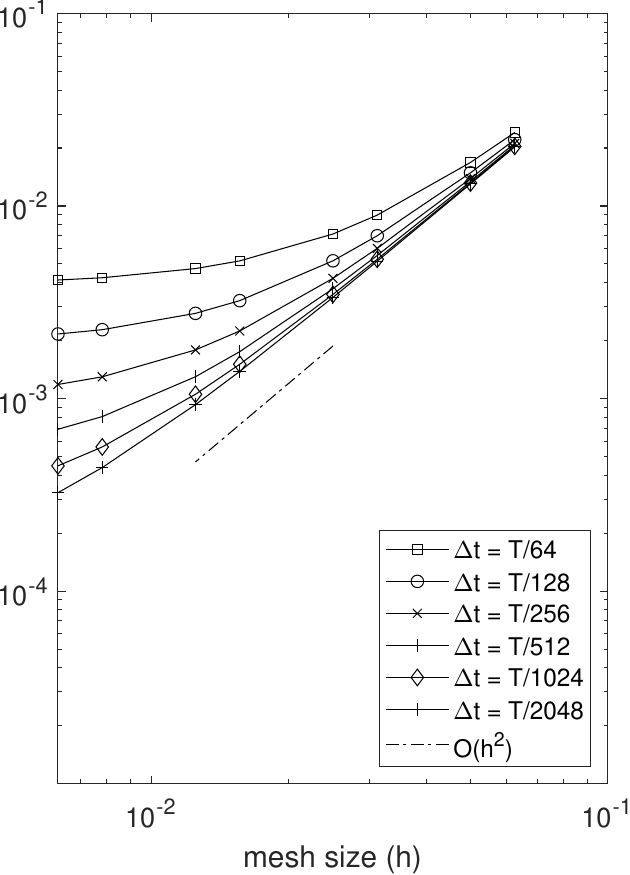}
\includegraphics[width=0.22\textwidth]{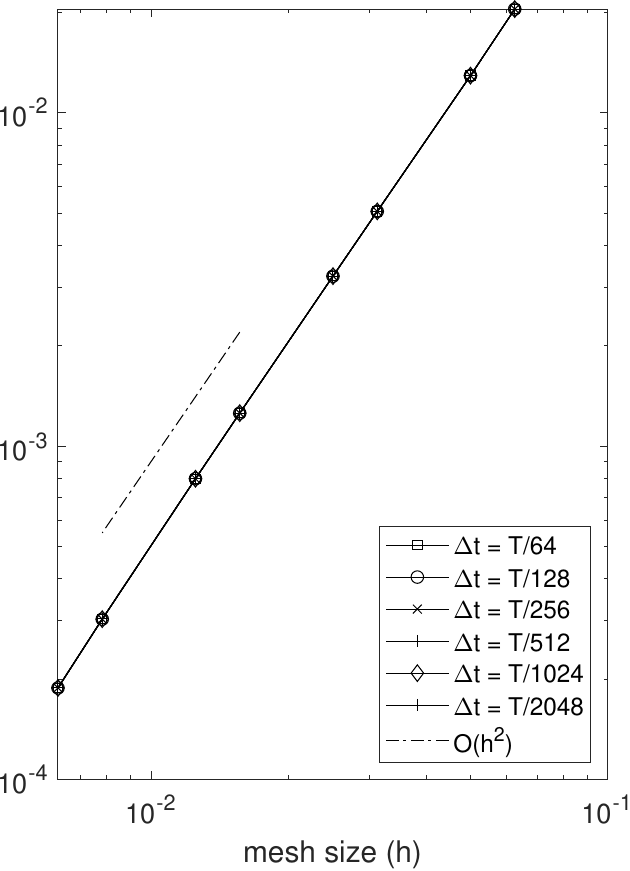}
\includegraphics[width=0.22\textwidth]{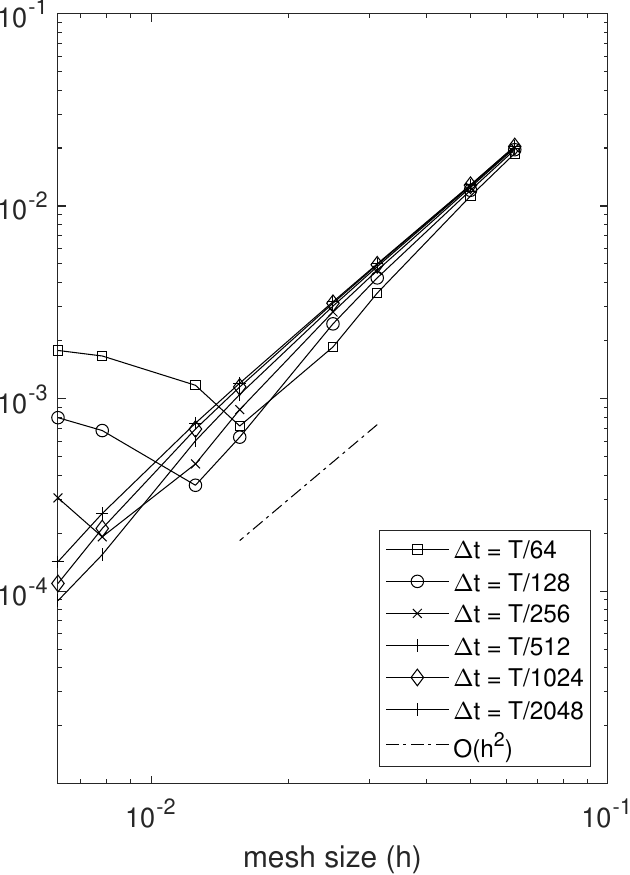}
\includegraphics[width=0.22\textwidth]{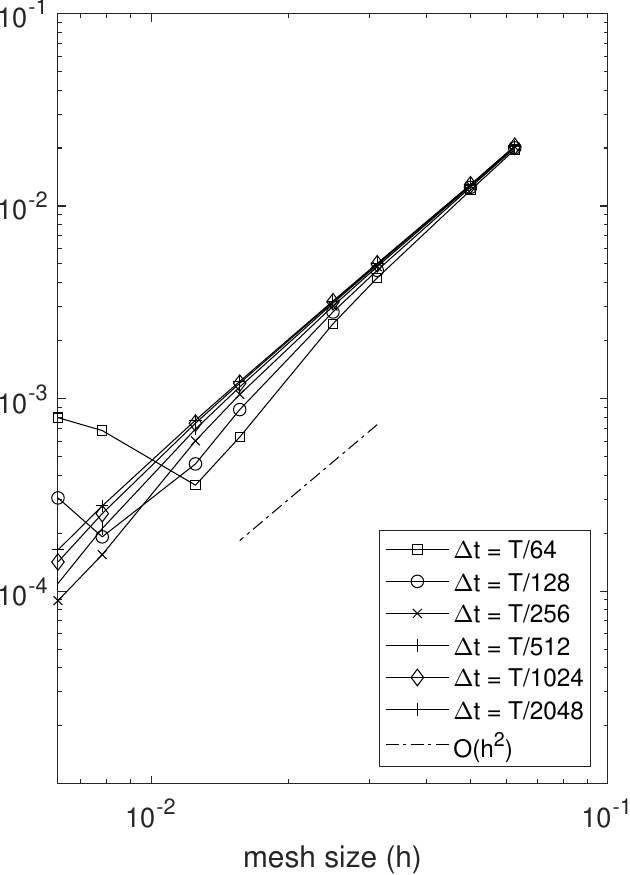}
\caption{Studying the spatial convergence of the error quantity $e_{\bX}$ for the parametrization $\bX$ in the expanding bubble experiment. From left to right: BGN scheme, Euler scheme, BDF2 scheme, Crank--Nicolson scheme.}
\label{fig:3}      
\end{figure}

\begin{figure}[h!]
\centering
\includegraphics[width=0.22\textwidth]{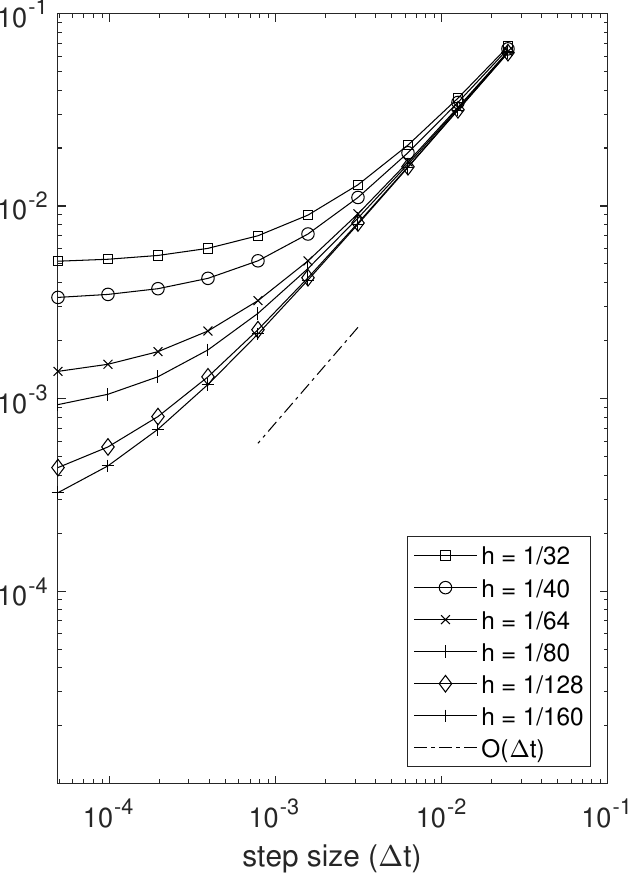}
\includegraphics[width=0.22\textwidth]{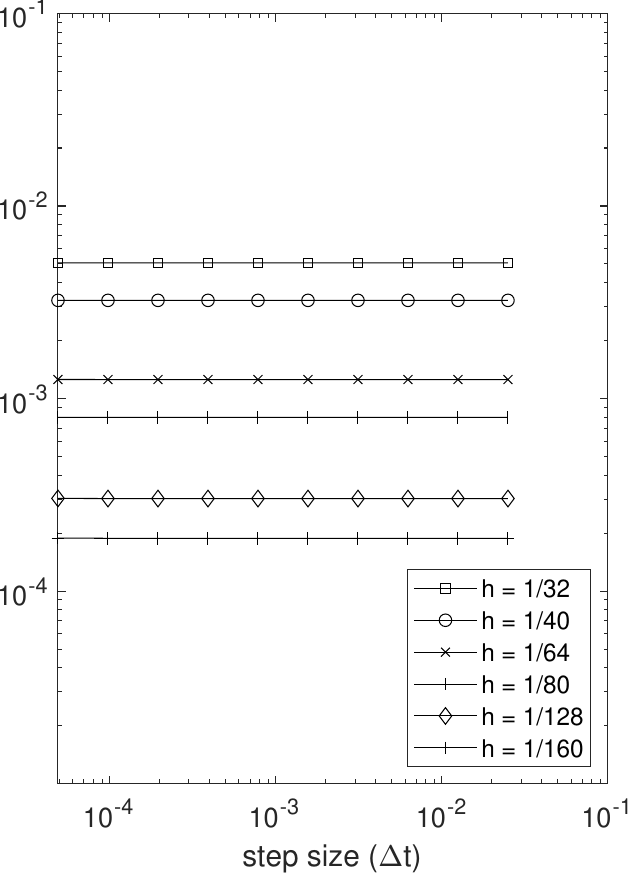}
\includegraphics[width=0.22\textwidth]{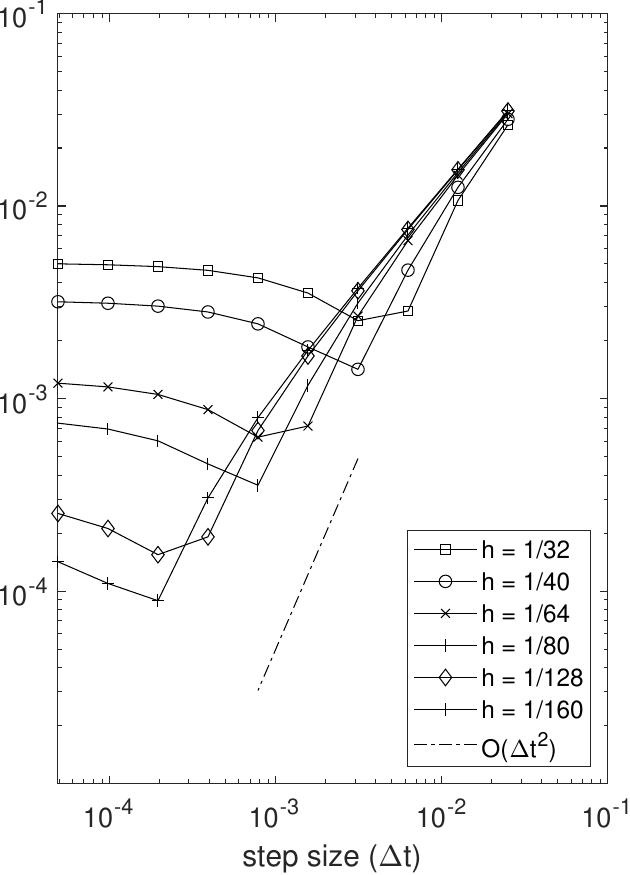}
\includegraphics[width=0.22\textwidth]{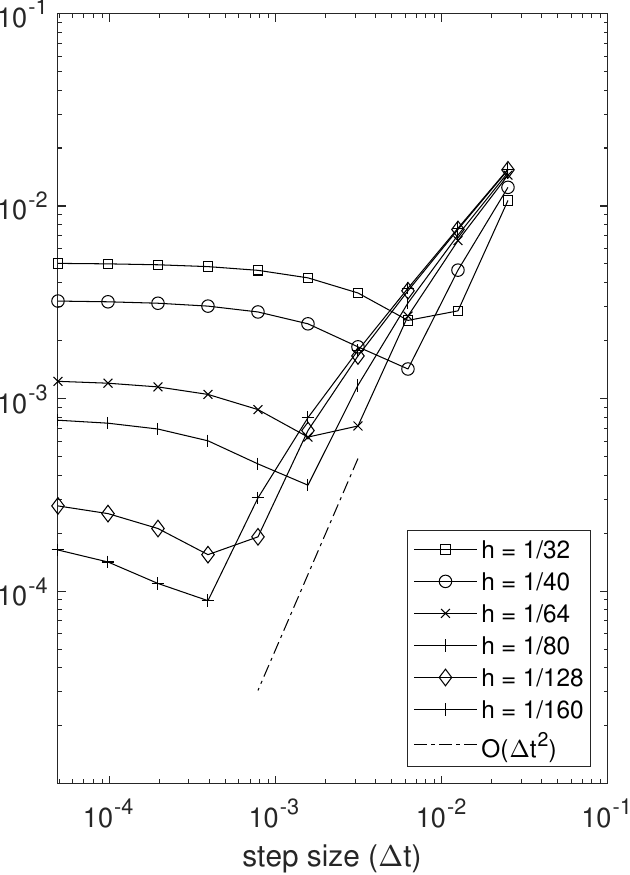}
\caption{Studying the temporal convergence of the error quantity $e_{\bX}$ for the parametrization $\bX$ in the expanding bubble experiment. From left to right: BGN scheme, Euler scheme, BDF2 scheme, Crank--Nicolson scheme.}
\label{fig:4}      
\end{figure}

Figures~\ref{fig:3} and~\ref{fig:4} show log-log plots of the interface error $e_{\bX}$ for the BGN scheme, the Euler scheme, the BDF2 scheme, and the Crank--Nicolson scheme with respect to spatial and temporal refinement, respectively. All four schemes demostrate second-order convergence in space, and the expected first- and second-order convergence in time, respectively. 
It is worth noting that the expected convergence rates become clearly visible only when the time step size is relatively large (e.g., $\Delta t \in \{T/64, \ldots, T/512\}$) in the spatial convergence tests, and when the mesh size is relatively coarse (e.g., $h \in \{1/32, \ldots, 1/128\}$) in the temporal convergence tests. For finer resolutions, the error curves tend to flatten out, which indicates that the discretization error from the other parameter becomes dominant. 
We also observed similar convergence behaviour for the error in the energy $\cE^m$ in our experiments, although these results are not shown here.
Among all schemes, the Euler method performs particularly well for larger time steps: the error remains nearly constant across a wide range of step sizes and is comparable in magnitude to the errors of the other schemes at smaller time steps. For very small time steps, the second-order methods yield the lowest errors, and we expect this trend to continue for even finer spatial and temporal discretizations.

\begin{figure}[h!]
\centering
\includegraphics[width=0.22\textwidth]{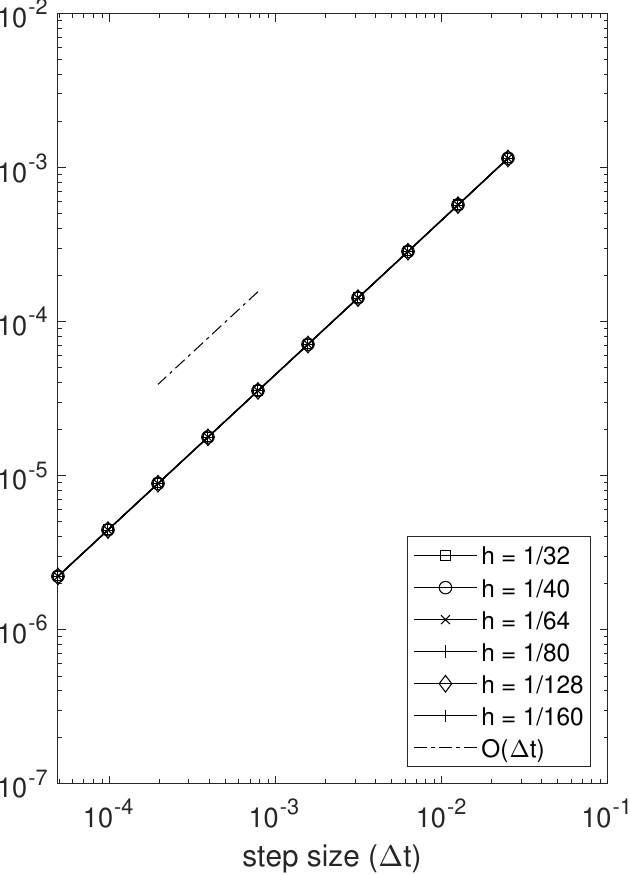}
\includegraphics[width=0.22\textwidth]{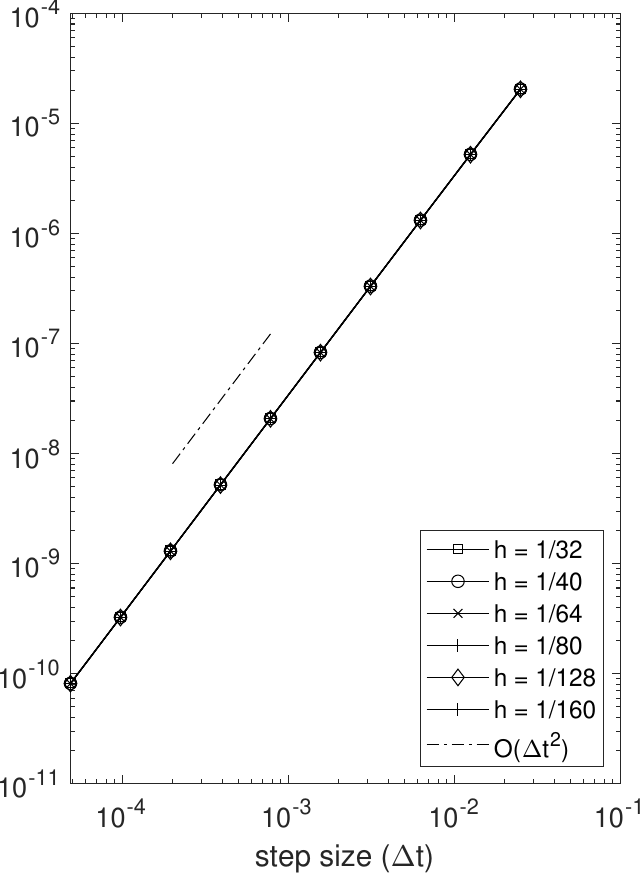}
\includegraphics[width=0.22\textwidth]{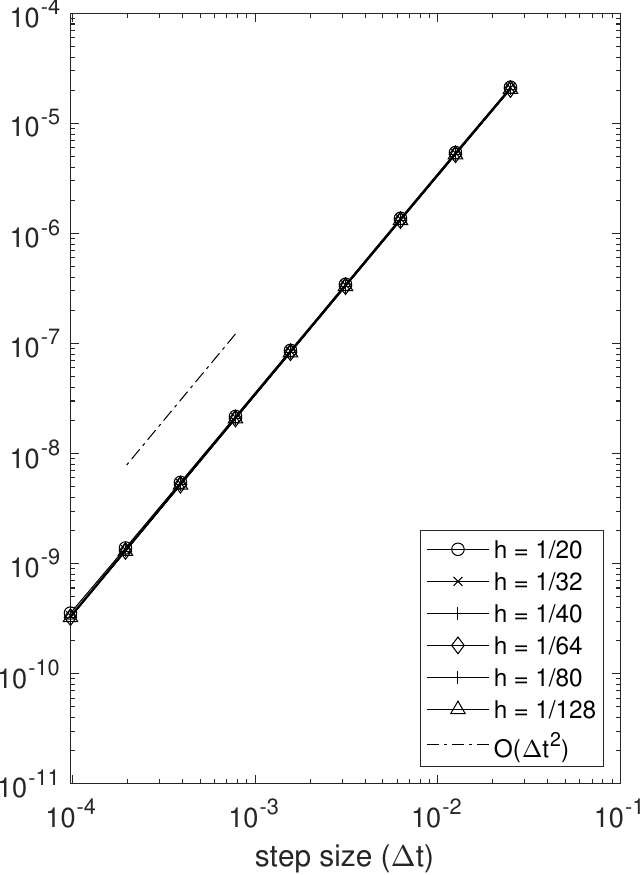}
\caption{Studying the temporal convergence of the error quantity $e_{\lambda,V}$ for the Lagrange multiplier $\lambda_{V}$ in the expanding bubble experiment. From left to right: Euler scheme, BDF2 scheme, Crank--Nicolson scheme. 
}
\label{fig:6}      
\end{figure}

In Figure~\ref{fig:6}, we present the temporal errors $e_{\lambda,V}$ in the Lagrange multiplicator $\lambda_{V}$ associated with the enforced volume identity, where the expected first- and second-order convergence is confirmed for the respective schemes.

\begin{figure}[h!]
\centering
\includegraphics[width=0.22\textwidth]{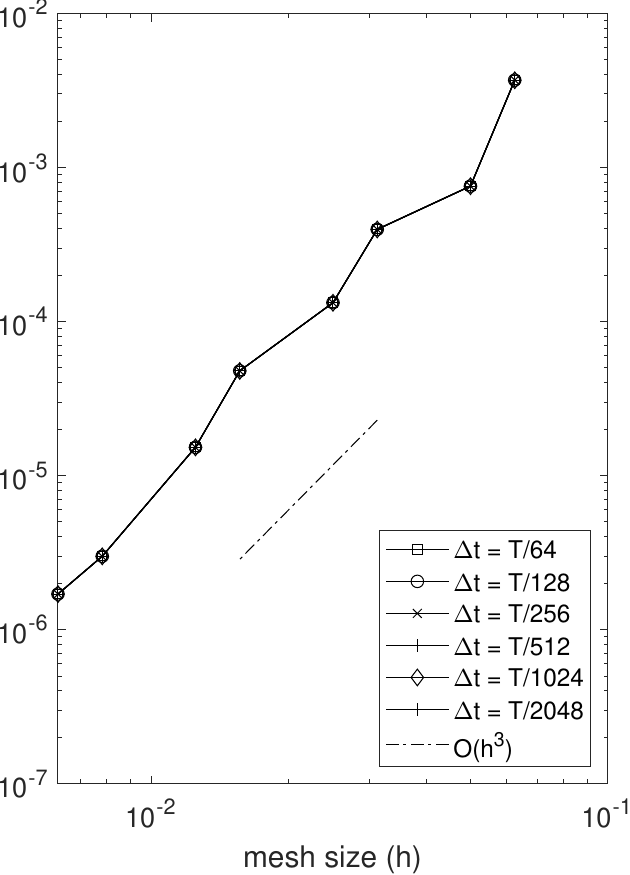}
\includegraphics[width=0.22\textwidth]{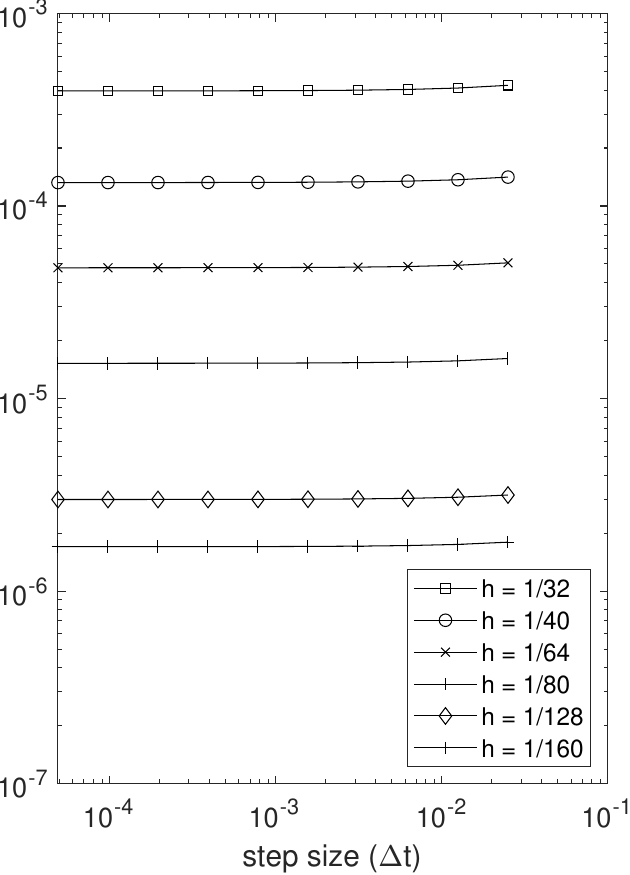}
\includegraphics[width=0.22\textwidth]{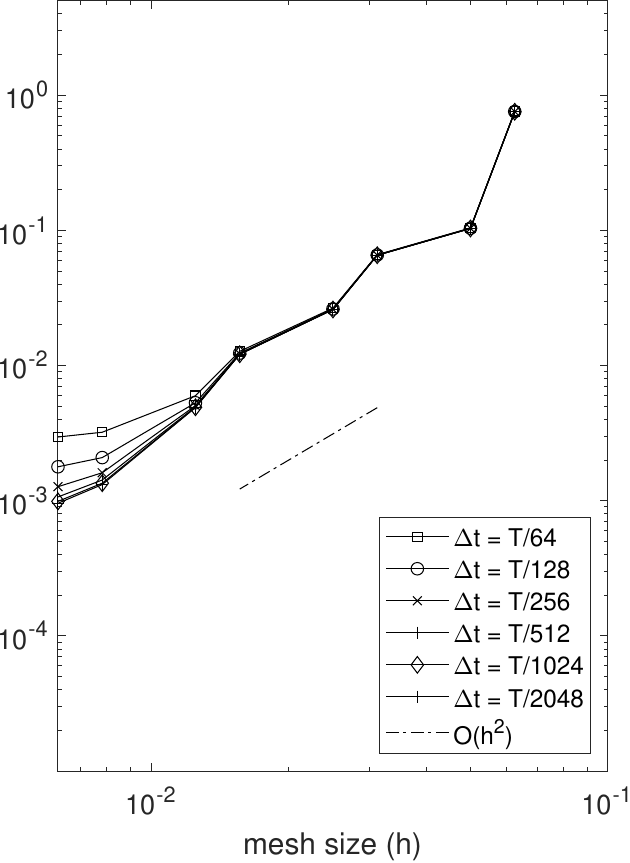}
\includegraphics[width=0.22\textwidth]{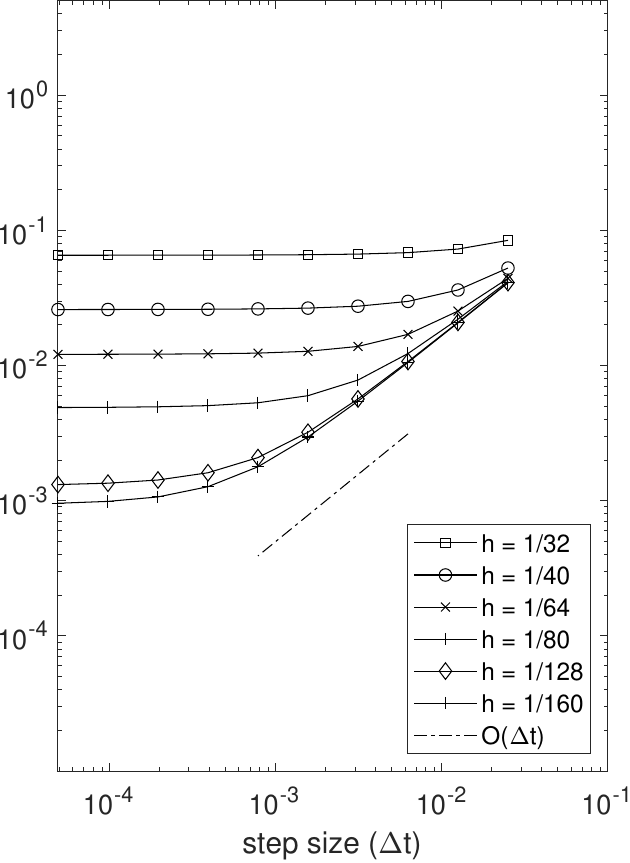}
\caption{From left to right: the spatial and temporal errors $e_{\bu}$ for the velocity $\bu$, and the spatial and temporal errors $e_{p}$ for the pressure $p$ in the expanding bubble experiment, computed with the Euler scheme. The error plots for the other schemes (BGN, BDF2, Crank--Nicolson) are almost identical, and therefore not shown.}
\label{fig:7}      
\end{figure}

Finally, Figure~\ref{fig:7} shows the velocity and pressure errors, $e_{\bu}$ and $e_p$. For simplicity, we report only the results for the Euler scheme, as the BGN, BDF2, and Crank--Nicolson schemes exhibit similar behaviour and are omitted for brevity. In the temporal refinement plots, the error remains nearly constant for all time steps, indicating that the spatial discretization error dominates. In the spatial refinement plots, we observe third-order convergence for the velocity and second-order convergence for the pressure.

\subsection{Retraction of an ellipse}
In this example, we consider an initial interface that approximates an ellipse. Due to surface tension forces, the interface would eventually retract into a circle. However, we terminate the simulation early, while the curvature is still non-uniform to avoid the issue of a singular system of equations within the Newton iterations, see Section~\ref{sec:Newton}.

As in the previous section, we now consider the fitted approach,  where the bulk mesh is aligned with the evolving interface using the strategies described earlier.

For this numerical experiment, we fix the parameters as follows: 
\begin{align*}
    \Omega=(0,1)^2, \quad T = 0.1, 
    \quad \mu_\pm = 2 \quad \gamma=5.
\end{align*}
The discretization parameters are chosen as:
\begin{align*}
    \Delta t \in \{T/4, T/8, \ldots, T/1024\}, \qquad
    K_{\Gamma^h} \in \{16,20,32,40,64,80,120\} .
\end{align*}

The initial parametric surface $\Gamma^0$ consists of $K_{\Gamma^h}$ vertices and is defined as a polygonal approximation of an ellipse with center $(0.5,0.5)^\top$ and horizontal and vertical axis of length $0.4$ and $0.1$, respectively.  
In particular, the vertices $\bq_k^0$, $k\in\{1,\ldots,K_{\Gamma^h}\}$, of $\Gamma^0$ are set to
\begin{align*}
    \bq_k^0 = 
    \begin{pmatrix}
        0.5 \\ 0.5
    \end{pmatrix}
    +
    \begin{pmatrix}
        0.4 \cos(\phi_k) \\ 0.1 \sin(\phi_k)
    \end{pmatrix},
    \quad 
    \phi_k = \frac{2\pi k}{K_{\Gamma^h}}, \quad k \in \{1,\ldots,K_{\Gamma^h}\}.
\end{align*}

\begin{figure}[h!]
\centering
\includegraphics[width=0.32\textwidth,trim={2cm 3.8cm 2cm 3.8cm},clip]{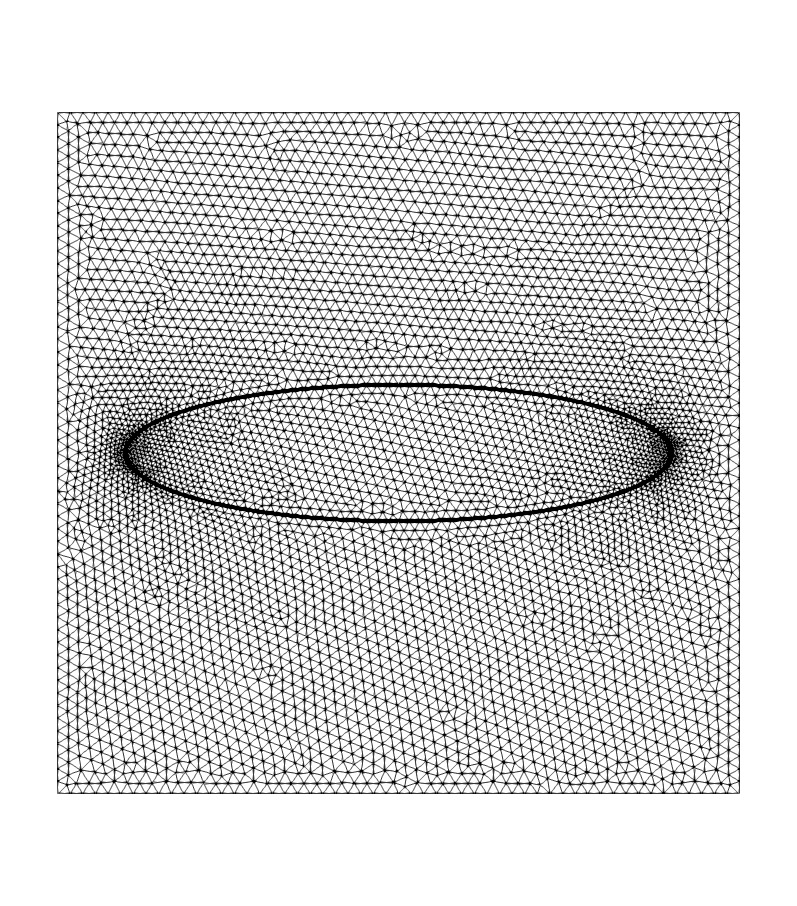} \,
\includegraphics[width=0.32\textwidth,trim={2cm 3.8cm 2cm 3.8cm},clip]{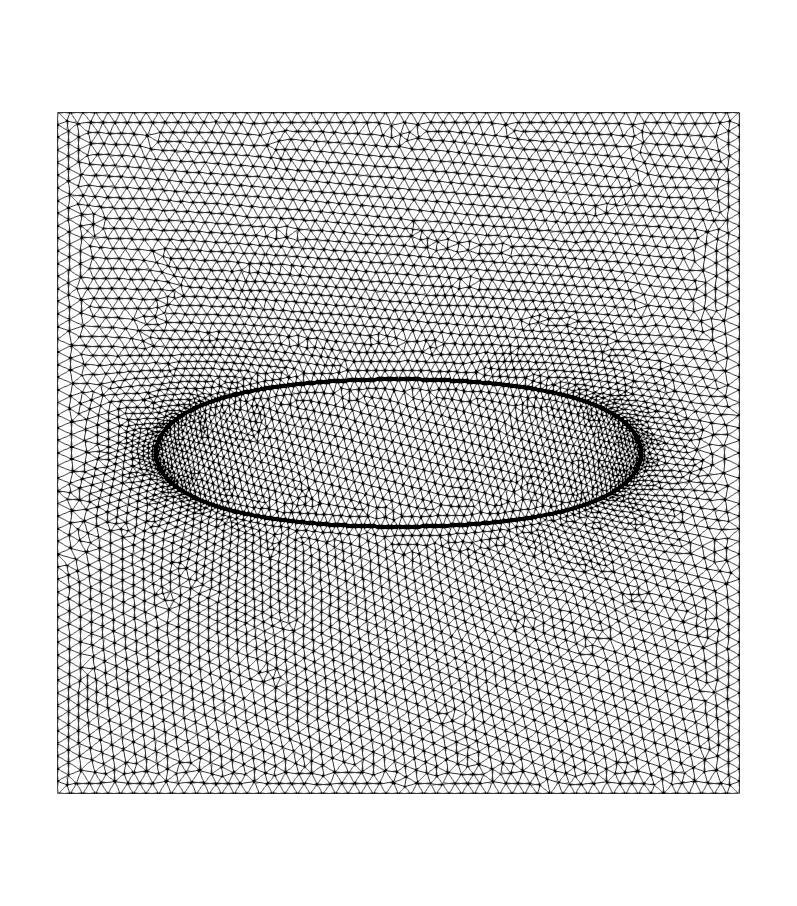}
\caption{
The interface and bulk meshes for the retracting ellipse at times $t\in\{0, 0.1\}$, using the structure-preserving Euler scheme with $K_{\Gamma^h} = 160$ and $\Delta t = T/2048$, where $T=0.1$.
}
\label{fig:retracting_bubbles}
\end{figure}

\begin{figure}[h!]
\centering
\includegraphics[width=0.23\textwidth]{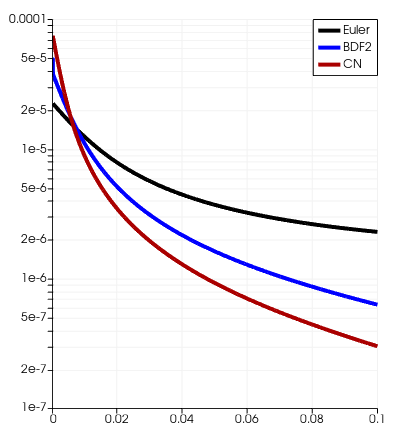} \,
\includegraphics[width=0.23\textwidth]{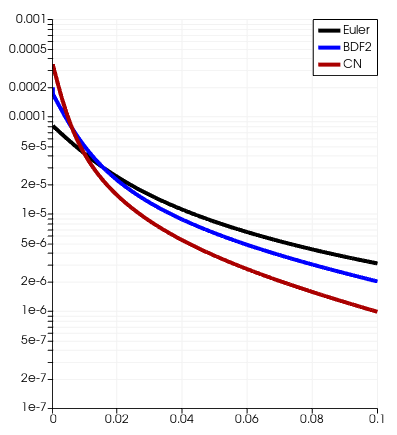} \,
\includegraphics[width=0.23\textwidth]{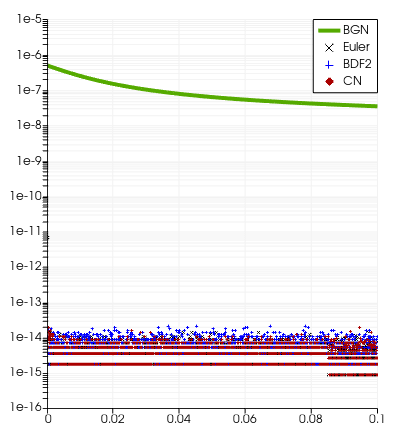} \, 
\includegraphics[width=0.23\textwidth]{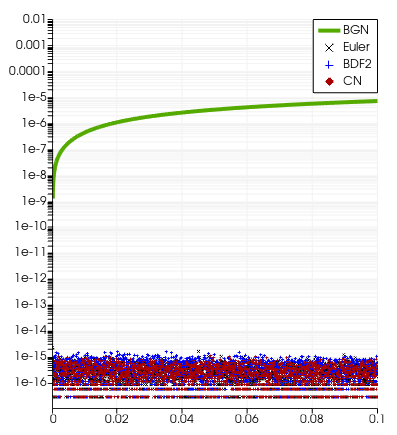}
\caption{
Numerical results for the retracting ellipse over the time interval $[0, 0.1]$, computed using the BGN scheme and the structure-preserving Euler, BDF2, and Crank--Nicolson (CN) schemes with $K_{\Gamma^h} = 160$ and $\Delta t = T/2048$, where $T=0.1$.
From left to right: time evolution of the magnitude of the Lagrange multipliers $\lambda_{\cE}$, $\lambda_V$, and the numerical residuals of the energy and volume identities for the respective scheme.
}
\label{fig:retracting_Lagr}
\end{figure}

Figure~\ref{fig:retracting_bubbles} shows the interface and the fitted bulk mesh at times $t \in \{0, 0.1\}$ for the retracting bubble experiment. The simulation was carried out using the structure-preserving Euler scheme with $K_{\Gamma^h} = 160$ and $\Delta t = T/2048$. Similar to before, since the meshes obtained from the BGN, BDF2, and Crank--Nicolson schemes are visually indistinguishable, they are not shown here.
In Figure~\ref{fig:retracting_Lagr}, we present the time evolution of the Lagrange multipliers $\lambda_{\cE}$ and $\lambda_V$ associated with the energy and volume identities, together with the numerical residuals of the corresponding identities, for all schemes. Note that the Lagrange multipliers are not included for the BGN scheme.
The results show that $\lambda_{\cE}$ and $\lambda_V$ decrease over time for all structure-preserving schemes, starting around $10^{-4}$ and dropping to values near $10^{-6}$, with slightly smaller values observed for the second-order methods. Both the energy and volume identities are preserved up to machine precision by the structure-preserving schemes. In contrast, the residuals for the BGN scheme from \cite{agnese_nurnberg_2016_stokes} remain at a range between $10^{-5}$ and $10^{-8}$ throughout the simulation.

As before, we study the error quantities defined in \eqref{def:errors}, but now computed with respect to the reference solution instead of an exact one. In addition, we introduce the error quantity $e_{\lambda,\cE}$ for the Lagrange multiplier $\lambda_{\cE}$, which enforces the energy identity:
\begin{align}
\label{def:errors2}
    e_{\lambda,\cE} &= \left( \Delta t \sum_{m=1}^{M} |\lambda_{\cE}^m|^2 \right)^{1/2}.
\end{align}
We will use the numerical solution obtained with $K_{\Gamma^h} = 160$ and $\Delta t = T/2048$ as a reference solution for the convergence tests.

\begin{figure}[h!]
\centering
\includegraphics[width=0.22\textwidth]{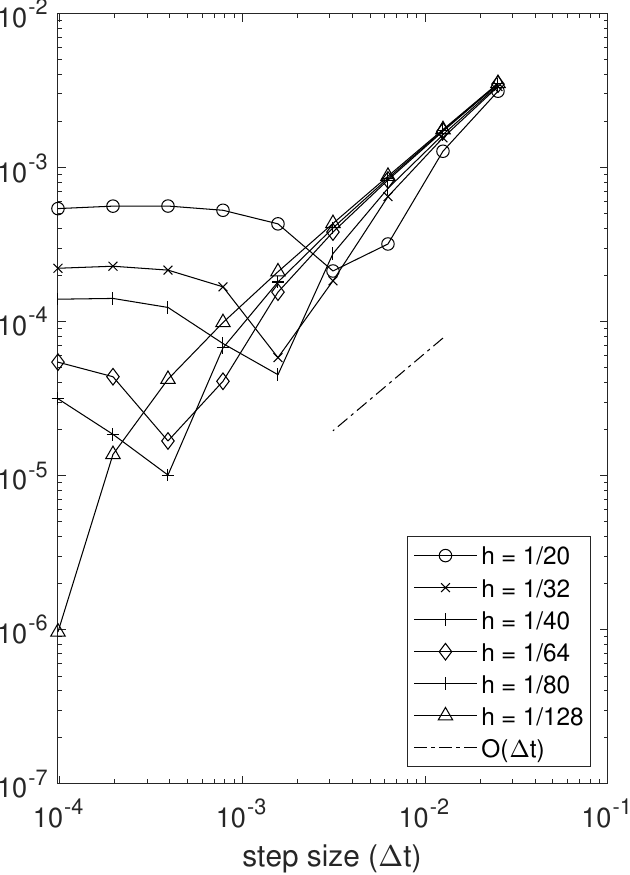}
\includegraphics[width=0.22\textwidth]{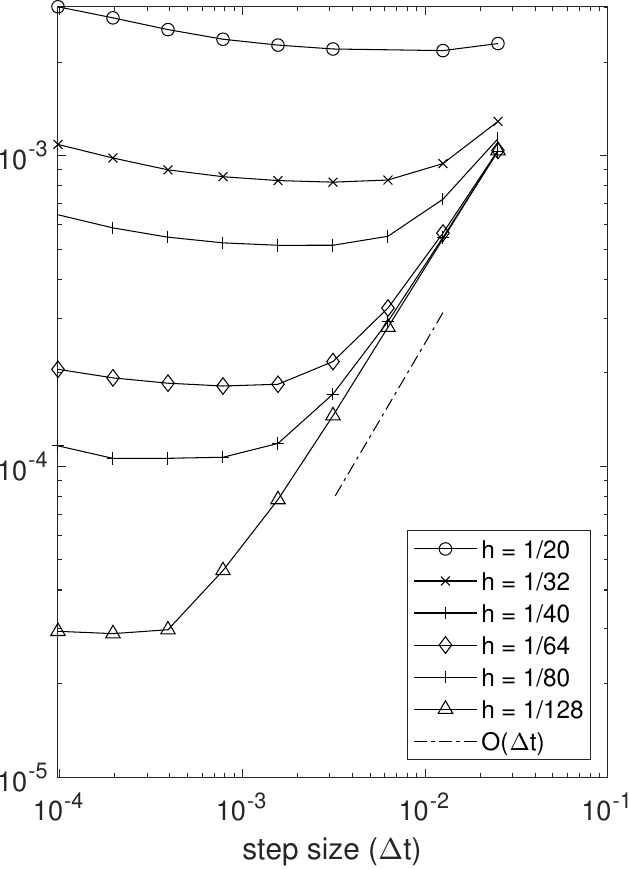}
\includegraphics[width=0.22\textwidth]{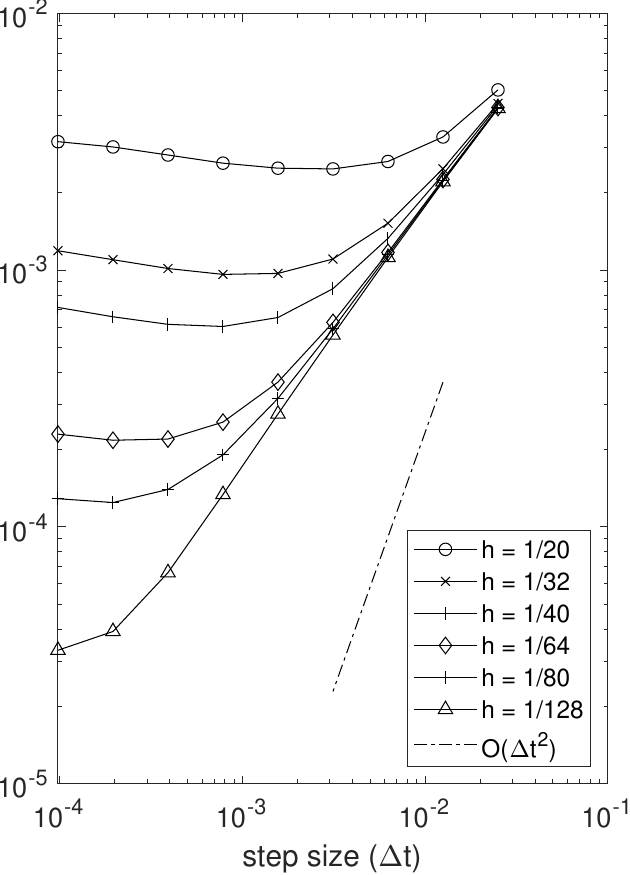}
\includegraphics[width=0.22\textwidth]{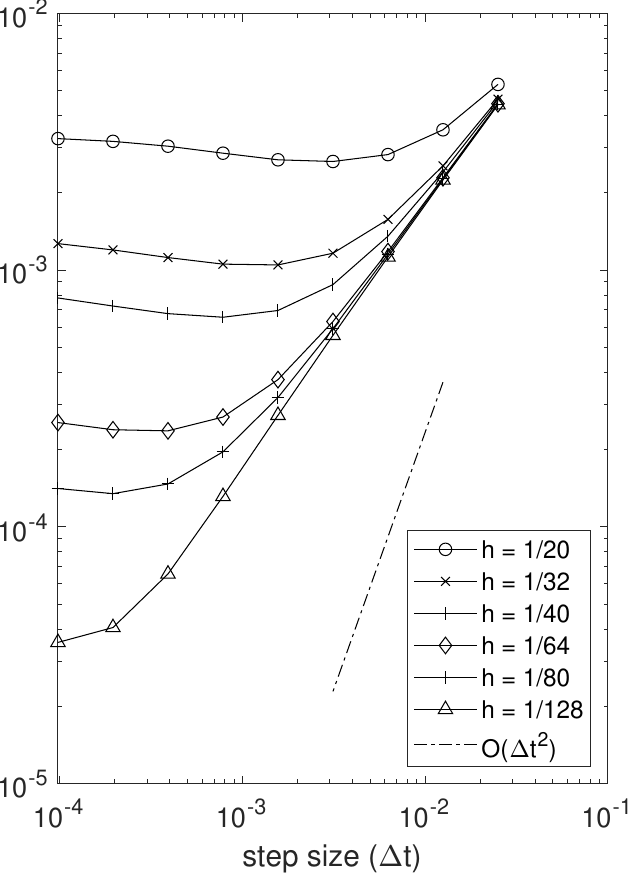}
\caption{Studying the temporal convergence of the error quantity $e_{\bX}$ for the parametrization $\bX$ in the retraction experiment. From left to right: BGN scheme, Euler scheme, BDF2 scheme, Crank--Nicolson scheme.}
\label{fig:4_retracting}      
\end{figure}

\begin{figure}[h!]
\centering
\includegraphics[width=0.22\textwidth]{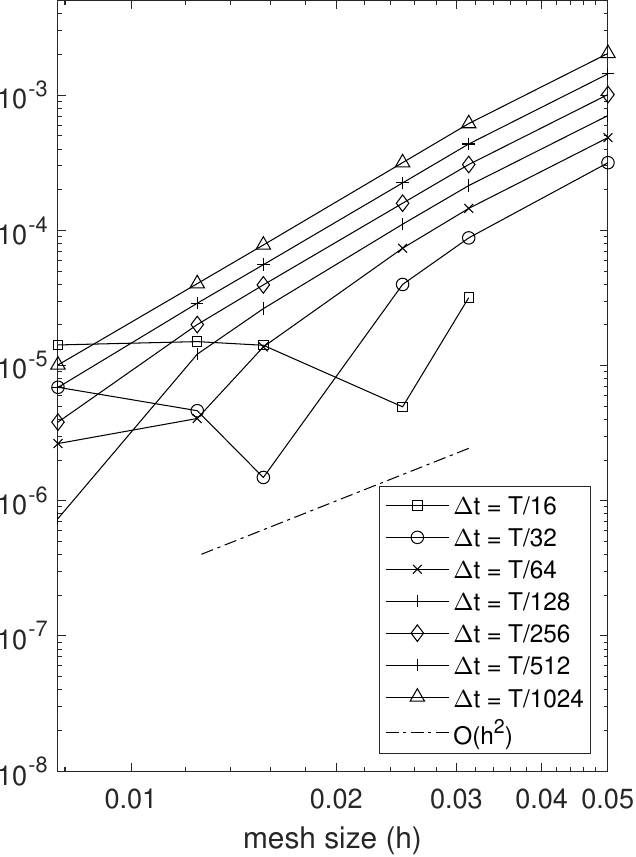}
\includegraphics[width=0.22\textwidth]{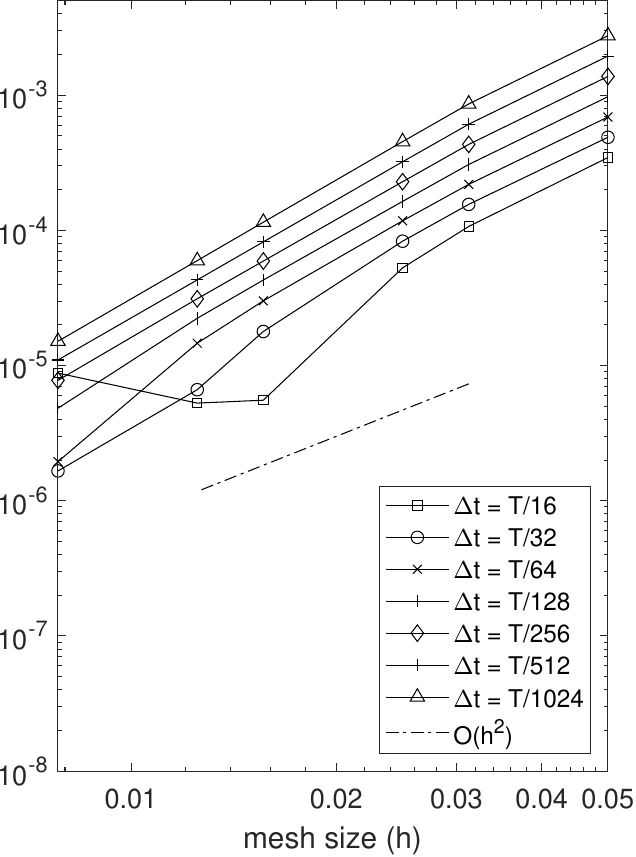}
\includegraphics[width=0.22\textwidth]{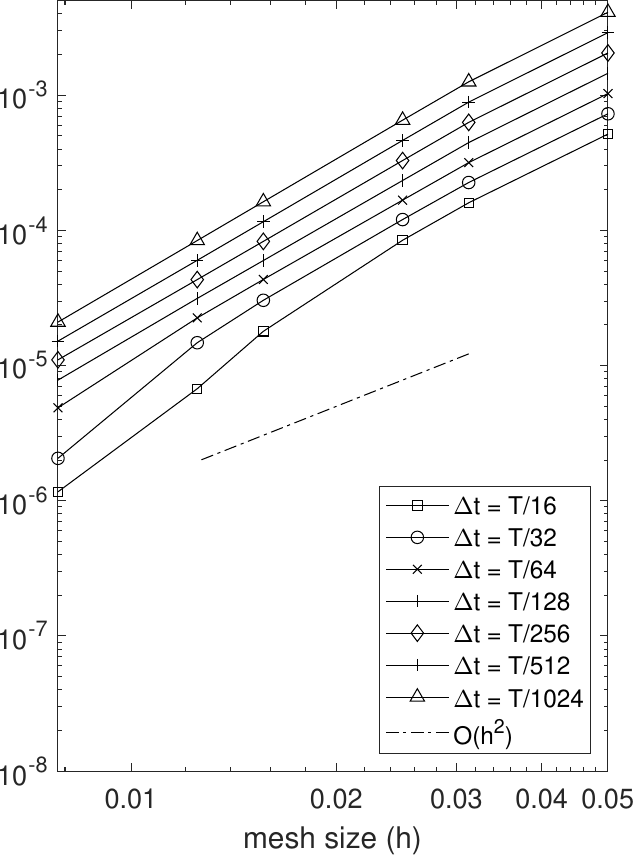}
\caption{Studying the spatial convergence of the error quantity $e_{\lambda,V}$ for the Lagrange multiplier $\lambda_{V}$ in the retraction experiment. From left to right: Euler scheme, BDF2 scheme, Crank--Nicolson scheme.}
\label{fig:5_retracting}      
\end{figure}

\begin{figure}[h!]
\centering
\includegraphics[width=0.22\textwidth]{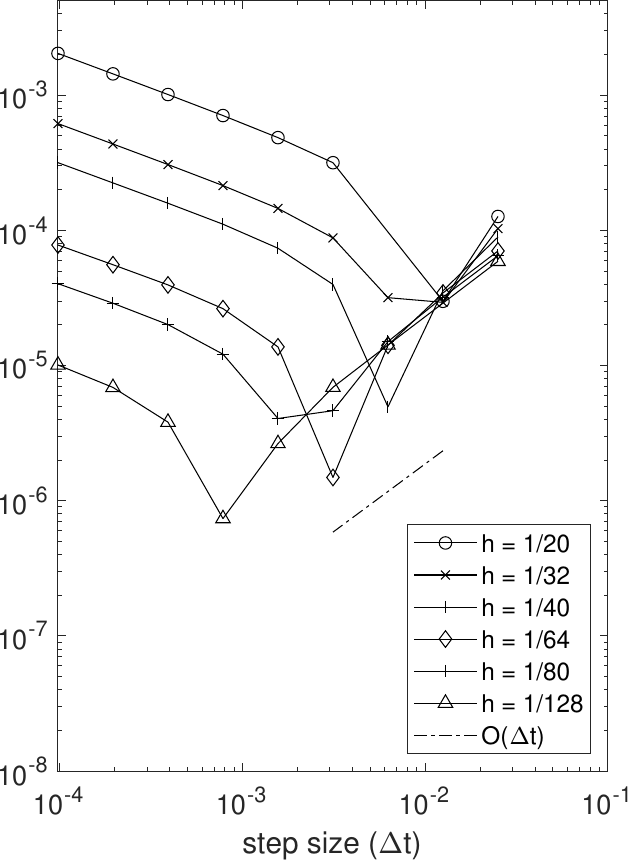}
\includegraphics[width=0.22\textwidth]{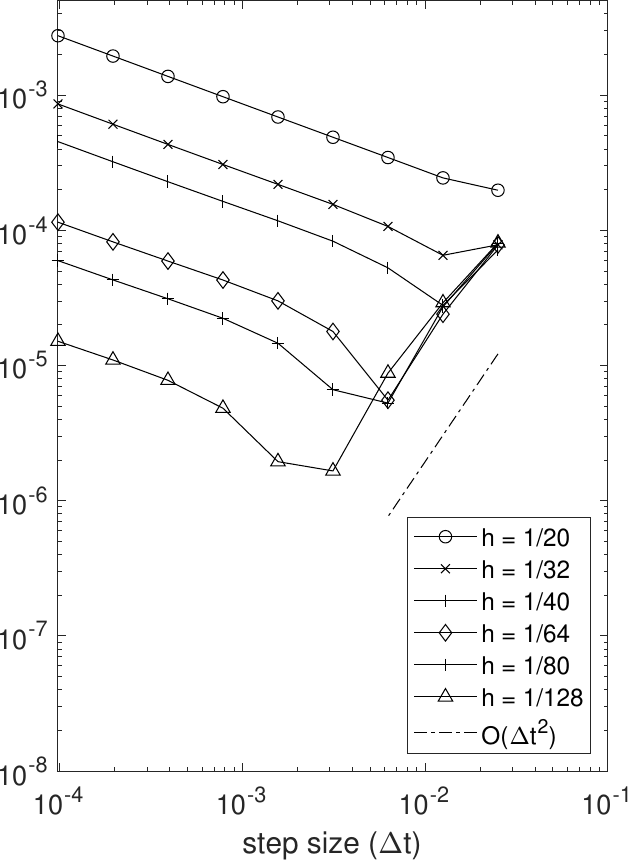}
\includegraphics[width=0.22\textwidth]{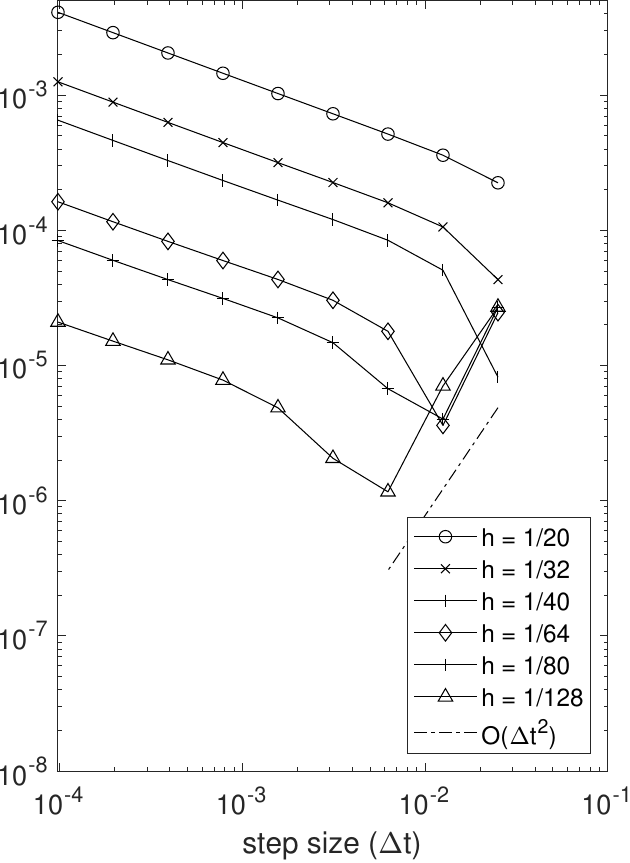}
\caption{Studying the temporal convergence of the error quantity $e_{\lambda,V}$ for the Lagrange multiplier $\lambda_{V}$ in the retraction experiment. From left to right: Euler scheme, BDF2 scheme, Crank--Nicolson scheme. 
}
\label{fig:6_retracting}      
\end{figure}

\begin{figure}[h!]
\centering
\includegraphics[width=0.22\textwidth]{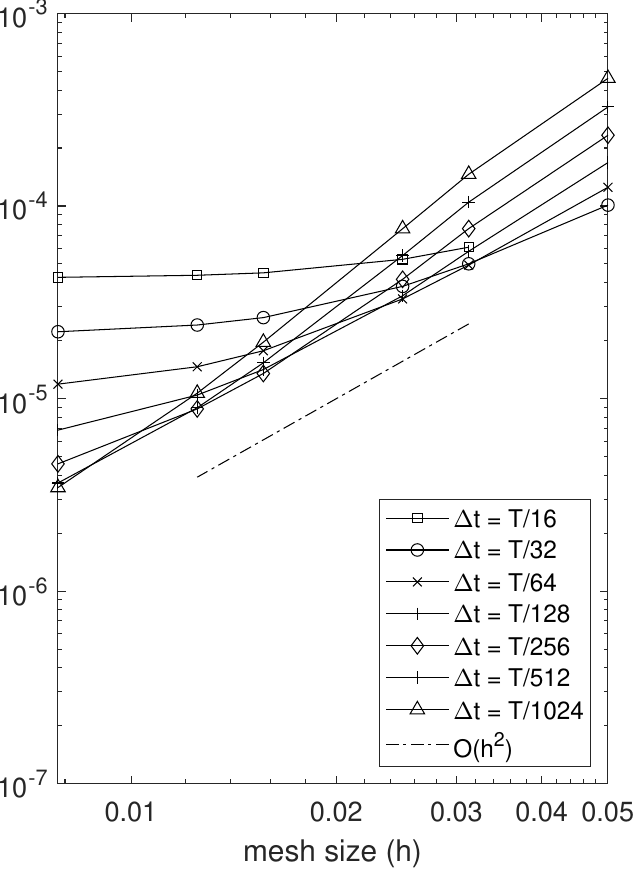}
\includegraphics[width=0.22\textwidth]{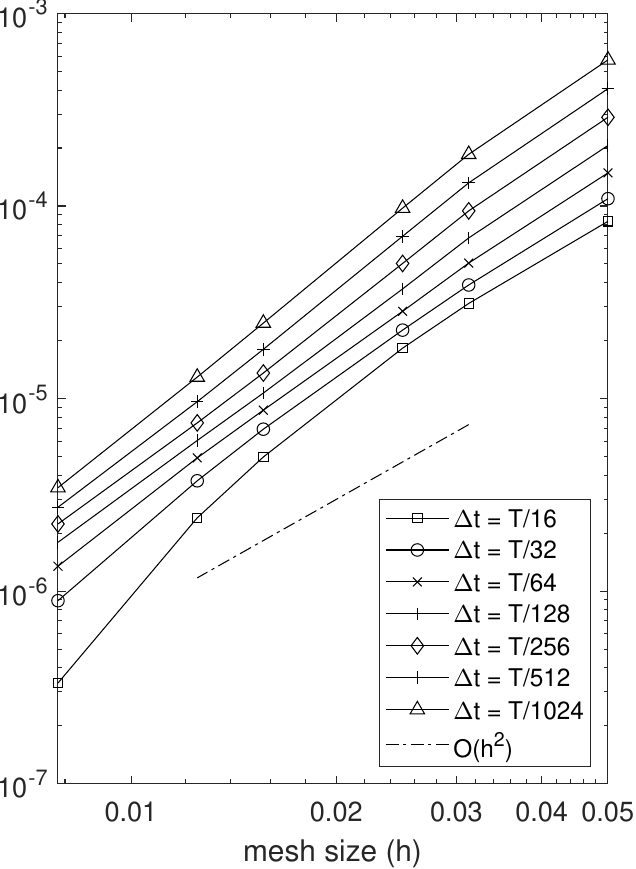}
\includegraphics[width=0.22\textwidth]{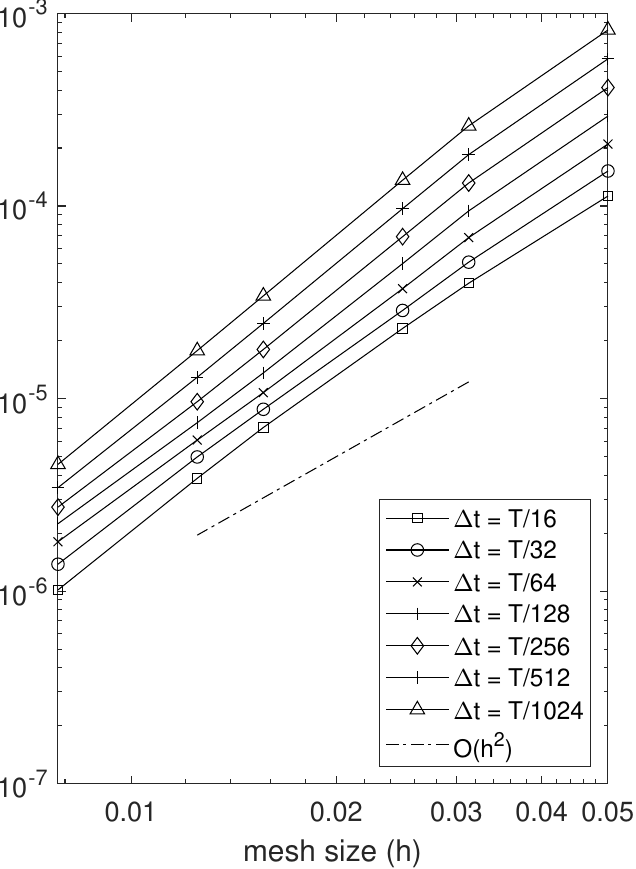}
\caption{Studying the spatial convergence of the error quantity $e_{\lambda,\cE}$ for the Lagrange multiplier $\lambda_{\cE}$ in the retraction experiment. From left to right: Euler scheme, BDF2 scheme, Crank--Nicolson scheme.}
\label{fig:7_retracting}      
\end{figure}

\begin{figure}[h!]
\centering
\includegraphics[width=0.22\textwidth]{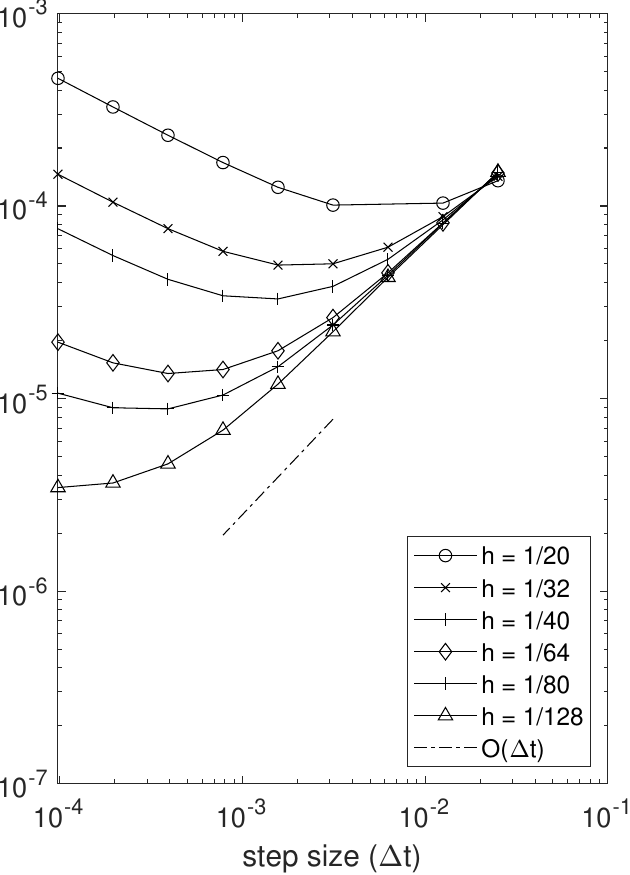}
\includegraphics[width=0.22\textwidth]{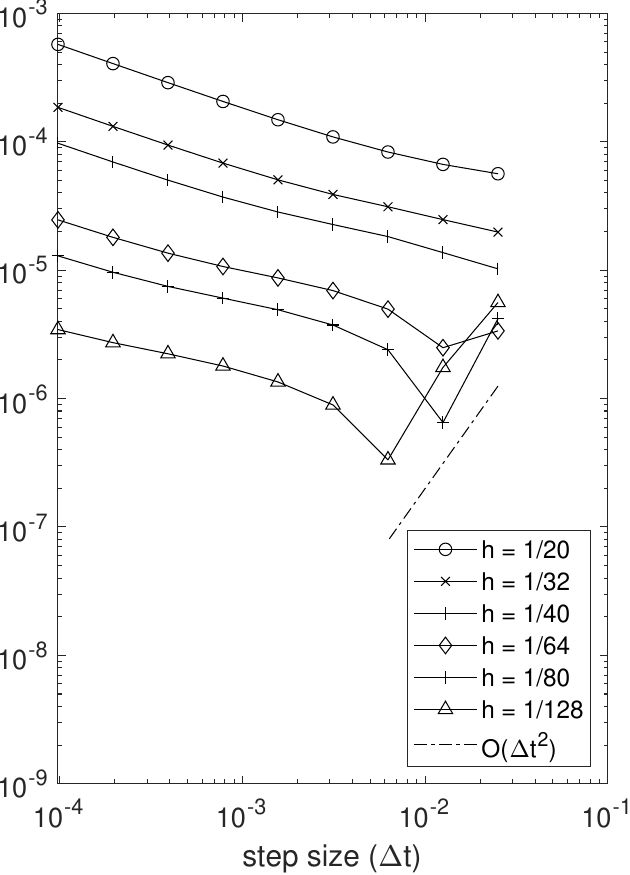}
\includegraphics[width=0.22\textwidth]{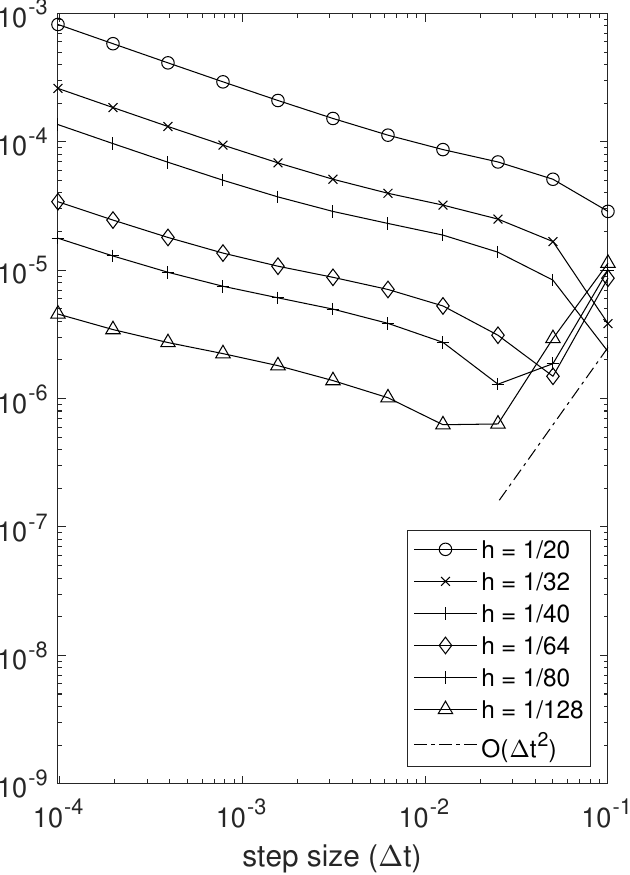}
\caption{Studying the temporal convergence of the error quantity $e_{\lambda,\cE}$ for the Lagrange multiplier $\lambda_{\cE}$ in the retraction experiment. From left to right: Euler scheme, BDF2 scheme, Crank--Nicolson scheme. 
}
\label{fig:8_retracting}      
\end{figure}

Our main focus in this experiment is on the temporal order of convergence of some selected quantities. We note that the spatial behaviour closely resembles that observed in the expanding bubble test case. 

Figure~\ref{fig:4_retracting} presents the interface error $e_{\bX}$ plotted against the time step size. The results align well with the expected convergence behaviour. However, the second-order methods do not fully achieve second-order error reduction rates, which may be attributed to the use of a numerical reference solution rather than an exact one.

For the Lagrange multipliers $\lambda_V$ and $\lambda_{\cE}$, we observe a different behaviour compared to the expanding bubble case, where the associated errors remained nearly constant under spatial refinement. In Figures~\ref{fig:5_retracting} and~\ref{fig:7_retracting}, we show the temporal convergence of the errors $e_{\lambda,V}$ and $e_{\lambda,\cE}$. For larger time steps, the error curves exhibit slopes consistent with the expected convergence orders of the respective schemes. At smaller time steps, however, the curves tend to flatten, indicating that spatial discretization errors dominate in this regime. Figures~\ref{fig:6_retracting} and~\ref{fig:8_retracting} display the spatial convergence of $e_{\lambda,V}$ and $e_{\lambda,\cE}$. Here, the observed behaviour is slightly better than second-order convergence in space.

\subsection{Viscous flow under gravitational forcing}

We extend the Stokes equation to include a gravitational force in a regime with very low Reynolds number \cite{salguero_2025_stokes_gravity}:
\begin{align*}
    -2\mu_\pm \nabla\cdot \mathbb{D}(\bu) + \nabla p = \rho_\pm \boldsymbol{f},
\end{align*}
where $\boldsymbol{f} = -0.981 \be_2$ with $ \be_2 = (0,1)^\top$, and $\rho_\pm$ denote the constant mass densities in the two phases $\Omega_\pm(t)$. As before, we define the piecewise constant density function by
\begin{align*}
    \rho(\cdot,t) =\rho_+ (1-\bm\chi_{\Omega_-(t)}) + \rho_- \bm\chi_{\Omega_-(t)}.
\end{align*}
This setting is used to perform a rising bubble experiment.

In this context, the energy and volume identities in \eqref{Physical-structures} remain valid, with the gravitational term $\bg$ replaced by $\rho \boldsymbol{f}$ in the energy equation.
We consider the numerical schemes \eqref{Lag:Euler}, \eqref{BGN}, \eqref{Lag:CN}, and \eqref{Lag:BDF2}, modified by substituting $\bg$ with $\rho \boldsymbol{f}$. The density $\rho$ is discretized analogously to \eqref{Discrete_mu}. Specifically, we define $\bg^m = \rho^m \boldsymbol{f}$ in \eqref{Lag:Euler} and \eqref{BGN}, and set $\widetilde{\bg}^{m+\frac12} = \widetilde{\rho}^{m+\frac12} \boldsymbol{f}$ in \eqref{Lag:CN} and $\widetilde{\bg}^{m+1} = \widetilde{\rho}^{m+1} \boldsymbol{f}$ in \eqref{Lag:BDF2}, where the densities $\widetilde{\rho}^{m+\frac12}$ and $\widetilde{\rho}^{m+1}$ are constructed using the predictor surfaces $\widetilde{\Gamma}^{m+\frac12}$ and $\widetilde{\Gamma}^{m+1}$, respectively.

We use the following parameters:
\begin{align*}
    \Omega = (0,1) \times (0,2), 
    \quad
    T = 3,
    \quad
    \mu_\pm = 2,
    \quad
    \rho_+ = 100,
    \quad
    \rho_- = 0.1,
    \quad
    \gamma = 0.1,
\end{align*}
with discretization parameters 
\begin{align*}
    K_{\Gamma^h} = 240,
    \quad \Delta t = 10^{-3}.
\end{align*}
The initial parametric surface $\Gamma^0$ is composed of $K_{\Gamma^h}$ elements and approximates an ellipse centered at $(0.5, 0.5)^\top$ with horizontal and vertical axes of length $0.3$ and $0.2$, respectively.
Specifically, the vertices $\bq_k^0$, $k \in \{1, \ldots, K_{\Gamma^h}\}$, of $\Gamma^0$ are defined by
\begin{align*}
    \bq_k^0 = 
    \begin{pmatrix}
        0.5 \\ 0.5
    \end{pmatrix}
    +
    \begin{pmatrix}
        0.3 \cos(\phi_k) \\ 0.2 \sin(\phi_k)
    \end{pmatrix},
    \quad 
    \phi_k = \frac{2\pi k}{K_{\Gamma^h}}, \quad k \in \{1,\ldots,K_{\Gamma^h}\}.
\end{align*}

In the unfitted approach, a uniform bulk mesh with mesh size $h = 0.02$ is used throughout the simulations. No remeshing or mesh smoothing is required, which offers a computational advantage: the bulk mesh remains fixed over time, while the interface can evolve independently of it.
In contrast, the fitted approach employs local mesh refinement near the interface. Specifically, elements within a distance of approximately $0.05$ from the interface are assigned a fine mesh size of $h \approx 0.006$, while elements farther away use a coarser mesh size of $h \approx 0.035$. This allows for higher resolution near the interface but requires the bulk mesh to be adapted at each time step, using the strategies described in \cite{agnese_nurnberg_2016_stokes}.

We perform the rising bubble experiment using fitted meshes for the BGN scheme \eqref{BGN}, the structure-preserving Euler scheme \eqref{Lag:Euler}, as well as the Crank--Nicolson scheme \eqref{Lag:CN} and the BDF2 scheme \eqref{Lag:BDF2}. Additionally, we perform computations on unfitted meshes for the BGN scheme and the Euler scheme. 

\begin{figure}
\centering
\includegraphics[width=0.23\textwidth,trim={6cm 0.5cm 6cm 0.5cm},clip]{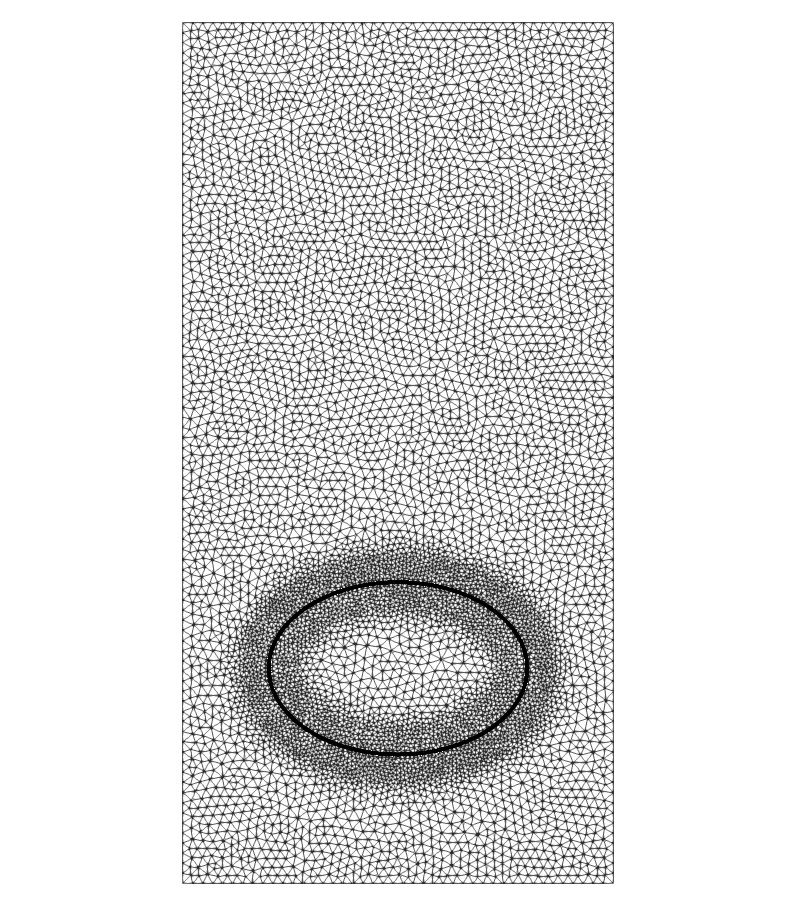} \,
\includegraphics[width=0.23\textwidth,trim={6cm 0.5cm 6cm 0.5cm},clip]{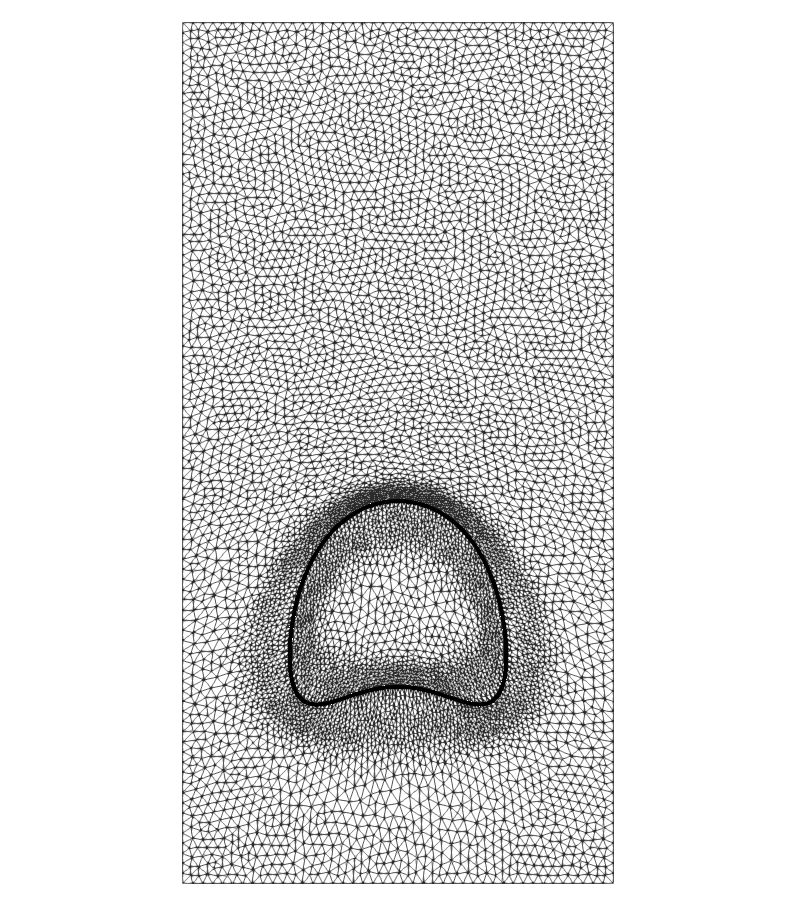} \,
\includegraphics[width=0.23\textwidth,trim={6cm 0.5cm 6cm 0.5cm},clip]{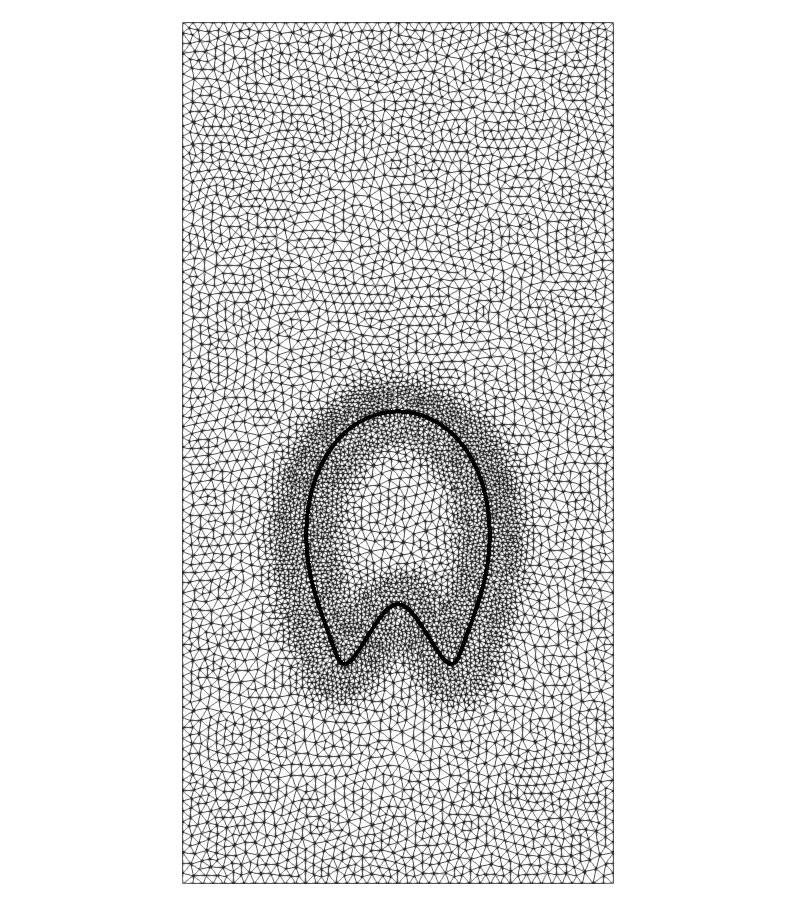} \,
\includegraphics[width=0.23\textwidth,trim={6cm 0.5cm 6cm 0.5cm},clip]{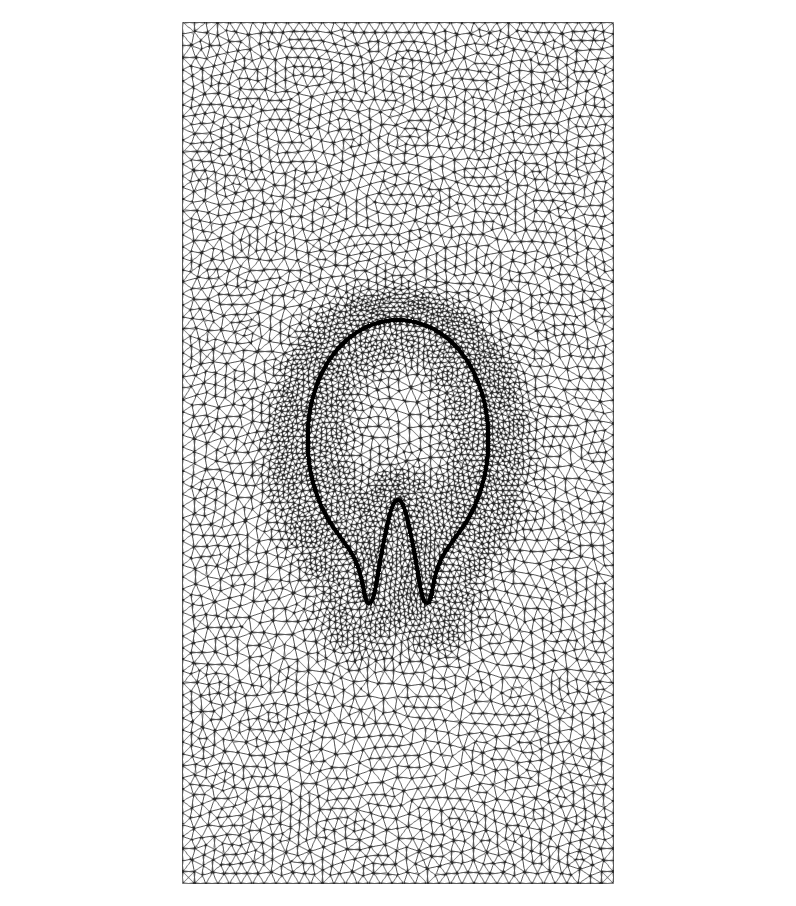}
\caption{
The interface evolution for the rising bubble at times $t \in \{0, 1, 2, 3\}$, computed using the Euler scheme with a fitted bulk mesh. Here, we use $K_{\Gamma^h} = 240$, $\Delta t = 10^{-3}$, and a local mesh size $h \in [0.006, 0.035]$.}
\label{fig:rising_bubbles_fitted}
\end{figure}

\begin{figure}
\centering
\includegraphics[width=0.23\textwidth,trim={6cm 0.5cm 6cm 0.5cm},clip]{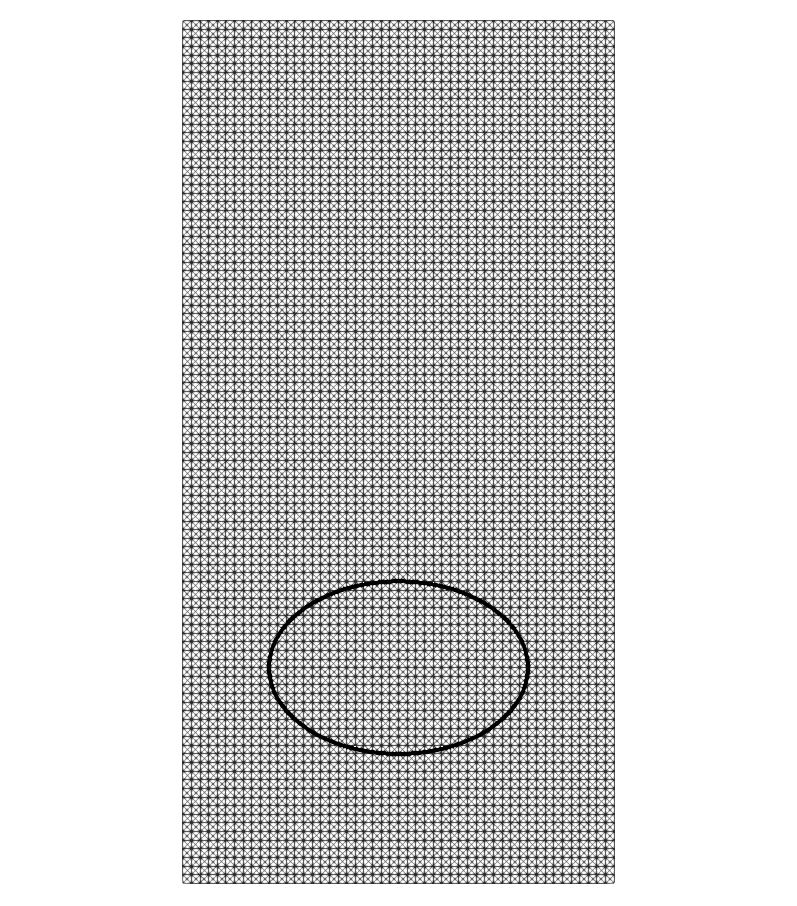} \,
\includegraphics[width=0.23\textwidth,trim={6cm 0.5cm 6cm 0.5cm},clip]{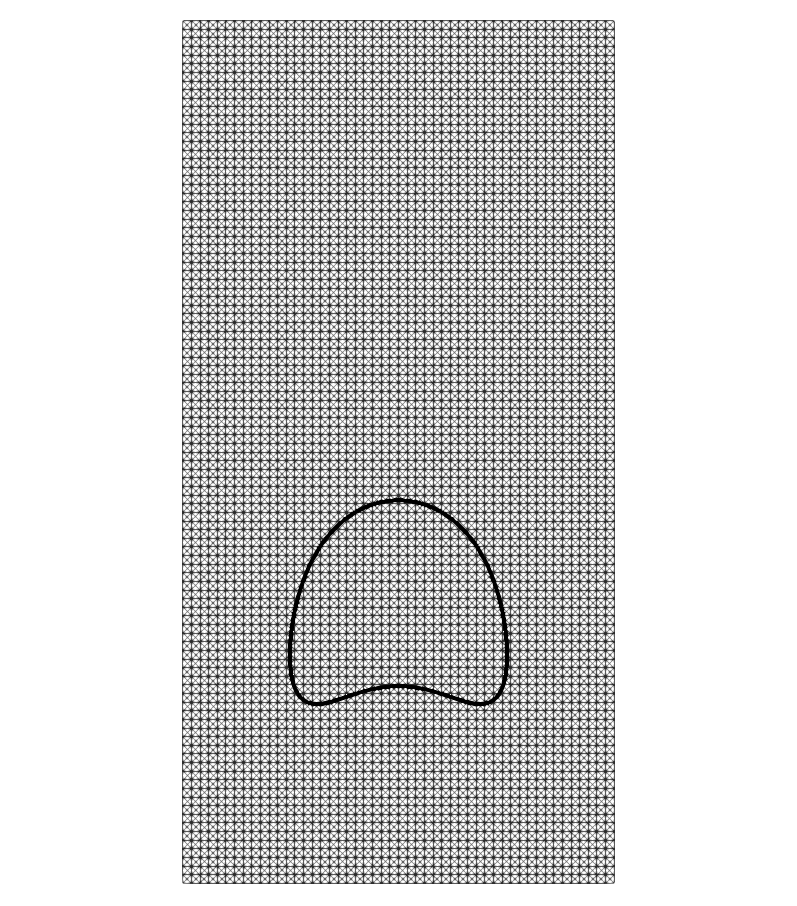} \,
\includegraphics[width=0.23\textwidth,trim={6cm 0.5cm 6cm 0.5cm},clip]{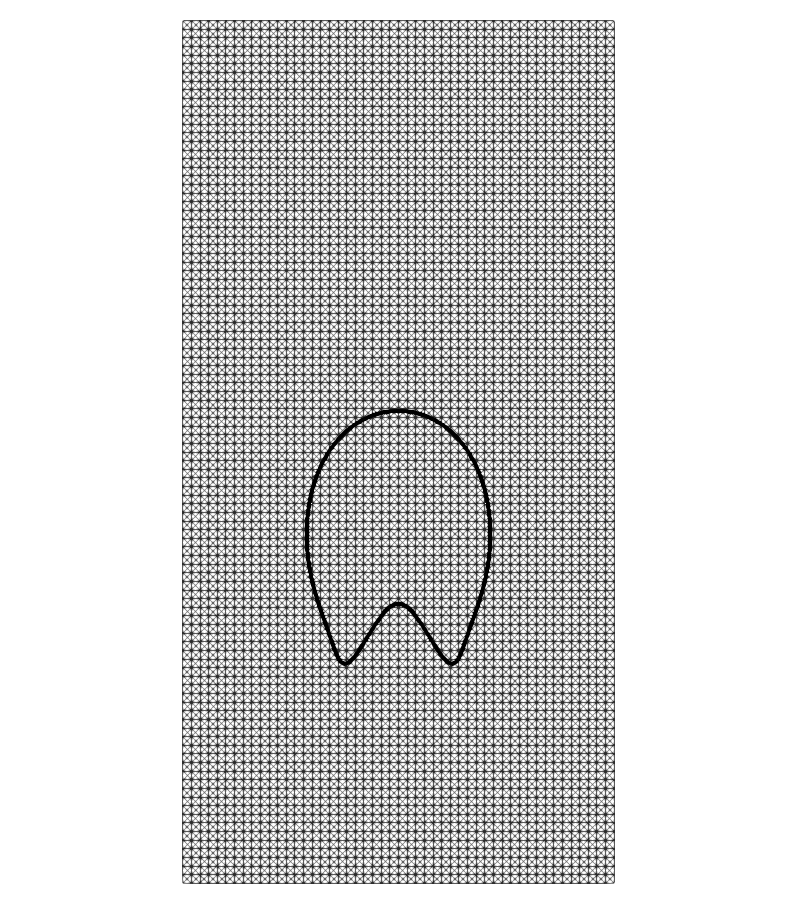} \,
\includegraphics[width=0.23\textwidth,trim={6cm 0.5cm 6cm 0.5cm},clip]{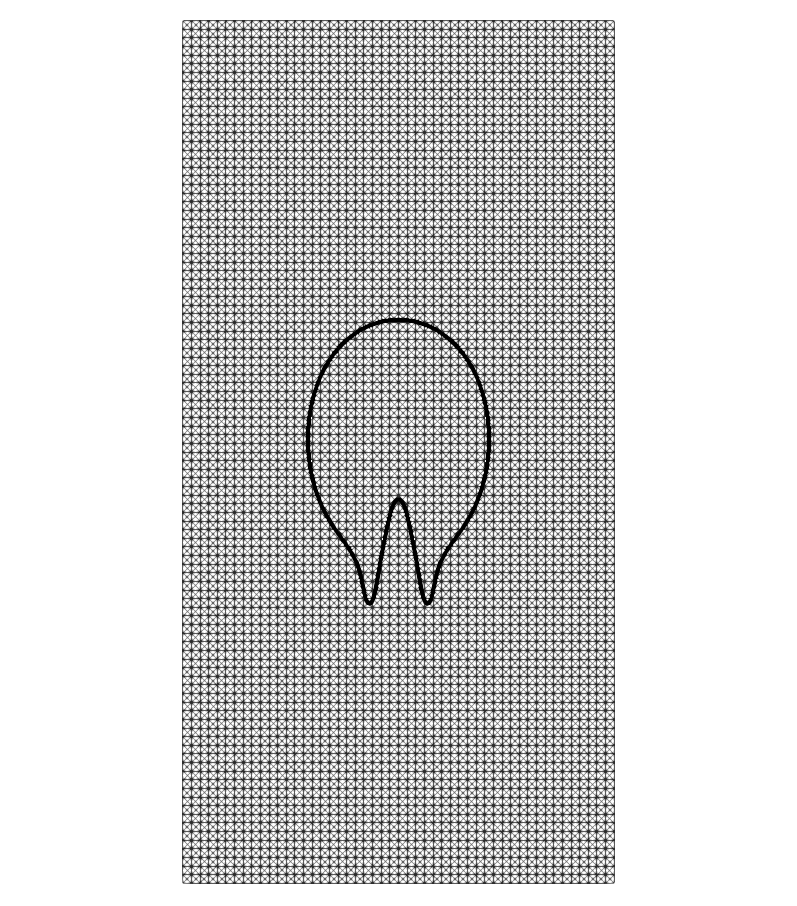}
\caption{
The interface evolution for the rising bubble at times $t \in \{0, 1, 2, 3\}$, computed using the Euler scheme with an unfitted bulk mesh. The parameters are $K_{\Gamma^h} = 240$, $\Delta t = 10^{-3}$, and a uniform mesh size $h = 0.02$.}
\label{fig:rising_bubbles_unfitted}
\end{figure}

Figures \ref{fig:rising_bubbles_fitted} and \ref{fig:rising_bubbles_unfitted} illustrate the interface evolution at times $t \in \{0, 1, 2, 3\}$ for the fitted and unfitted mesh approaches, respectively. For clarity, only the meshes for the structure-preserving Euler schemes are shown, as the interface dynamics for the other methods are visually indistinguishable in this setting.

\begin{figure}
\centering
\includegraphics[width=0.24\textwidth]{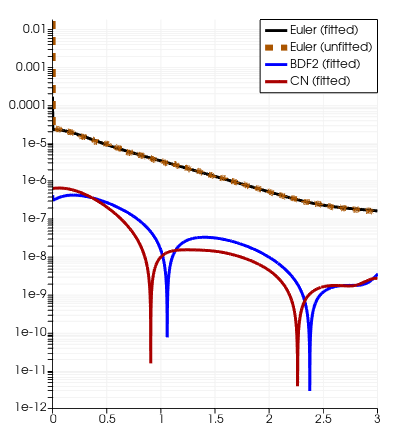} 
\includegraphics[width=0.24\textwidth]{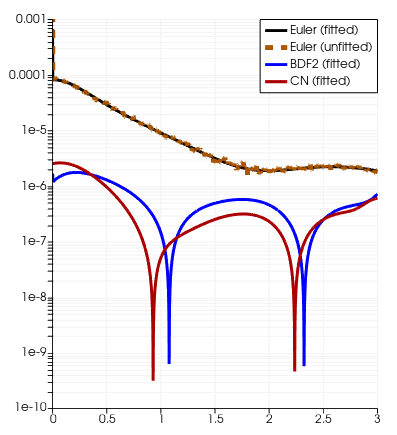} 
\includegraphics[width=0.24\textwidth]{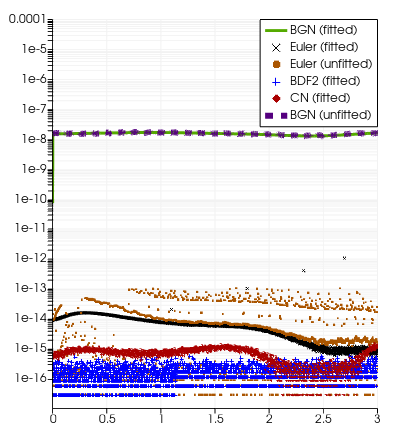} 
\includegraphics[width=0.24\textwidth]{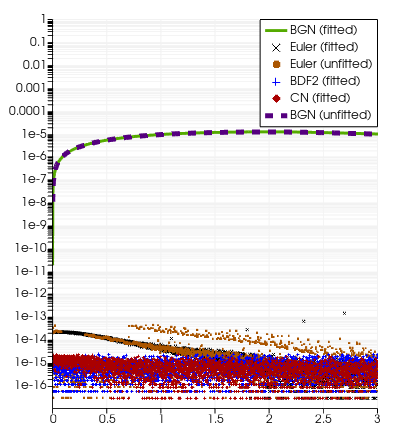}
\caption{Numerical results for the rising bubble experiment over the time interval $[0, 3]$, computed using the BGN scheme and the structure-preserving Euler, BDF2, and Crank--Nicolson (CN) schemes. For the fitted bulk meshes, local mesh sizes $h \in [0.006, 0.03]$ are used with $K_{\Gamma^h} = 240$ and $\Delta t = 10^{-3}$. Additionally, the BGN and structure-preserving Euler schemes are applied on unfitted bulk meshes with uniform mesh size $h = 0.02$, also using $K_{\Gamma^h} = 240$ and $\Delta t = 10^{-3}$.
From left to right: time evolution of the magnitudes of the Lagrange multipliers $\lambda_{\mathcal{E}}$ and $\lambda_V$, followed by the numerical residuals of the energy and volume identities for each respective scheme.}
\label{fig:rising_Lagr}
\end{figure}

In Figure \ref{fig:rising_Lagr}, we plot the time evolution of the Lagrange multipliers $\lambda_{\cE}$ and $\lambda_{V}$, along with the numerical residuals of the energy and volume identities. For both the fitted and unfitted BGN schemes, the residuals are of similar magnitude, approximately $10^{-8}$ for the energy and between $10^{-7}$ and $10^{-5}$ for the volume identity. For all structure-preserving schemes, these residuals are at the level of machine precision.

The Lagrange multipliers $\lambda_{\cE}$ obtained from the fitted and unfitted Euler schemes are nearly identical and decay from about $10^{-5}$ to $10^{-7}$ over time, and $\lambda_{V}$ shows a similar trend, decreasing from around $10^{-4}$ to $10^{-6}$. In contrast, the second-order methods yield significantly smaller quantities of Lagrange multipliers, with $\lambda_{\cE}$ ranging from $10^{-6}$ to $10^{-12}$ and $\lambda_{V}$ from $10^{-5}$ to $10^{-10}$.

\section{Conclusion}

We propose a novel approach for two-phase Stokes flow based on a Lagrange multiplier method. The resulting equations are discretized by a parametric finite element method and several different type of time discretization methods, including the Euler method, the Crank--Nicolson method and the second-order backward differentiation formula (BDF2). The employment of Lagrange multipliers enables the fully discrete schemes to both decrease the energy and preserve the volume. In practice, the nonlinear schemes are solved by the Newton 
method with an efficient decoupling technique. Extensive numerical experiments are presented to demonstrate the desired temporal accuracy and structure-preserving properties. Similar Lagrange multiplier approaches can be extended to simulate other complex interface problems involving, e.g., the Navier--Stokes system, multiphase problems, a coupling to equations modeling surfactants and biomembranes \cite{BGN_2015_navierstokes, BGN2015c, BGN2015b, BGNbio17, Garcke-Nurnberg-Zhao2023}, which will be part of our future research.


\bibliographystyle{abbrv}
\bibliography{ref}

\end{document}